\documentclass[11pt]{amsart}
\usepackage[T1]{fontenc}
\usepackage[utf8]{inputenc}
\usepackage{lmodern}
\usepackage[a4paper,margin=27mm]{geometry}
\usepackage{amsmath,amssymb,amsthm,mathtools}
\usepackage{booktabs,longtable,array}
\usepackage{xcolor,graphicx,tikz}
\usepackage[colorlinks=true,linkcolor=blue!40!black,citecolor=blue!40!black,urlcolor=blue!40!black,breaklinks=true]{hyperref}
\usepackage{microtype}

\numberwithin{equation}{section}
\numberwithin{figure}{section}
\newcommand{\F}{\mathbb F_2}
\newcommand{\CP}{\mathbb C P}
\newcommand{\HP}{\mathbb H P}
\newcommand{\OP}{\mathbb O P}

\DeclareMathOperator{\ev}{ev}

\newtheorem{theorem}{Theorem}[section]
\newtheorem{proposition}[theorem]{Proposition}
\newtheorem{lemma}[theorem]{Lemma}
\newtheorem{corollary}[theorem]{Corollary}
\theoremstyle{definition}

\theoremstyle{remark}
\newtheorem{remark}[theorem]{Remark}
\makeatletter
\@namedef{subjclassname@2020}{\textup{2020} Mathematics Subject Classification}
\makeatother
\title[Closed geodesics on compact rank-one symmetric spaces]{Infinitely many closed geodesic images on simply connected compact rank-one symmetric spaces}
\author{Zhiqi Chen}
\address[Zhiqi Chen]{School of Mathematics and Statistics, Guangdong University of Technology, Guangzhou 510520, P. R. China}
\email{chenzhiqi@nankai.edu.cn}
\author{Shaoqiang Deng}
\address[Shaoqiang Deng]{School of Mathematical Sciences and LPMC, Nankai University, Tianjin 300071, P. R. China}
\email{dengsq@nankai.edu.cn}
\author{Hui Zhang}
\address[Hui Zhang]{School of Mathematics, Southeast University, Nanjing 211189, P. R. China}
\email{huizhang@mail.nankai.edu.cn}
\date{15 October 2026\quad Independent-review revision for author review}
\hypersetup{pdftitle={Infinitely many closed geodesic images on simply connected compact rank-one symmetric spaces},pdfauthor={Zhiqi Chen, Shaoqiang Deng, Hui Zhang},pdfsubject={Technical revision with an explicit controlled finite-diagram cut construction}}
\begin{document}
\raggedbottom
\begin{abstract}
We prove that every smooth Riemannian metric on the underlying smooth
manifold of a simply connected compact rank-one symmetric space, namely
$S^n$ $(n\ge2)$, $\CP^m$ $(m\ge1)$, $\HP^m$ $(m\ge1)$, or
$\OP^2$ admits infinitely many geometrically distinct prime closed
geodesics. Here a prime closed geodesic is a nonconstant closed geodesic
which is not an iterate of a shorter one, and geometric distinctness
means that the images in the manifold are pairwise different. Thus
iterating a geodesic, changing its starting point, or reversing its
orientation does not count as a new geodesic. No nondegeneracy assumption
on the metric is imposed.
These manifolds have bounded loop-space Betti numbers, so the usual
homological growth criterion does not yield this conclusion. We use
the length filtration together with rotation and reversal of loops.
A controlled cut construction gives a uniform filling estimate in the
based loop space, with a fixed increase in length independent of degree,
which we transfer to free loops. We then calculate global equivariant
classes associated with odd and even iterates. Assuming finitely many
prime geodesic images forces incompatible requirements at a class's
first critical level. The proof treats spheres and projective spaces
within a common framework and does not assume a sphere infinitude theorem.
\end{abstract}
\maketitle
\setcounter{tocdepth}{1}
\tableofcontents

\section{Introduction}\label{sec:introduction}

\subsection{The problem and the limitation of ordinary homology}

The purpose of this paper is to study geometric multiplicity of closed
geodesics on spheres and simply connected projective spaces with an
arbitrary Riemannian metric. These manifolds are natural test cases:
their loop spaces have homology in arbitrarily high degrees, but the
ranks in individual degrees remain bounded. The question is therefore
how to use this topology to distinguish a new geodesic from an iterate
of one already present.

A closed geodesic on a Riemannian manifold is a nonconstant periodic solution
of the geodesic equation. Its iterates are again closed geodesics, so an
existence theorem immediately produces infinitely many parametrized critical
points of the energy. The substantive multiplicity question asks for
infinitely many \emph{geometrically distinct} closed geodesics: their images in
the manifold must be different. In particular, changing the initial point,
reversing orientation, or traversing the same curve several times does not
produce another image. A closed geodesic is \emph{prime} if it is not an
iterate of a shorter closed geodesic.

The interaction between variational methods and periodic dynamics is already
central in Birkhoff's work~\cite{Birkhoff1917}. The existence of at least one
closed geodesic on every closed Riemannian manifold of positive dimension was
proved by Lyusternik and Fet~\cite{LF1951}. Passing from existence to geometric
infinitude requires information that survives iteration. Bott's iteration
theory~\cite{Bott1956} describes the Morse indices of the iterates of one
closed geodesic. Gromoll and Meyer~\cite{GM1969} combined iteration with local
Morse theory to obtain a topological criterion: if a closed simply connected
manifold has free-loop Betti numbers unbounded over a field, then every
Riemannian metric on it has
infinitely many geometrically distinct closed geodesics. The point of this
criterion is that finitely many prime geodesics and all their iterates cannot
account for unbounded homology in individual degrees.

Vigu\'e-Poirrier and Sullivan~\cite{VPS1976} made this criterion effective for
simply connected manifolds. If the rational cohomology algebra of a closed
simply connected manifold requires at least two generators, its free loop
space has unbounded rational Betti numbers. This settles the infinitude
problem for a large topological class. It also isolates an important
obstruction to the method: a truncated polynomial algebra has only one
generator. Spheres and the complex, quaternionic, and Cayley projective
spaces belong to this exceptional class. Their loop spaces have nonzero
homology in arbitrarily large degrees, but their rational Betti numbers
remain bounded. Infinite total homology alone therefore does not distinguish
new geodesic images from iterates of old ones.

Two further developments illustrate the different information needed in
these cases. On the two-sphere, Bangert~\cite{Bangert1993}, together with
Franks~\cite{Franks1992} or Hingston~\cite{Hingston1993}, proved geometric
infinitude for every Riemannian metric. These arguments use features of
surface dynamics and of the Birkhoff annulus construction. De Philippis,
Marini, Mazzucchelli, and Suhr~\cite[Theorem~1.5]{DPMMS2022} extended the
infinitude conclusion to reversible Finsler two-spheres. In another
direction, Rademacher~\cite[Theorems~1.1--1.3]{Rademacher1994} proved
infinitude for strongly bumpy metrics on compact simply connected manifolds
and showed that these metrics form a residual set. Such generic results do
not give the assertion for an arbitrary metric: infinitely many geodesics
for approximating metrics may coalesce or acquire increasingly large
periods in a limit. Degenerate metrics must therefore be included in the
argument itself.

\subsection{The main statement}

This paper concerns the underlying smooth manifolds of the simply connected
compact symmetric spaces of rank one. Thus the manifold is one of
\[
 S^n\ (n\ge2),\qquad \CP^m\ (m\ge1),\qquad
 \HP^m\ (m\ge1),\qquad \OP^2.
\]
Here $\CP^m$ and $\HP^m$ have real dimensions $2m$ and $4m$,
respectively, and $\OP^2$ is the $16$-dimensional Cayley projective plane.
We use the abbreviation CROSS for a compact rank-one symmetric space.
The standard smooth structure is understood throughout.

\begin{theorem}\label{thm:main}
Let $M$ be diffeomorphic to one of
\[
 S^n\ (n\ge2),\qquad \CP^m\ (m\ge1),\qquad
 \HP^m\ (m\ge1),\qquad \OP^2.
\]
Let $g$ be a smooth Riemannian metric on $M$. Then $g$ has infinitely many prime closed geodesics with pairwise distinct
images.
\end{theorem}

The assertion concerns every smooth Riemannian metric on each of these
manifolds. In particular, the symmetric-space structure is used to
understand the topology, while the metric in the statement may have
arbitrary curvature and degenerate closed geodesics. The identities
$\CP^1=S^2$ and $\HP^1=S^4$ are included. A prime geodesic need not be
embedded; the conclusion distinguishes the images of the curves.

For the symmetric metric, closed geodesics occur in large families. Those
families disappear under a general perturbation, but the topology of the
loop space remains. Our aim is to use more of that topology than the ranks
of its ordinary homology groups. Two additional structures are relevant:
the filtration by length, and the action on loops by rotation and reversal.
The first records the amount of length needed to represent or fill a cycle.
The second distinguishes odd and even iterates. The argument combines a
uniform filling estimate with this distinction between iterates.

\subsection{Relation to the sphere manuscript and scope of the contribution}

OpenAI's September 2026 preprint~\cite[Theorem~1.1]{OpenAI2026}
states the same infinitude conclusion for every smooth Riemannian metric
on the standard sphere $S^n$, $n\ge2$, including degenerate metrics.
We cite that work explicitly as an antecedent for the sphere case and
as a source of the finite-cut and doubling strategy. Accordingly, the
sphere conclusion itself is not presented here as a new claim of
priority. Relative to that sphere theorem, the additional manifolds
in Theorem~\ref{thm:main} are $\CP^m$ and $\HP^m$ for $m\ge2$,
and $\OP^2$; the cases $\CP^1$ and $\HP^1$ are already spheres.
The title refers to the underlying smooth manifolds: the metric under
consideration is arbitrary and need not be symmetric.

The independence asserted here is independence of the proof from the
results of~\cite{OpenAI2026}. Neither its sphere theorem nor its
intermediate lemmas are assumed as mathematical inputs. Instead, we
develop the controlled based-loop filling argument, transfer it to free
loops, construct the required global equivariant classes for both spheres
and projective spaces, and prove the finite-prime obstruction in this
paper. The projective cases require their own global multiplication and
restriction calculations, including the complex projective extension
and the Cayley-plane characteristic-class calculation in
Section~\ref{sec:projective-topology}. Thus the proof is intended to be
checkable using the arguments given here and the published or accepted
inputs identified below. This logical independence does not assert an
independent discovery of the finite-cut and doubling strategy; its
source remains~\cite{OpenAI2026}.

\subsection{Loop products, critical levels, and related work}

String topology provides operations that retain information beyond Betti
numbers. The Chas--Sullivan homology product admits the published homotopy-theoretic
construction of Cohen--Jones~\cite{CJ2002} and the geometric chain
construction of Hingston--Wahl~\cite{HW2023}. Concatenation makes it
compatible with addition of lengths.
Goresky and Hingston~\cite{GH2009} constructed a relative cohomology product
with a complementary filtration inequality and related these products to
closed-geodesic critical levels. Hingston and Rademacher~\cite{HR2013}
used the two products to prove resonance for loop homology of spheres of
dimension greater than two:
critical levels in high degree have a common asymptotic ratio to degree.
This is precisely the kind of information that ordinary ranks omit, but
resonance alone is compatible with a finite collection of iterated
geodesics whose average indices have the corresponding ratios to length.

Recent work further clarifies the role of filtrations. Ginzburg, G\"urel,
and Mazzucchelli~\cite{GGM2024} study barcode entropy for the energy
functional, relating the persistence of homological classes to the dynamics
of the geodesic flow. Nabutovsky and Rotman~\cite{NR2013} give quantitative
homotopy constructions on loop spaces. Block, Manin, and
Weinberger~\cite{BMW2026} place controlled fillings in the broader setting
of persistent homology of function spaces. In the present argument the
relevant question is a uniform one: can a globally null cycle be filled
after adding a length independent of its degree?

We approach the uniform filling question through the based loop space,
whose nonzero homology groups over $\mathbb F_2$ are one-dimensional.
This simple algebra has a useful geometric consequence. A homology
power measures the length needed to restore a degree, while its
cohomological test measures the length removed by a cut. Comparing
these two quantities allows a filling problem in high degree to be
reduced to one of finitely many low degrees. To return the resulting
filling to the original cycle, we construct a homotopy parametrized by
a square. Both the identity and its additive length bound are needed;
they are proved together in Lemma~\ref{lem:based-cut-square}.

The cap constructions are related to the geometric chain methods
of \cite{HW2023}, and the boundary identity has an unfiltered
counterpart in the open-string algebra of \cite{CO2026}. The role of
our construction is to retain the required support estimate throughout
the homotopy. After transferring the filling estimate to free loops,
we combine it with equivariant cohomology. Here a second issue arises:
local Morse data must give actual global classes with specified
multiplication and restriction properties. We establish these
properties before applying the finite-prime criterion of
Section~\ref{sec:finite-prime}. Thus the analytic estimates, the global
calculation, and the local comparisons have separate roles in the
proof. The relation to~\cite{OpenAI2026} is specified above, and the
acknowledgements record the additional sources and computational assistance.

\subsection{Loop spaces, filtrations, and symmetries}\label{subsec:conventions}

We introduce the notation needed to describe the proof. Based loops
are convenient for concatenation and cutting; free loops are the
variational space whose critical points are closed geodesics.

Except where another coefficient ring is explicitly stated, homology and
cohomology have coefficients in $\F=\mathbb F_2$. Write $S^1=\mathbb R/\mathbb Z$.
For a fixed smooth Riemannian metric $g$ on the closed manifold $M$, set
\[
 \Lambda M=H^1(S^1,M),\qquad
 \Omega_qM=\{\gamma\in\Lambda M:\gamma(0)=q\},\qquad n=\dim M.
\]
These spaces carry the strong Sobolev $H^1$ topology. The inclusion
$M\hookrightarrow\Lambda M$ identifies a point with its constant loop.
The evaluation map is $\ev_0(\gamma)=\gamma(0)$, and
$SM=\{(x,v)\in TM:g_x(v,v)=1\}$ denotes the unit tangent bundle.
For $\gamma\in\Lambda M$ define
\begin{equation}\label{eq:intro-energy}
 \ell(\gamma)=\int_0^1|\dot\gamma(t)|_g\,dt,
 \qquad E(\gamma)=\int_0^1|\dot\gamma(t)|_g^2\,dt,
 \qquad F(\gamma)=\sqrt{E(\gamma)}.
\end{equation}
Thus $\ell\le F$, with equality for constant-speed loops. If $X$ is a
subspace of $\Lambda M$, the notations $X^{F<a}$ and $X^{\ell<a}$ mean the
corresponding open sublevels. A superscript $X^a$ always refers to the
filtration explicitly fixed in that section. Constants in uniform
estimates may depend on $M$ and $g$, but not on homological degree or
filtration level, unless a different dependence is stated.

A closed curve can be parametrized from a different initial point or
with the opposite orientation. These changes give the action of
$G=O(2)$ on the free loop space. We shall use reflection and
half-period translation separately, since a loop fixed by the latter
is the double of a shorter loop. We write
\[
 (\rho\gamma)(t)=\gamma(-t),\qquad
 (\tau\gamma)(t)=\gamma(t+\tfrac12),\qquad
 R=\langle\rho\rangle,\quad K=\langle\rho,\tau\rangle.
\]
Then $R\cong\mathbb Z/2$ and $K\cong(\mathbb Z/2)^2$.
To retain these symmetries in cohomology we use the Borel construction.
It replaces a group action by a space over its classifying space and
provides natural restriction maps when the acting group is reduced.
For a compact group $H$ and an $H$-invariant pair $(X,A)$, set
\[
 H_H^*(X,A)=H^*(EH\times_HX,EH\times_HA;\F),
\]
where $EH$ is a contractible free $H$-space and $BH=EH/H$.
The coefficient classes are
\[
 H^*(BG;\F)=\F[w,U],\quad |w|=1,\quad |U|=2,
 \qquad H^*(BR;\F)=\F[w].
\]
Here $w$ and $U$ are the first and second Stiefel--Whitney classes of the
standard real plane representation; its restriction to $R$ has $U=0$.
The notation $w$ also denotes its pullback to an equivariant pair.
Evaluation at zero is $R$-equivariant for the trivial action on $M$.
Consequently, a specified class $\eta\in H^d(M;\F)$ gives an evaluation
class $h=\ev_0^*\eta$ acting on reflection-equivariant cohomology. The
choice of $\eta$ and its degree will be stated in each case. Evaluation
at a fixed parameter is not asserted to be $G$-equivariant.

\subsection{The mechanism of the proof}

We describe the three steps before entering their constructions. The first
is an estimate for cycles, the second is a computation of equivariant
cohomology, and the third compares that computation with the consequences
of having only finitely many prime geodesics.

\emph{Uniform filling and spectral concentration.}
For $X=\Omega_qM$ let $A=\varnothing$, and for $X=\Lambda M$ let $A=M$.
Write $X^a=X^{F<a}$, where $a>0$. Using Proposition~\ref{hyp:controlled-cut-realization}, the analytic
part gives numbers $D,C\ge0$ and $\lambda>0$ such that, in every degree
$j\ge0$,
\begin{equation}\label{eq:intro-death}
 \ker\{H_j(X^a,A)\longrightarrow H_j(X,A)\}
 \longrightarrow H_j(X^{a+D},A)
 \quad\text{is the zero map},
\end{equation}
and the map $H_j(X^a,A)\to H_j(X,A)$ is zero when
$a<\lambda j-C$ and surjective when $a>\lambda j+C$.
We call these properties uniform death and spectral resonance,
respectively. Uniform death bounds the lifetime of a class which
eventually vanishes; spectral resonance locates classes which survive in
the full space. The estimates are needed together.

On based loops, Pontryagin multiplication restores length while a
cohomologically tested cut removes it. The one-dimensional homology
groups allow the two costs to be compared up to a fixed additive error.
The square of two cut parameters supplies a chain homotopy that returns
the cut-and-filled cycle to the original one. Crucially, the two output
loops together traverse each input only once. This yields an additive
estimate, rather than an estimate with a multiplicative loss in the
original length. A finite-cover Mayer--Vietoris argument then transfers
both properties to free loops. Resonance controls the otherwise
uncontrolled preimage in the Mayer--Vietoris sequence.

\emph{Odd and even equivariant classes.}
The symmetric metric is used only to calculate the unfiltered
equivariant cohomology of $(\Lambda M,M)$. The calculation separates
classes arising from odd symmetric iterates from those arising from even
ones. The odd classes occur in pairs related by multiplication by the
evaluation class $h$. One summand in each pair has lifts to
$G$-equivariant cohomology annihilated by a fixed power of $U$; its
$h$-multiple supplies the paired summand. The even classes have a different
lifting behavior. These are statements about global classes and maps,
not merely about the local groups of a Morse filtration. In complex
projective space, an extension at the next critical index requires a
separate argument; for the Cayley plane, the required characteristic
class is computed from its actual geodesic-space bundle.

\emph{Finite cuts and doubling.}
Suppose there are only finitely many prime images. Index iteration and
local Morse theory then place the local groups near finitely many lines
in the degree--length plane, with uniform bounds on their dimensions.
The analytic estimates restrict how classes at a finite sublevel can
disappear in the full loop space. Combining this control with restriction
to the fixed locus of $\tau$ gives a uniform bound for the total number
of transient reflection-equivariant classes. The two Gysin sequences for
$G$, $R$, and $K$ compare a regular cut with its double. Equality in the
resulting rank inequality forces multiplication by $h$ to be injective on
the visible part of the odd summand having full-group lifts. At the first
appearance of an appropriate odd
class, however, the evaluation class $h$, whose degree is at least two,
vanishes on the local reflection model. This gives
the contradiction. The local argument includes degenerate critical orbits
and geodesics with zero mean index, using the ordinary normal-stabilization
comparison proved in Proposition~\ref{prop:local-thom-realization}.

Section~\ref{sec:based-analytic} proves the based-loop estimates and compares
length with the strong-$H^1$ energy filtration.
Section~\ref{sec:finite-cover} passes to free loops.
Sections~\ref{sec:sphere-topology} and~\ref{sec:projective-topology}
establish the sphere and projective equivariant calculations.
Section~\ref{sec:finite-prime} proves the obstruction to finitely many
prime images, and Section~\ref{sec:main-proof} assembles these results.

\paragraph{Sources and scope.}
All external theorems invoked as inputs to the proofs are taken from
published papers or, in the case of \cite{BMW2026}, a manuscript
explicitly accepted for publication. The preprint~\cite{OpenAI2026}
is cited for the antecedent sphere result and the source of the strategy;
no result from it is used as an unproved premise. Author-hosted and arXiv copies are used to locate the
statements in those papers. The based Pontryagin algebra, the controlled
cut square, and the finite-prime criterion are proved below; the local
comparison is reduced to published index-pair results. The distinction
between these mathematical inputs and the historical sources of the
strategy is explained in the acknowledgements.

\section{Length estimates on the based loop space}
\label{sec:based-analytic}

The main result of this section is
Proposition~\ref{prop:based-uniform-persistence}. It controls two kinds
of homology classes in a length sublevel. A class which eventually
vanishes must vanish after a bounded increase in length. A class which
survives has its first occurrence within a bounded distance of a line
in the degree--length plane. In both statements the bound is
independent of degree.

The difficulty is the first assertion. Known filling estimates in a
fixed degree allow a constant which grows with that degree. Our method
is to cut a high-degree cycle down to a bounded degree, fill the cut
cycle, and restore the removed degree by concatenation. There are two
questions to answer. Does this procedure fill the original cycle, and
does restoring the degree cost more length than the cut saved? The
cut square answers the first question. The comparison of Pontryagin
powers with Goresky--Hingston powers answers the second. We develop
these ingredients before giving the filling construction.

All chains and cochains in this section have coefficients in $\mathbb F _2$.
Let $M$ be a closed simply connected Riemannian manifold of dimension $n\ge2$,
let $q\in M$, and write $\Omega=\Omega_qM$ for the Hilbert manifold of $H^1$
based loops. We use the functions
\[
 \ell(\gamma)=\int_0^1|\dot\gamma(t)|\,dt,
 \qquad F(\gamma)=\left(\int_0^1|\dot\gamma(t)|^2\,dt\right)^{1/2}.
\]
The energy is $E=F^2$. Until the subsection on length and Hilbert-energy
filtrations, superscripts denote open length sublevels.
For a nonzero homology class $a$ define its length minimax by
\[
 c(a)=\inf\{A:a\in\operatorname{im}
 [H_*(\Omega^{<A})\longrightarrow H_*(\Omega)]\}.
\]
For a positive-degree cohomology class $\alpha$ write $t(\alpha)$ for the
supremum of the levels $L$ for which $\alpha$ has a cocycle representative
vanishing on $C_*(\Omega^{<L})$.
Thus $c(a)$ measures the first level at which $a$ is represented, whereas
$t(\alpha)$ measures the last level on which $\alpha$ can be made to
vanish. Homology classes are represented by finite singular cycles; cochains are
ordinary singular cochains and are not required to have finite support.
All chain diagrams used in the constructions are finite. A claim
at an infimum or supremum is used with an arbitrarily small positive
slack; no attainment of a minimax value is assumed.

\subsection{What is required of the cut construction}

A loop which returns near its basepoint can be cut into two based
loops by inserting short joining arcs. We test one piece with a
cohomology class and retain the other. If the test vanishes on short
loops, a nonzero contribution forces the discarded piece to be long;
this is the source of the length saving. The following proposition
states the two properties needed later: a one-input chain map and a
compatible two-input homotopy.

The compatibility is required on a finite collection of chains,
including their faces and the fillings used in the argument. We call
this collection a finite diagram. Making this choice before subdividing
ensures that repeated inputs have exactly the same realization.

\begin{proposition}[Controlled finite-diagram cut realization]
\label{hyp:controlled-cut-realization}
For the given metric there is $K\geq0$ with the following property.
Let $T$ be a normalized cocycle of degree $p>0$ which vanishes on
$C_*(\Omega^{<L})$, with $L>0$, and put $Q=n+p-1$.
On every prescribed finite diagram of chains, their faces, and the
fillings and parameter homotopies used in the argument, the geometric
cuts described below admit compatible realizations satisfying:
\begin{enumerate}
\item $D_T^L,D_T^R$ are degree $-Q$ chain maps on the diagram and
send a chain supported below $a$ to a chain supported below $a-L+K$.
For positive-degree cocycles $S$ and cycles $z$ their homological
pairing is $\langle S,D_T^Lz\rangle=\langle T\circledast S,z\rangle$,
where $\circledast$ is the usual relative based GH product.
\item For cycles $x,y$ of degrees $j,N>p$, supported below $a,b$,
there is a square chain of degree $j+N-Q+1$, including its refinement
and timing corrections, with boundary
\[
 D_T^R(y)*x+y*D_T^L(x)
\]
and support below $a+b-L+K$.
\end{enumerate}
The constant is independent of $p,j,N,L,a,b$ and of the finite
diagram. Finite positive slacks may be absorbed by increasing it once.
Every square requested in a given finite diagram must use the same
chosen cuts on repeated inputs. A single realization is chosen after
all inputs of that diagram have been specified. No comparison with a
previously chosen realization on a smaller diagram, and no single
subdivision valid for all singular chains, is required.
\end{proposition}
\begin{proof}
Part~(1) is proved by the supported-cap construction in
Proposition~\ref{prop:one-input-cuts}. Part~(2), including its
simultaneous choice of those same cuts, is proved independently in
Lemma~\ref{lem:based-cut-square}. Both constructions use a fixed
joining-arc bound $\rho$ and length oscillation less than $1$.
Thus one constant, for example $K=2\rho+2$, works for the entire
finite diagram, in all degrees.
\end{proof}

\subsection{A cut operation on finite chain families}

We now prove the proposition. The construction has three stages.
First, a Thom cap restricts attention to loops returning near the
basepoint, where cutting and closing the two pieces are defined.
Second, a cohomological test selects the piece to discard. Finally,
a square of cut parameters relates the left and right versions.
The details of subdivision serve these stages: they make the caps
well defined and preserve the boundary identities at the same time.

We use normalized singular chains, the standard cross product, and
Alexander--Whitney cup and cap products. Normalization matters here:
the vanishing of two boundary faces of the parameter square will use
actual degenerate simplices. To make subdivision precise,
write $C_*^{\mathrm{un}}(X)$ for ordinary singular chains and
$C_*(X)=N_*(X)$ for their quotient by degenerate singular simplices.
Write
$\pi:C_*^{\mathrm{un}}(X)\to C_*(X)$ for normalization and choose its
natural Moore chain section $\iota$. Equivalently, $\iota$ selects
the representative in the Moore subcomplex
$\bigcap_{i=1}^k\ker d_i\subset C_k^{\mathrm{un}}(X)$ in degree $k$,
where $d_i$ is the $i$th face map. Lift each finite input diagram by $\iota$,
subdivide its ordinary chains factorwise, perform the ordinary shuffle
and cap construction, and apply $\pi$ to every output. Normalized
cochains are pulled back by $\pi$. The maps $\iota,\pi$, ordinary
subdivision $\mathrm{Sd}$, and its prism $P$ preserve the image of each
parent simplex. For $k\geq1$ put
\[
 P_k=\sum_{r=0}^{k-1}\mathrm{Sd}^{r}P,\qquad
 \mathcal S_k=\pi\mathrm{Sd}^{k}\iota,\qquad
 \mathcal P_k=\pi P_k\iota.
\]
Telescoping gives
$\partial P_k+P_k\partial=\mathrm{Sd}^{k}-\mathrm{id}$ and hence
$\partial\mathcal P_k+\mathcal P_k\partial=\mathcal S_k-\mathrm{id}$.
Every summand stays in the image of its parent simplex. Thus $\mathcal P_k$ is the homotopy which restores an input after all
$k$ subdivisions. The section $\iota$ is used only as a natural
chain map preserving the image of each simplex; it need not commute
with subdivision or be multiplicative. Products are formed on the
ordinary lifts, and normalization is applied after the cap operations. We will use this cumulative correction at the end
of the square construction. Performing subdivision before the final
normalization is necessary because ordinary subdivision need not
preserve the singular degeneracy subcomplex.
Fix a normalized diagonal Thom
cocycle $U$ of degree $n$ supported in an arbitrarily small tubular
neighborhood of the diagonal of $M\times M$. Its restriction to each
normal fibre represents the generator of
$H^n(D^n,S^{n-1};\mathbb F_2)$. Write
\[
 U_L(z)=U(z,q),\qquad U_R(z)=U(q,z).
\]
These are normalized point Thom cocycles. For $z$ sufficiently close to $q$,
$[z,q]$ denotes the unique short geodesic from $z$ to $q$.
We shall also use a representative for which endpoint vanishing follows
directly from support. Inside a fixed outer
diagonal tube, move its zero section to the graph of a sufficiently small
diffeomorphism $f$ with $f(q)\ne q$, using a normal-bundle isotopy which
is fixed near the tube boundary. The transported Thom cocycle represents
the same class relative to the complement of the outer tube and may be
supported near that graph, away from $(q,q)$. Its two point restrictions
are supported near $f^{-1}(q)$ and $f(q)$, respectively, still within the
small return ball about $q$. This is a representative of the outer-tube
Thom class, not a cochain supported on the exact original diagonal.
All endpoint terms then vanish by support. Comparing with the undisplaced
representative uses a normalized primitive of degree $n-1$; its pullback
to the constant endpoint is zero since $n>1$. Thus this stronger choice
does not change the induced absolute operation.
For a returning loop and a cut parameter define
\[
 L(\gamma,t)=\gamma|_{[0,t]}\#[\gamma(t),q],\qquad
 R(\gamma,t)=[q,\gamma(t)]\#\gamma|_{[t,1]}.
\]
Restriction to a variable interval, followed by linear rescaling, is
$H^1$-continuous, including collapsed intervals: the energy of the rescaled
prefix is $t\int_0^t|\dot\gamma|^2$, which tends to zero as $t\to0$.
In each of $L$ and $R$, give the two displayed pieces equal time
slots of length $1/2$. The binary product $*$ likewise gives each
input one half of the interval. These allocations are independent
of input lengths, including when a piece is constant. Other slot
allocations, such as the three equal slots in the raw square, will
be stated explicitly. Changing these fixed slots gives length-preserving
timing homotopies. This choice ensures that, when an output omits one input
loop, its parametrization also omits that input. An allocation of time
proportional to input lengths would lose this independence and would
not support the degeneracy argument below. The comparison with energy
is made only after the chain identities have been established.

Let $T$ be a normalized cocycle of degree $p>0$ vanishing below length $L>0$,
and put $Q=n+p-1$. Define $D_T^L$ by crossing with the cut interval, capping
first with the pullback of $U_L$ and then with $T\circ L$, and retaining
the $R$ factor. Define $D_T^R$ using $U_R$, testing $R$ and retaining $L$.
Both lower degree by $Q$.
The degree shift consists of the one parameter supplied by the cut
interval, followed by caps of degrees $n$ and $p$:
$j+1-n-p=j-Q$. The superscript on $D_T^L$ or $D_T^R$ specifies the
piece tested by $T$; the other piece is the output.

Here and below a cut is defined on any prescribed finite collection of
chains, faces, and parameterized homotopies, enlarged when necessary to
include global fillings used only to identify homology classes. Choose
an inner support tube and a strictly larger return tube. After all
parameter factors have been specified, refine each input-chain domain
and each interval factor separately, and then use the standard product
shuffle triangulation. A common sufficiently fine factorwise mesh exists
by uniform continuity on the finite compact family; using common
refinements on input faces and interval endpoints preserves all boundary
identities. No arbitrary final barycentric subdivision of the product
simplices is permitted: that could destroy the one-factor-at-a-time
steps used in the normalized degeneracy argument below. Every parent simplex
which can contribute to the Thom cap is then contained in the larger
return tube. Define a supported cup product there and extend it by zero
only after multiplication by the Thom cocycle. This is the ordinary
small-chain model for relative cup and cap products; on faces it obeys
the cap boundary identity. Lemma~\ref{lem:supported-cap-chain} below proves the local cap boundary
identity in this common small-chain complex. The same convention can
be imposed on finitely many previously specified parameter homotopies.
The local cap identity alone does not identify the square edges
with chosen products; the required shuffle and timing comparisons
are proved in Lemmas~\ref{lem:shuffle-omitted-factor}
and~\ref{lem:based-cut-square}. No universal subdivision of all singular chains is asserted.

At $t=0,1$, the pullback of the point Thom cochain factors through the
constant map to $q$. It is zero because it has positive degree and is
normalized. Thus endpoint terms present no obstruction to an absolute chain
map. Proposition~\ref{prop:one-input-cuts} gives these one-input
chain maps independently of the square construction. Indeed, the cap boundary formula has only the input
boundary, the two interval endpoints, and the coboundaries of the
cochains as possible terms. The latter vanish because $U$ and $T$ are
cocycles, and the endpoint terms have just been shown to vanish.
Their positive-output scalar pairings agree with the relative based
Goresky--Hingston coproduct \cite[Proposition~9.4]{GH2009}.
We denote its dual product by
$\circledast$; it raises total cohomological degree by $n-1$:
\[
 H^a(\Omega,\{q\})\otimes H^b(\Omega,\{q\})
 \xrightarrow{\ \circledast\ }H^{a+b+n-1}(\Omega,\{q\}),
 \qquad a,b>0.
\]
Here $\{q\}$ denotes the constant based loop. Since $\Omega$ is
connected, positive-degree relative and absolute cohomology are
canonically identified. If $S$ is a positive-degree cocycle, then
\begin{equation}\label{eq:based-1}
\langle S,D_T^L z\rangle=\langle T\circledast S,z\rangle.
\end{equation}
This statement only compares the relative coproduct on two positive
factors; no choice of an absolute lift from another model is imported.

The two output lengths sum to at most the input length plus twice the
length of a short joining arc. On a contributing parent simplex a
nonzero $T$-face contains a tested loop of length at least $L$.
Continuity of length and sufficient subdivision make every tested loop
on that parent have length at least $L-\varepsilon$. Consequently, with
a fixed additive allowance $K$ independent of $p$ and chain degree,
\begin{equation}\label{eq:based-2}
D_T^{L/R}\bigl(C_j(\Omega^{<a})\bigr)
       \subset C_{j-Q}(\Omega^{<a-L+K}).
\end{equation}
The tube radii and subdivisions may depend on the finite family and its
length bounds. Their dependence does not increase $K$: each output uses
only the fixed number of short arcs displayed above, and the support of
a sum of chains is a union.

There is one limit to the pairing argument which will matter for
fillings. It determines outputs tested by positive-degree classes,
where relative and absolute cohomology agree. A degree-zero output
requires a separate argument. In
Subsection~\ref{subsec:recover-unit} we recover that output from the
square identity and a nonzero Pontryagin power. This explains why the
square is needed even after the one-input operation has been identified.

\subsubsection*{Making the cap a chain map}
The cutting maps exist only near a return. We must therefore justify
that the cap actually produces a chain in this region before applying
them. The next lemma supplies precisely this local-to-global step.

\begin{lemma}[Supported caps on a common small-chain complex]
\label{lem:supported-cap-chain}
Let $e:P\to M\times M$ be continuous, let $V_0$ be closed and
contained in an open return tube $V_1$, and let
$U\in C^n(M\times M,(M\times M)\setminus V_0)$ be a normalized
cocycle. Put $P_1=e^{-1}(V_1)$. Suppose $A,B:P_1\to\Omega$ are
continuous, and $T\in C^p(\Omega)$ is a normalized cocycle.
Let $C_*^{\mathrm{sm}}(P)$ be the ordinary singular subcomplex
generated by simplices contained either in $P_1$ or in
$e^{-1}((M\times M)\setminus V_0)$. There is a degree $-n-p$
chain map
\[
 \mathcal K_{U,T}:C_*^{\mathrm{sm}}(P)\longrightarrow
                  C_{*-n-p}^{\mathrm{un}}(\Omega),\qquad
 \mathcal K_{U,T}(c)=B_*\bigl(A^*T\cap(e^*U\cap c)\bigr).
\]
In particular, this formula does not require extending $A$, $B$,
or $A^*T$ outside $P_1$.
\end{lemma}
\begin{proof}
First we must justify the indicated codomain of the inner cap.
For a singular $k$-simplex $\sigma$, a nonzero term of
$e^*U\cap\sigma$ has coefficient
$U(e\sigma|[0,\ldots,n])$. If $\sigma$ is contained in
$e^{-1}((M\times M)\setminus V_0)$, this coefficient is zero.
Every surviving simplex therefore has its entire parent simplex
in $P_1$, and its back face is a singular simplex of $P_1$.
Consequently $e^*U\cap c$ is an actual chain in $P_1$.

The ordinary cap boundary formula, with mod-two coefficients, gives
\[
 \partial(e^*U\cap c)=e^*U\cap\partial c.
\]
This equality holds in $C_*^{\mathrm{un}}(P_1)$: both sides lie
there by the preceding observation, and its inclusion into
$C_*^{\mathrm{un}}(P)$ is injective. Capping there with the closed
cochain $A^*T$ and applying the chain map $B_*$ gives
$\partial\mathcal K_{U,T}(c)=\mathcal K_{U,T}(\partial c)$.
Thus cancellation across the boundary of the return tube is checked
in a common chain complex, not inferred from a termwise extension
by zero. Postcomposition with normalization is also a chain map.

For a finite family of product domains, the preimages of
$V_1$ and of the complement of $V_0$ form an open cover. A common
sufficiently fine factorwise subdivision followed by the shuffle
triangulation lies in this small-chain complex. Include every face
and every already specified parameter homotopy before choosing the
subdivision; the same subdivision on a repeated factor makes its
face occurrences identical. This is a finite construction and does
not assert a universal subdivision of the whole singular complex.
\end{proof}

\begin{proposition}[The one-input cut maps]
\label{prop:one-input-cuts}
For a fixed finite family of chains closed under boundary, the
one-input cuts can be chosen to be degree $-Q$ chain maps with
the support estimate~\eqref{eq:based-2}. Their positive-output
pairing is~\eqref{eq:based-1}. These assertions are proved directly, independently of the
square construction.
\end{proposition}
\begin{proof}
Choose one subdivision order $k$ for the entire finite family and
its product with the cut interval, and set
$c_z=\mathrm{Sd}^k\iota z\times\mathrm{Sd}^k[I]$.
Use Lemma~\ref{lem:supported-cap-chain} with
$P=\Omega\times I$, $e(\gamma,t)=(\gamma(t),q)$,
$A=L$, $B=R$ to define $D_T^Lz=\pi\mathcal K_{U,T}(c_z)$.
For $D_T^R$ use $e(\gamma,t)=(q,\gamma(t))$, $A=R$, $B=L$.
Since the same $k$ is used on faces, $\partial c_z$ is the sum
of $c_{\partial z}$ and its two parameter-end chains. On each
end chain $e$ is constant at $(q,q)$. A normalized cochain of
positive degree pulls back to zero under a constant map. Both
end contributions therefore vanish, giving
\[
 \partial D_T^{L/R}z=D_T^{L/R}\partial z.
\]
The degree is $1-n-p=-Q$.

On the return tube the two output lengths have sum at most the
input length plus $2\rho$, where $\rho$ bounds the short joining
arc. Choose the common mesh also to make the tested length oscillate
by less than $1$ on each parent product cell. A nonzero evaluation
by $T$ means that the tested face is not wholly below $L$.
The retained length is therefore below $a-L+2\rho+1$ for an
input chain below $a$. This argument is independent of the degree
and of the number of simplices. An increase of this constant
absorbs strict-cut slacks.

For a positive-degree cocycle $S$, the cap--cup adjunction gives
\[
 \langle S,D_T^Lz\rangle
 =\bigl\langle e^*U\smile A^*T\smile B^*S,c_z\bigr\rangle.
\]
The right side denotes a supported cup on the common small-chain
complex: $A^*T$ and $B^*S$ are evaluated only after a nonzero
Thom coefficient has placed the parent simplex in the return tube.
It is not a cup of globally extended cochains on $\Omega\times I$.
The Thom pullback is zero on both interval ends. Moreover
$T\times S$ is a relative cocycle for
$(\Omega\times\Omega,\{q\}\times\Omega\cup\Omega\times\{q\})$,
since both degrees are positive and the cochains are normalized.
Thus the expression is the scalar pairing of the relative Thom
construction, with the interval providing the degree shift by one.
On the exact return locus the map $(A,B)$ is restriction to the
two subintervals, up to the specified constant time slots. In the
tube its short joining paths are the usual stick-insertion
extension of that cut. Relative Thom naturality and the cap--cup
adjunction therefore identify the expression with the based GH
product of \cite[Proposition~9.4]{GH2009}. The small-chain version
of these homological naturality statements is described in
\cite[Appendix~A]{HW2023}; it is used here only for this scalar
comparison, not for the square identity.
Subdivision preserves the relative input class. Changing the slots
preserves the relative target on the interval ends, so the associated
timing homotopy does not change this pairing. This proves
\eqref{eq:based-1}.
The argument identifies this positive-output pairing only; it
makes no inference about a degree-zero output from relative
cohomology.
\end{proof}

\subsubsection*{Boundary faces of a product of chains}
To construct the square we shall take products of two input chains and
two intervals. On its boundary, some maps cease to depend on one input.
This has two useful consequences: one edge becomes a product with the
unused input, and another edge vanishes after normalization. The
following elementary shuffle calculation makes both statements exact
at chain level.

\begin{lemma}[Shuffle caps and an omitted factor]
\label{lem:shuffle-omitted-factor}
Work over $\mathbb F_2$ with the ordinary shuffle cross product,
normalized cochains of degrees $n,p$, and normalized output chains.
Let $z$ be a chain of pure degree $j$ in one factor $X$, and let
$c$ be a chain in a product $Y$ of the other factors. Assume all
product simplices lie in the small-chain complex of
Lemma~\ref{lem:supported-cap-chain}.
\begin{enumerate}
\item Suppose $e$ and $A$ depend only on $Y$, and
$B(x,y)=m(b(y),x)$ for a continuous map $m$.
Then the supported double cap on $z\times c$, followed by $B$,
is exactly the shuffle product of the supported double cap on
$c$, followed by $b$, with $z$, followed by $m$ and normalization.
The analogous identity holds with $m(x,b(y))$.
\item Suppose $e$ and $B$ depend only on $Y$, although $A$ may
depend on both factors. If $j>p$, the normalized supported double
cap on $z\times c$ is zero.
\end{enumerate}
In part~(1), use the same triangulated $Y$ factors on the two sides;
no subsequent subdivision of a whole product simplex is allowed.
\end{lemma}
\begin{proof}
It suffices to consider a product of individual simplices. A shuffle
simplex is specified by a monotone lattice path, one step at a time
in one of the factors. Its projection to any collection of factors
is the affine simplicial map specified by the projected vertices.
If a step belongs to an omitted factor, the two projected vertices
coincide and this affine map factors through the corresponding
simplicial degeneracy. This is an actual degeneracy of singular
simplices, not merely equality of two vertices of an arbitrary
singular map.

The first $n$ steps form the front face evaluated by $U$, the next
$p$ steps form the face evaluated by $T$, and the remaining steps
form the output back face. In part~(1), an $X$-step in either of
the first two blocks makes the corresponding cochain evaluation
zero. Thus all $X$-steps occur in the back block. Removing them
specifies a unique shuffle in the $Y$ factors, with its two capped
front blocks, and the residual interleaving specifies precisely a
shuffle of its back face with the full $X$ simplex. Conversely,
every such pair of shuffles gives one surviving product shuffle.
The cochain coefficients and the output maps agree term by term.
The support condition causes no extra terms: since $e$ ignores
$X$, the relevant return set is pulled back from $Y$. This proves
part~(1), including its version with the output order reversed.

For part~(2), there are $j$ steps in $X$, but only $p$ slots in
the middle block. Since $j>p$, an $X$-step occurs in the first or
last block. In the first block its omission by $e$ makes the
normalized $U$ evaluation zero. In the last block its omission by
$B$ makes the output a degenerate singular simplex. Every summand
therefore vanishes after normalization. The argument applies to
subdivided factor simplices as well: their dimensions, and the
factorwise shuffle description, are unchanged.
\end{proof}

The shuffle calculation controls the algebraic boundary. We next
control lengths on the same simplices. A cochain is evaluated on a
face, whereas the retained output comes from another face of its
parent simplex. Uniform continuity on a sufficiently fine finite mesh
is what transfers the lower bound on the tested piece to that output.

\begin{lemma}[Finite meshes and the support estimate]
\label{lem:finite-mesh-support}
Consider finitely many product parameter domains, with their faces and
specified timing homotopies. Suppose their evaluation maps are continuous,
their cutting maps $A,B$ are continuous on an outer return tube, and
each output pair satisfies
\[
 \ell(A)+\ell(B)\leq a+b+2\rho
\]
on a two-input domain (replace $a+b$ by $a$ on a one-input domain).
Here the lengths of the joining arcs are at most $\rho$.
Choose a closed inner support tube contained in the outer tube.
For every $\varepsilon>0$ there are compatible factorwise refinements
such that every product cell meeting the support lies in the outer
tube and the oscillation of $\ell(A)$ on that cell is less than
$\varepsilon$. A nonzero test by $T$ on a face of such a cell implies
that its retained output has length at most
$a+b-L+2\rho+\varepsilon$.
\end{lemma}
\begin{proof}
The preimages of the closed inner tube are compact. Each has a
neighborhood whose closure maps into the outer tube; on that closure
the cutting maps and their lengths are uniformly continuous.
Choose a finite cover on which the oscillation of $\ell(A)$ is less
than $\varepsilon$, together with open sets missing the inner tube.
A sufficiently fine product mesh refines this cover. Use a common
number of subdivisions on repeated input simplices and their faces,
and a common interval subdivision on repeated parameter edges.
Subdivision commutes with the boundary, so these choices agree on
faces. Apply the shuffle triangulation only after these factorwise
refinements. Timing intervals are included before choosing the mesh.

If the tested face were entirely in $\Omega^{<L}$, its evaluation by
$T$ would be zero. Thus a nonzero tested face contains a point with
$\ell(A)\geq L$. At every point of its parent cell,
$\ell(A)\geq L-\varepsilon$. Subtracting this from the displayed
sum gives the retained-output bound. This proves a support estimate
for contributing terms; it does not construct a chain map or prove
that the cap constructions glue across the cells.
\end{proof}

It is also useful to record how a chain homotopy changes when its
boundary operations are replaced by chain-homotopic ones. This
observation explains the form of the correction terms. In the main
construction we make all choices simultaneously; the lemma by itself
does not supply a comparison between independently chosen cuts.

\begin{lemma}[Changing a realization in the square identity]
\label{lem:cut-change-correction}
Suppose degree $-Q$ chain maps $D_0^{L/R},D_1^{L/R}$ are related by
degree $1-Q$ maps $C_T^{L/R}$ satisfying
\[
 D_1^{L/R}+D_0^{L/R}=\partial C_T^{L/R}+C_T^{L/R}\partial.
\]
If $H_1$ has boundary $D_1^R(y)*x+y*D_1^L(x)$ for cycles $x,y$,
then
\[
 H_0=H_1+C_T^R(y)*x+y*C_T^L(x)
\]
has boundary $D_0^R(y)*x+y*D_0^L(x)$.
If each $C_T^{L/R}$ sends chains below $a$ to chains below $a-L+K$,
and $H_1$ has the corresponding two-input support estimate,
their sum is supported below $a+b-L+K$ after taking the maximum
of the three constants.
\end{lemma}
\begin{proof}
Over $\mathbb F_2$, the boundary of a product is the sum of the two
product boundaries. Since $x,y$ are cycles,
$\partial C_T^{L/R}(x)=(D_1^{L/R}+D_0^{L/R})(x)$, and similarly
for $y$. The terms involving $D_1$ cancel in the asserted formula.
The support of a sum is contained in the union of the supports;
the bounds on different summands are therefore not added.
The same argument applies to any finite sequence of such homotopies
with a common bound. This lemma is an algebraic statement, not an
assertion that these comparison homotopies exist. The filling proof
below avoids needing them by choosing its whole finite diagram at once.
\end{proof}

\begin{remark}[The meaning of a uniform constant]
Fix an outer tube with joining-arc bound $\rho$ and take the mesh
oscillation below $1$. Lemma~\ref{lem:finite-mesh-support} gives
$2\rho+1$ for the raw cut and square terms, independently of the
number and dimension of their simplices. A product prism restoring
an input has the bound of its cut factor plus the unchanged input
length. Lemma~\ref{lem:cut-change-correction} handles changes of
realization if controlled comparison chains have been constructed.
The simultaneous realization is provided by the explicit timing
homotopies and product prisms in Lemma~\ref{lem:based-cut-square}.
Compactness supplies the mesh, while the shuffle and cap identities
supply compatibility; these are distinct parts of that proof.
\end{remark}

\subsubsection*{The square and its correction terms}
We can now put the preceding ingredients together. Cut each of two
loops once and exchange their remaining pieces. The two cut parameters
form a square. Two boundary edges vanish by the omitted-factor lemma;
the other two give the desired products after their time slots have
been adjusted. Subdivision prisms then restore the original inputs.
Because the two reconnected outputs together traverse each input only
once, the resulting length bound is additive.

\begin{lemma}[The based cut square]
\label{lem:based-cut-square}
Choose the one-input cuts on the common finite diagram as in
Proposition~\ref{prop:one-input-cuts}. Let $x,y$ be cycles of degrees $j,N$, supported below $a,b$. If
$j>p$ and $N>p$, there is a chain $H_T(x,y)$ such that
\begin{equation}\label{eq:based-3}
\partial H_T(x,y)=D_T^R(y)*x+y*D_T^L(x),
 \qquad
 \operatorname{supp}H_T(x,y)\subset\Omega^{<a+b-L+K}.
\end{equation}
Here $*$ is based concatenation, and a larger universal additive
allowance may be used for the finite list of timing and cap comparisons.
\end{lemma}

\begin{figure}[tbp]
\centering
\begin{tikzpicture}[x=0.9cm,y=0.9cm,>=stealth,font=\small]
\definecolor{cutblue}{RGB}{20,100,150}
\definecolor{cutorange}{RGB}{175,70,20}
\node[anchor=west,font=\small\bfseries] at (-.4,3.1) {(a) Cut and reconnect};
\coordinate (q) at (0,0);
\coordinate (xs) at (4.4,.65);
\coordinate (yt) at (4.4,-.65);
\draw[cutblue,thick,->] (q) .. controls (0,2.6) and (3.7,2.7) .. (xs)
 node[pos=.49,above] {$\xi|_{[0,s]}$};
\draw[cutorange,thick,->] (xs) .. controls (2.7,1.7) and (1.3,1.1) .. (q)
 node[pos=.51,below] {$\xi|_{[s,1]}$};
\draw[cutorange,thick,->] (q) .. controls (1.4,-1.2) and (2.8,-1.6) .. (yt)
 node[pos=.46,above] {$\eta|_{[0,t]}$};
\draw[cutblue,thick,->] (yt) .. controls (3.7,-2.6) and (0,-2.6) .. (q)
 node[pos=.49,below] {$\eta|_{[t,1]}$};
\draw[cutblue,thick,dashed,->] (xs) .. controls (4.95,.3) and (4.95,-.3) .. (yt);
\draw[cutorange,thick,dashed,->] (yt) .. controls (4,-.3) and (4,.3) .. (xs);
\fill (q) circle (1.6pt) node[left] {$q$};
\fill (xs) circle (1.6pt) node[above right] {$\xi(s)$};
\fill (yt) circle (1.6pt) node[below right] {$\eta(t)$};
\node[align=left,anchor=west] at (4.95,0) {short\\arcs};
\node[cutblue,anchor=west] at (-.2,-3.1) {$A$: tested};
\node[cutorange,anchor=west] at (2.9,-3.1) {$B$: retained};
\draw[black!20] (6.2,-3.25)--(6.2,2.8);
\node[anchor=west,font=\small\bfseries] at (6.6,3.1) {(b) Parameter square};
\draw[thick] (7,-1.55)--(7,1.65)--(10.5,1.65);
\draw[thick,dashed] (10.5,1.65)--(10.5,-1.55)--(7,-1.55);
\foreach \x in {7,10.5} {\foreach \y in {-1.55,1.65} {\fill (\x,\y) circle (1.4pt);}}
\node[above] at (8.75,1.85) {$y*D_T^L(x)$};
\node[rotate=90,above] at (6.8,.05) {$D_T^R(y)*x$};
\node[rotate=90,below] at (10.7,.05) {$0\quad(j>p)$};
\node[below] at (8.75,-1.75) {$0\quad(N>p)$};
\node at (8.75,1.28) {$t=1$};
\node at (8.75,-1.15) {$t=0$};
\node[rotate=90] at (7.35,.05) {$s=0$};
\node[rotate=90] at (10.15,.05) {$s=1$};
\node[align=center,text width=2.2cm] at (8.75,.25) {cap with $U$,\\then $T(A)$};
\node[cutorange] at (8.75,-.6) {retain $B$};
\node[align=center,font=\footnotesize] at (8.75,-2.7)
 {Dashed sides vanish in\\normalized chains.};
\node at (5.3,-3.95) {$\ell(A)+\ell(B)\leq\ell(\xi)+\ell(\eta)+2\rho$};
\end{tikzpicture}

\caption{The two based outputs and the boundary of the parameter square.
The blue output is tested by $T$ and the orange output is retained.
The zero-edge assertions require the strict degree inequalities.}
\label{fig:based-cut-square}
\end{figure}
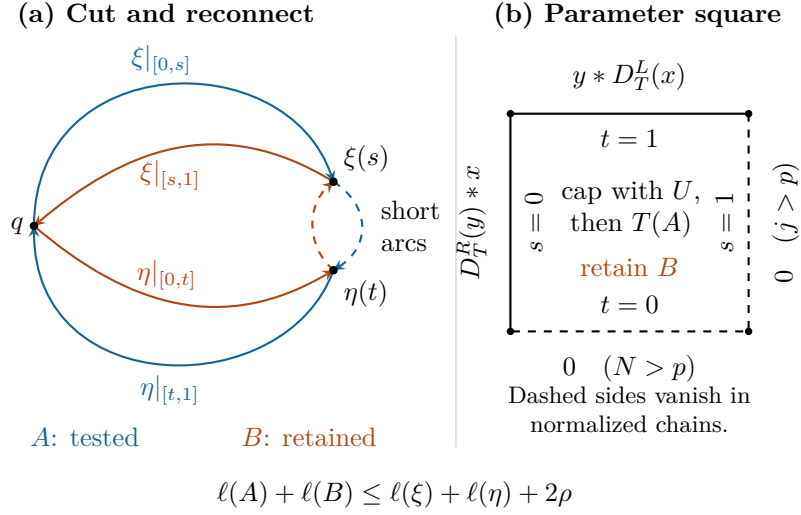

\begin{proof}
We construct all terms on one finite diagram. In particular, the
cuts used on the edges will be exactly those of
Proposition~\ref{prop:one-input-cuts}, with one common subdivision
order, rather than independently chosen representatives of the
same homology operation.

\emph{The geometric maps and the raw boundary.}
On $P=\Omega\times\Omega\times I_s\times I_t$ set
\[
 e(\xi,\eta,s,t)=(\xi(s),\eta(t)).
\]
On its outer return tube define
\[
 A=\xi|_{[0,s]}\#[\xi(s),\eta(t)]\#\eta|_{[t,1]},\qquad
 B=\eta|_{[0,t]}\#[\eta(t),\xi(s)]\#\xi|_{[s,1]}.
\]
Assign each of these three pieces one third of the total time,
including a constant piece when a restriction interval collapses.
These maps are continuous into the $H^1$ based loop space.
Restriction and rescaling are continuous at collapsed intervals
by the energy estimate above; the short joining geodesic depends
smoothly on its endpoints in the outer tube.

Specify the two edge timing homotopies described below before
choosing a mesh. Include all input faces, the one-input diagrams,
and all squares requested in the given finite family. Choose one
$k$ so that their factorwise $k$-fold subdivisions, followed by
shuffle triangulations, satisfy both smallness and length
oscillation less than $1$ in Lemma~\ref{lem:finite-mesh-support}.
Write $\widetilde x=\mathrm{Sd}^k\iota x$,
$\widetilde y=\mathrm{Sd}^k\iota y$, and $J=\mathrm{Sd}^k[I]$.
Set
\[
 H^{\rm raw}=\pi\mathcal K_{U,T}
       (\widetilde x\times\widetilde y\times J_s\times J_t).
\]
It has degree $j+N+2-n-p=j+N-Q+1$. The input chains are cycles
and $\partial J=[1]+[0]$. Lemma~\ref{lem:supported-cap-chain}
therefore identifies its boundary with the sum of its four actual
edge caps; it has no interior boundary term.

\emph{The two zero edges.}
At $s=1$, both $e$ and $B$ ignore $\xi$. The last piece of $B$
is a constant path with a fixed time slot, also independent of
$\xi$. Part~(2) of Lemma~\ref{lem:shuffle-omitted-factor}, with
$j>p$, makes this edge zero. At $t=0$, $e$ and $B$ ignore
$\eta$, including its constant prefix, so the same argument with
$N>p$ makes that edge zero. These are equalities of normalized
chains. No cancellation in homology or transversality assertion
is being substituted for the required chain equality.

\emph{The edge at $s=0$.}
Here $A$ consists of a constant piece, $[q,\eta(t)]$, and the
suffix of $\eta$, each initially assigned time $1/3$. Vary the
three times linearly, with parameter $\tau\in[0,1]$, from
$(1/3,1/3,1/3)$ to $(0,1/2,1/2)$.
For $B$, whose pieces are the prefix of $\eta$, $[\eta(t),q]$,
and $\xi$, vary the times from $(1/3,1/3,1/3)$ to
$(1/4,1/4,1/2)$. Thus at $\tau=1$ the maps are exactly
\[
 A=R(\eta,t),\qquad B=L(\eta,t)*\xi,
\]
where a one-input cut uses two equal slots and $*$ assigns one
half to each of its two factors. Only a constant piece has its
slot shrunk to zero, and all nonconstant slots stay bounded away
from zero. Consequently the timing homotopy is $H^1$-continuous.
It preserves both output lengths.

Use these maps, and $e=(q,\eta(t))$ independent of $\tau$, in
Lemma~\ref{lem:supported-cap-chain} on
$\widetilde x\times\widetilde y\times J_t\times J_\tau$.
Denote the resulting normalized chain by $C_0$. At $t=0,1$,
$e$ is constant at $(q,q)$, so the normalized positive-degree
Thom cocycle pulls back to zero. Input boundaries vanish. The
only boundary terms of $C_0$ are therefore the raw $s=0$ edge
and the edge with the displayed final timing. Part~(1) of
Lemma~\ref{lem:shuffle-omitted-factor} identifies the latter
\emph{exactly} as
\[
 D_T^R(y)*\mathcal S_kx.
\]
Normalization commutes with this cross product because shuffle
products with a degenerate factor are degenerate. The $y$ and
$t$ factors, their subdivision order, and the cochains are exactly
those used to define the chosen $D_T^R(y)$.

\emph{The edge at $t=1$.}
The three pieces of $A$ are the prefix of $\xi$, $[\xi(s),q]$,
and a constant piece. Vary their times from $(1/3,1/3,1/3)$ to
$(1/2,1/2,0)$. The pieces of $B$ are $\eta$, $[q,\xi(s)]$,
and the suffix of $\xi$; vary their times to $(1/2,1/4,1/4)$.
The final maps are $A=L(\xi,s)$ and $B=\eta*R(\xi,s)$.
The same supported cap on the edge times its timing interval
gives a chain $C_1$ whose boundary is the raw $t=1$ edge plus
\[
 \mathcal S_ky*D_T^L(x).
\]
The $s$-endpoint terms vanish by the same constant-evaluation
argument. In particular, no unconstructed comparison between two
independent cut operators is used on either edge.

\emph{Restoring the input chains.}
It follows that $H'=H^{\rm raw}+C_0+C_1$ satisfies
\[
 \partial H'=D_T^R(y)*\mathcal S_kx+
                    \mathcal S_ky*D_T^L(x).
\]
Define
\[
 H_T(x,y)=H'+D_T^R(y)*\mathcal P_kx+
                       \mathcal P_ky*D_T^L(x).
\]
Since $x,y$ are cycles and the chosen cuts are chain maps, the
identity for $\mathcal P_k$ cancels the two subdivided-input terms
and proves the boundary formula in~\eqref{eq:based-3}.

\emph{Uniform support and simultaneous choices.}
At every point of the square, and throughout both timing
homotopies,
\[
 \ell(A)+\ell(B)\leq\ell(\xi)+\ell(\eta)+2\rho.
\]
On a contributing product cell, the tested face contains a loop
of length at least $L$, since $T$ vanishes below $L$. The common
mesh makes every tested loop in that cell have length greater
than $L-1$. Thus every retained output in $H^{\rm raw},C_0,C_1$
has length less than $a+b-L+2\rho+1$.
The two restoration prisms have the same bound: each cut factor
saves $L-(2\rho+1)$, and $\mathcal P_k$ stays in the image of its
uncut input. A sum takes the union of supports. One fixed increase
of $K=2\rho+1$ absorbs any strict-cut slack.

All timing formulas are fixed before subdivision. A finite number
of requested squares adds only finitely many compact parameter
domains, so the same $k$ works for all of them and for the
one-input cuts. Repeated inputs therefore receive identical cuts.
The radius $\rho$ and the oscillation allowance are fixed
independently of $p,j,N$ and of the number of requested squares.
This completes the simultaneous finite-diagram construction.
\end{proof}

\begin{remark}
The construction proves the supported square directly from singular
caps and explicit timing homotopies. The unfiltered open-string
relation of Cieliebak--Oancea
\cite[Definition~6.1 and Theorem~9.2]{CO2026}, with locators in
arXiv:2209.02972v2, is consistent with the boundary formula but is
not an input to its proof. Nor is an unfiltered comparison used to
infer its metric support bound.
\end{remark}

\subsection{Pontryagin classes and Goresky--Hingston powers}

The cut construction applies to a general simply connected manifold.
To obtain a degree-independent estimate, we now use the particular
based-loop algebra of a CROSS. We need classes whose powers have two
compatible interpretations: concatenation restores degree, and the
GH product tests the amount removed by a cut.

The unfiltered operations depend on the diagonal Thom class and the
based loop space. Common smaller tubes and homotopies of the joining
paths compare two metric choices on each finite chain diagram. We may
therefore use the symmetric metric to identify the products. Their
length minimax values will always refer to the given metric.

For $M=\mathbb CP^m,\mathbb HP^m$, $m\ge2$, put $d=2,4$, $n=dm$,
$e=d-1$, and $\nu=n+d-2$. For $\mathbb OP^2$ put
$(n,e,\nu)=(16,7,22)$. The based Pontryagin algebra is
\begin{equation}\label{eq:based-4}
H_*(\Omega;\mathbb F _2)=\Lambda(u_e)\otimes\mathbb F _2[y_\nu].
\end{equation}
The product in \eqref{eq:based-4} is induced by concatenation of based loops; we use
$*$ for it at chain level and juxtaposition for homology classes. The
subscripts record homological degrees. Thus $u\in H_e(\Omega)$,
$y\in H_\nu(\Omega)$, $u^2=0$, and the nonzero groups are spanned by
$y^k$ and $uy^k$ for $k\ge0$. The algebra $\Lambda(u_e)$ is the
exterior algebra, including the relation $u^2=0$ in characteristic two.
The relation $u^2=0$ and the nonvanishing of every power of $y$ will
both be used. We therefore give an algebra calculation, including the
multiplicative extensions, rather than only counting the graded groups.

\begin{lemma}[The based Pontryagin algebra]\label{lem:based-algebra}
For the projective spaces above, equation~\eqref{eq:based-4} is an
isomorphism of algebras over $\mathbb F_2$.
\end{lemma}
\begin{proof}
Write
\[
 A=H^*(M;\mathbb F_2)=\mathbb F_2[x]/(x^{m+1}),\qquad |x|=d,
\]
where $(d,m)=(2,m),(4,m)$ with $m\ge2$, or $(8,2)$ for the
Cayley plane. Thus $e=d-1$ and $\nu=(m+1)d-2$.
We use the multiplicative cobar spectral sequence supplied by Adams'
published cobar theorem~\cite{Adams1956}; a modern multiplicative
formulation is given in \cite[Section~3]{MR2024}. Filter the cobar
complex by word length. In word-length and internal
degrees $(s,t)$ its second page is
\[
 E_2^{s,t}=\operatorname{Cotor}_{H_*(M;\mathbb F_2)}^{s,t}
                 (\mathbb F_2,\mathbb F_2)
       =\operatorname{Ext}_{A}^{s,t}(\mathbb F_2,\mathbb F_2),
 \qquad E_\infty^{s,t}\Longrightarrow H_{t-s}(\Omega M;\mathbb F_2).
\]
The product converges to the Pontryagin product. One may use a
1-reduced chain model, since $M$ is simply connected. Every
nonempty desuspended factor then has positive degree, so word length
is bounded in each total degree; this gives convergence with a finite
filtration in each degree. Finite type identifies the coalgebra
cobar calculation with the displayed Ext calculation.

We compute its algebra, including the square of the odd resolution
class. A free resolution of the augmentation module has one generator
$z_s$ in every resolution degree, with
\[
 |z_{2k}|=k(m+1)d,\qquad |z_{2k+1}|=k(m+1)d+d,
\]
where these are internal degrees, and differential
\[
 \partial z_{2k+1}=x z_{2k},\qquad
 \partial z_{2k}=x^m z_{2k-1}\quad(k\ge1).
\]
It is exact because the annihilators of $x$ and $x^m$ in $A$
are $(x^m)$ and $(x)$, respectively. Applying
$\operatorname{Hom}_A(-,\mathbb F_2)$ gives zero differential.
The two-step shift of this periodic resolution represents a class
$\bar y$ of bidegree $(2,(m+1)d)$, and its powers identify all even
resolution classes. A one-step lift of the class $\bar u$ of
bidegree $(1,d)$ is given on $P_{s+1}\to P_s$ by multiplication
by $1$ for even $s$ and by $x^{m-1}$ for odd $s$. These maps
commute with the resolution differential. Their square, followed
by augmentation, is zero since $m\ge2$. Hence $\bar u^2=0$.
The one-step and two-step lifts commute, and $\bar u\bar y^k$
evaluates nontrivially on $z_{2k+1}$. Consequently
\[
 E_2=\Lambda(\bar u)\otimes\mathbb F_2[\bar y].
\]

A higher differential has bidegree $(r,r-1)$ for $r\ge2$ and
lowers total degree by one. No such differential can leave
$\bar y$: there is no group of total degree $\nu-1$. None can
leave $\bar u$: total degree $e-1$ is absent, except for the unit
when $e=1$, and the unit has word length zero rather than $1+r$.
Multiplicativity now makes every differential zero, page by page,
since these classes generate the second page.

Finally there are no hidden algebra extensions. The total degrees
$k\nu$ and $e+k\nu$ are all distinct; each has one associated
graded generator. Choose lifts $y,u$. The products $y^k,uy^k$
are nonzero because their leading classes are nonzero. The group
in degree $2e$ is zero, since $0<2e<\nu$, so $u^2=0$.
The difference $uy+yu$ has zero associated graded class and there
is no other nonzero filtration quotient in its total degree;
hence it is zero. These products therefore give the asserted
algebra, not only its associated graded vector space. No formality
assumption on the mod-two cochains of $M$ is used: possible higher
cobar differentials have been ruled out by the preceding degree argument.
\end{proof}

Both generators are primitive for the ordinary diagonal: $u$ has the
least positive degree and $\nu>2e$ leaves no possible reduced coproduct
for $y$. Hence the ordinary class $\beta\in H^e(\Omega)$ dual to $u$
satisfies
\begin{equation}\label{eq:based-5}
\beta\cap(uy^k)=y^k.
\end{equation}
For completeness, the diagonal of $uy^k$ has exactly one summand with
first tensor factor of degree $e$, namely $u\otimes y^k$. Pairing
$\beta$ with that first factor is the ordinary cap product in \eqref{eq:based-5}.
This cap product does not use the Goresky--Hingston multiplication.

We write $\theta\in H^e(\Omega)$ for the class dual to
$u$ when using this product; as an ordinary cohomology class it is
the same class as $\beta$. The two symbols distinguish the two
operations in which that class will be used.

For the reference symmetric metric, let $\lambda_r$ be the Morse index
at the critical level of the $r$th iterates of prime geodesics,
$r\ge1$. In the sphere case use $e=n-1$ and $\nu=2(n-1)$ in this
calculation, as in the algebra described below.
Goresky--Hingston \cite[(13.1.1) and Section~14.9]{GH2009} give the
same based and free indices
$\lambda_r=r\lambda_1+(r-1)(n-1)$ and based critical spheres
$S^{n-1}$. Here $\lambda_1>0$; the zero-index case in
\cite[Proposition~13.2]{GH2009} is real projective space.
The first based local bottom
survives: the constants split off, and all subsequent indices exceed
$\lambda_1+1$. It is therefore the first positive based homology class,
so \eqref{eq:based-4}, or the sphere algebra below, identifies $\lambda_1=e$.
Consequently $\lambda_r=e+(r-1)\nu$.

We use the symmetric-metric Morse perfection over $\mathbb F_2$.
This coefficient assertion follows by repeating the proof of
\cite[Sections~13.3--13.5]{GH2009}, rather than by inferring it
from integral perfection. Evaluation splits the constants and lifts
the fundamental class of the first critical manifold to the first
sublevel. Its absolute loop product carries each successive local
Thom group isomorphically to the next: under the bundle identification
this is intersection with the critical fundamental class. All Thom
maps and that coefficient-one intersection are valid over $\mathbb F_2$.
Starting at the first crossing, surjectivity from sublevel homology
to the local group therefore propagates by induction, so every
connecting map vanishes. Field duality gives cohomological perfection.
For based loops the local sphere groups have exactly the degreewise
ranks of \eqref{eq:based-4} and its sphere analogue below, which also forces perfection.
Thus no change-of-coefficients assertion about perfection is needed.

Goresky--Hingston \cite[Proposition~9.4 and Section~14.9]{GH2009}
identify the restriction of the free GH product with the based GH
product and its associated graded map with
$H^*(SM)[T]\to H^*(S^{n-1})[T]$.
Their Section 14.6 identifies the relevant Thom class with a product
of Thom classes by an actual bundle isomorphism; this calculation
also holds over $\mathbb F_2$. The bottom coefficient is the
restriction $1\mapsto1$.
It follows that the iterated product $\theta^{\circledast r}$, with
$r\ge1$, satisfies
\begin{equation}\label{eq:based-6}
\theta^{\circledast r}\text{ is dual to }uy^{r-1},\qquad
 |\theta^{\circledast r}|=e+(r-1)\nu=r\nu-(n-1).
\end{equation}
Multiplication by this bottom class carries the first top class, dual
to $y$, to the top class dual to $y^{r+1}$, again with coefficient one.
Every degree in question has rank one, so associated-graded detection
proves these exact scalar identities: there is no second class in
that degree with which the detected class could be mixed.

For $M=S^n$, $n\ge2$, put $e=n-1$, $\nu=2e$, $u=v$, and $y=v^2$ in
$H_*(\Omega S^n;\mathbb F _2)=\mathbb F _2[v]$.
Here $u$ is not exterior. Formula \eqref{eq:based-6} follows from the same calculation
in \cite[Section~14.9]{GH2009},
and \eqref{eq:based-5} follows from
$\beta\cap v^{2k+1}=(2k+1)v^{2k}=v^{2k}$ over $\mathbb F _2$.
The sphere cases include $\mathbb CP^1$ and $\mathbb HP^1$.

We return to the given metric. The algebraic calculation now provides
classes whose critical lengths can be compared. The key point is that
each relevant degree has only one nonzero class: a cohomology test
becomes nonzero at exactly the level at which its dual homology class
can first be represented.

Set $c_k=c(y^k)$, $c_0=0$, $t_r=t(\theta^{\circledast r})$,
$H=c(u)$, and $Y=c_1$. Ordinary cap and multiplication by $u$ give
\begin{equation}\label{eq:based-7}
c_k\le c(uy^k)\le c_k+H.
\end{equation}
Indeed, capping a representative of $uy^k$ with $\beta$ does not change
its loop-space support, whereas concatenating representatives of $u$
and $y^k$ adds their lengths. Taking infima gives the two inequalities.
Field duality in the rank-one degree of \eqref{eq:based-6} gives
$t_r=c(uy^{r-1})$. To see the filtration point, a cohomology class
vanishes on a subspace in cohomology exactly when it lifts to relative
cohomology and hence admits a cocycle vanishing there. Over a field,
its restriction is zero exactly when it pairs trivially with the image
of homology. In the rank-one degree under consideration the latter
image is either zero or the whole group. Subadditivity of $c_r$ then
implies
\begin{equation}\label{eq:based-8}
c_{r-1}\le t_r\le c_{r-1}+H,\qquad c_r-t_r\le Y.
\end{equation}
The GH product makes $t_r$ superadditive. Here the filtered based
construction is the based-return version of \cite[Proposition~9.5
and Section~10.1, (10.1.1)]{GH2009}, just as Proposition~9.4 replaces
free loops by based loops. Explicitly, when the two tested cocycles
vanish below $L_1,L_2$, a nonzero term has two outputs of total length
at least $L_1+L_2$; the input length differs from that total by at most
the arbitrarily small return-arc allowance. Shrinking the return tube
and taking the supremum of vanishing thresholds gives the claimed
inequality without a fixed loss. Fekete's lemma and \eqref{eq:based-8}
therefore give a common limit. More explicitly, the subadditive
sequence $c_r$ has slope $\inf_r c_r/r$, and the superadditive sequence
$t_r$ has slope $\sup_r t_r/r$. The inequalities
$c_{r-1}\le t_r\le c_{r-1}+H$ show that these slopes agree; call
their common value $\mu$. Now $c_r\ge r\mu$ and $t_r\le r\mu$.
Combining each of these inequalities with $c_r-t_r\le Y$ gives
\begin{equation}\label{eq:based-9}
\mu=\inf_{r\ge1}\frac{c_r}{r}=\sup_{r\ge1}\frac{t_r}{r},
 \qquad r\mu\le c_r\le r\mu+Y,\quad r\mu-Y\le t_r\le r\mu.
\end{equation}
Since short based loops contract in a convex ball, $t_1>0$, and hence
$0<t_1\le\mu\le Y$.

\subsection{Recovering the unit and reducing the filling degree}
\label{subsec:recover-unit}

We now apply the cut construction to fillings. Cutting by a suitable
GH power lowers the degree by a multiple of $\nu$; multiplication by
$y^r$ restores it. The estimates just proved show that the restoration
cost exceeds the saving by at most a fixed constant. It remains to
show that the constructed chain has the original cycle as boundary.
The square reduces this question to the degree-zero output of the
right cut, which must represent the unit.

All the cuts in this subsection use the simultaneous realization of
Proposition~\ref{hyp:controlled-cut-realization}. We first recover the
unit, then recall the fixed-degree filling estimate, and finally
combine the two constructions.

Fix $r\ge1$, and let $T$ represent $\theta^{\circledast r}$; its degree
is $p=r\nu-(n-1)$ and the degree lowered by the cut is $Q=r\nu$.
When homology classes occur as inputs to a chain operation, choose
cycle representatives and use the same symbols for those representatives.
By \eqref{eq:based-1}, \eqref{eq:based-6}, and the
bottom-times-top coefficient,
\begin{equation}\label{eq:based-10}
[D_T^L(y^{r+1})]=y.
\end{equation}
Apply the square to $x=y^{r+1}$ and $y_{\rm input}=y^r$.
The degree inequalities hold since $p=r\nu-(n-1)<r\nu$.
Writing $[D_T^R(y^r)]=a\cdot1$ in $H_0(\Omega)=\mathbb F _2$,
equation \eqref{eq:based-3} and \eqref{eq:based-10} give $a y^{r+1}=y^{r+1}$.
Thus
\begin{equation}\label{eq:based-11}
[D_T^R(y^r)]=1.
\end{equation}
The equality is in $H_0(\Omega)$, which is one-dimensional because
$M$ is simply connected. It identifies the unit directly from the
square, so it does not require an additional choice of absolute lift
of a relative operation.

\subsubsection*{Fixed-degree fillings}
We next require an additive filling estimate in each of finitely many
degrees. The estimate must shorten a \emph{prescribed} filling relative
to its boundary. Existence of an unrelated short disk would not give
the cell-by-cell compatibility needed for a singular filling chain.

For fixed-degree filling we use the relative disk theorem of
Nabutovsky--Rotman \cite[Theorem~7.5]{NR2013}, printed p.~56 of the
author manuscript, and the finite-cell induction in
Block--Manin--Weinberger \cite[Theorem~4.1(ii)]{BMW2026}, printed
p.~11 of arXiv:2505.16907v2. We spell out the prescribed-class
qualification needed by that induction. Put $d_0=\operatorname{diam}M$.
Choose $S_0\ge0$ such that every based loop $\gamma$ of length at most
$2d_0$ contracts through based loops of length at most
$\ell(\gamma)+S_0$. Thus $S_0$ is an upper bound for the short-loop
depth in \cite[Definition~7.4]{NR2013} (equivalently use
\cite[Theorem~7.3]{NR2013} with $k=1$). Define
\[
 B_0=5d_0+S_0,\qquad A_0=6d_0+2S_0,
 \qquad T_m=A_0m-d_0-S_0.
\]
Such a finite $S_0$ exists for every fixed metric here. The regularized
parametrization below, with $\delta=d_0$, moves each loop of length at
most $2d_0$ without changing length to a loop with Lipschitz constant
at most $3d_0$. These loops belong to the compact set of based loops
with Lipschitz constant at most $3d_0$, in the uniform topology.
Each of its loops has a finite-length contraction: polygonize it,
contract the polygon by a finite piecewise smooth disk using simple
connectedness, and restore the short-geodesic collar. A sufficiently
small uniform neighborhood is joined to its center by pointwise short
geodesics. The differential of this interpolation is uniformly bounded
on a compact diagonal neighborhood, so its length is bounded by a fixed
multiple of the two input lengths. That neighborhood therefore inherits
a finite contraction bound. A finite cover gives one common bound;
enlarging $S_0$ gives the required length-plus-$S_0$ estimate.
This argument does not assert compactness of all parametrized
bounded-length loops in the strong $H^1$ topology.

Here the integer $m$ denotes the dimension of a boundary sphere and is
unrelated to the dimension parameter of a projective space.
For $m\ge1$, a boundary map $f:S^m\to\Omega^{<a}$ can first be shortened
to length $T_m+\varepsilon$, through a homotopy of length at most
$a+A_0m+\varepsilon$. Extend this boundary motion across a \emph{given}
filling disk by the homotopy extension property. Its new boundary lies
below $T_{m+1}$ after choosing $\varepsilon<A_0$, since
$T_{m+1}-T_m=A_0$. The relative disk clause of
\cite[Theorem~7.5]{NR2013} shortens
this particular disk, fixing that boundary, to length
$T_{m+1}+\varepsilon$. Restoring the boundary cylinder gives a filling
in the prescribed relative homotopy class of length at most
\begin{equation}\label{eq:based-12a}
\max\{a,B_0\}+A_0m+\varepsilon.
\end{equation}
For $m=0$, apply the $S^0$ clause of \cite[Theorem~7.5]{NR2013}
to each endpoint together with
the constant loop, and concatenate the two endpoint tracks. This
moves both endpoints to the constant loop with cost at most
$a+2d_0+S_0+\varepsilon$; extend that
motion across the \emph{given} interval. Apply the relative disk theorem
in dimension one to the resulting interval with constant endpoints,
then restore its endpoint collars. This preserves its prescribed path
class and gives the safe bound
\[
 \max\{a+2d_0+S_0,B_0\}+\varepsilon\le a+B_0+\varepsilon.
\]
Existence of some short endpoint path alone would not suffice here,
since $\pi_1\Omega M=\pi_2M$ need not vanish.

In a finite $D$-dimensional relative cell induction, first move the free
vertices to the constant loop. The one-cells then require at most $B_0$
extra length, and an $r$-cell with $r\ge2$ requires at most $A_0(r-1)$
more. Thus the final map has length below
$a+B_0+A_0\binom D2+\varepsilon$ and represents the original relative
class. A finite chain and its prescribed boundary have the requisite
finite cell presentation. First use the natural normalization splitting
to lift normalized chains to ordinary singular chains; its added
degeneracies stay in the images of their parent simplices. This is a
chain map, so cycles and prescribed boundaries remain exact. Over
$\mathbb F_2$, pair equal canceling faces of the finite filling and
attach its remaining faces to a finite presentation of the boundary.
A mapping cylinder, if necessary, makes that boundary a subcomplex
without increasing dimension beyond $D$. The maps on all identified
faces agree. This gives a finite CW pair, which is all the cell
induction requires; a manifold realization is unnecessary.
For a globally null $b$-cycle use $D=b+1$ and absorb the positive slack
in $C_b=B_0+A_0\binom{b+1}2+1$. It bounds below $a+C_b$.
The finitely many auxiliary approximation errors in a cell induction
are chosen with total less than the displayed slack. Thus the slack
is not multiplied by the number of cells.
This estimate concerns the given Riemannian metric; no bilipschitz
change of metric is used to transfer an additive estimate. Only the
finitely many constants $C_0,\ldots,C_{\nu-1}$ will be needed.

\subsubsection*{Uniform fillings in all degrees}
We now reduce any remaining degree modulo $\nu$. The essential estimate
is $c_r-t_r\le Y$: the length added by a representative of $y^r$
exceeds the length removed by its dual cut by a constant independent
of $r$.

Let $z$ be a globally null cycle of degree $j\ge\nu$ below $a$.
Choose $r=\lfloor j/\nu\rfloor$ and $b=j-r\nu\in\{0,\ldots,\nu-1\}$.
Then $p=r\nu-(n-1)<r\nu\le j$.
Choose $0<\varepsilon<t_1/4$, a representative of $y^r$ below
$c_r+\varepsilon$, and a cocycle $T$ vanishing on
$C_*(\Omega^{<t_r-\varepsilon})$. Since superadditivity gives
$t_r\geq rt_1>0$, the threshold $L=t_r-\varepsilon$ is positive,
as required in Proposition~\ref{hyp:controlled-cut-realization}.

\smallskip\noindent\emph{Fixing the entire input diagram before cutting.}
Choose also a finite cycle representing $y^{r+1}$ and a finite
global filling $B$ with $\partial B=z$. Request the cuts on these
chains and their faces and request the two squares with input pairs
$(z,y^r)$ and $(y^{r+1},y^r)$. This is a finite request made before
any cut output is constructed. Apply
Proposition~\ref{hyp:controlled-cut-realization} once to this diagram.
Its second square and the positive-output pairing prove
$[D_T^R(y^r)]=1$ for precisely this realization, by the computation
in \eqref{eq:based-10}--\eqref{eq:based-11}. Its chain-map identity
on $B$ proves $D_T^Lz=\partial D_T^LB$.
No bound on $B$ or on the auxiliary $y^{r+1}$ cycle is needed:
their outputs are used only to identify global homology classes.
The support bounds below are applied to $z$, $y^r$, and their square.
All later low-degree fillings are taken after this choice and are
never used as inputs to a cut. Thus no enlargement or comparison
of already chosen cut operators is hidden in the argument.

\smallskip\noindent\emph{Filling the cut cycle.}
Since the cut is a chain map,
$D_T^Lz$ is globally null, has degree $b$, and lies below
$a-t_r+K+\varepsilon$. Fill it using $C_b$ and concatenate the filling
on the left with $y^r$. In other words, if $W$ is the resulting
low-degree filling, then
\[
 \partial W=D_T^Lz,\qquad
 \partial(y^r*W)=y^r*D_T^Lz.
\]
The degree of $y^r*W$ is $r\nu+(b+1)=j+1$, as required for a
filling of $z$. Equation \eqref{eq:based-8} bounds its support by
\begin{equation}\label{eq:based-12}
a+c_r-t_r+K+C_b+2\varepsilon
 \le a+Y+K+\max_{0\le b<\nu}C_b+2\varepsilon.
\end{equation}
If the displayed input level for the cut is nonpositive, that cut chain
is zero and this filling step is unnecessary.
The square gives a chain, within the same type of bound, joining this
boundary to $D_T^R(y^r)*z$. By \eqref{eq:based-2}, the zero-cycle $D_T^R(y^r)$ has
uniformly bounded length support, and by \eqref{eq:based-11} its difference from the
constant unit is globally null. More precisely, \eqref{eq:based-2} places this
zero-cycle below $c_r-t_r+K+2\varepsilon\le Y+K+2\varepsilon$.
The degree-zero filling theorem therefore gives a chain $V_r$ with
\[
 \partial V_r=D_T^R(y^r)+1
\]
and a length bound independent of $r$. The boundary of $V_r*z$ is
$D_T^R(y^r)*z+1*z$. Adding this chain to the square and to $y^r*W$
leaves precisely $1*z$ as boundary; all intermediate terms occur
twice and cancel over $\mathbb F_2$. Finally, the length-preserving
timing prism from the fixed-slot product $1*z$ to $z$ restores the
original cycle as boundary. Each summand has support below $a$ plus
one constant independent of $j$ and $a$.
Degrees $j<\nu$ are covered by the same finite list of $C_b$.

\subsection{Length and the Hilbert energy filtration}

The preceding construction gives length bounds, because concatenation
adds lengths. Morse theory is applied instead to $F=\sqrt E$ on the
strong $H^1$ loop space. Although $\ell\leq F$, this inequality alone
does not turn a length-controlled filling into an energy-controlled
one. We therefore change the speed of each loop so that $F$ is as
close to its length as needed. A positive regularization makes this
change continuous even when the loop has constant subarcs. The
homotopy back to the original parametrization will preserve the
boundary data of a filling.

Fix $\delta>0$. For $\gamma\in\Omega$, set
$v(t)=|\dot\gamma(t)|$ and $L=\ell(\gamma)$, and, for $s\in[0,1]$, define
\[
 a_{\delta,s}(t)=(1-s)+s\frac{v(t)+\delta}{L+\delta},\qquad
 \phi_{\delta,s}(t)=\int_0^t a_{\delta,s}(r)\,dr,
 \qquad R_{\delta,s}\gamma=\gamma\circ\phi_{\delta,s}^{-1}.
\]
Since $\int_0^1a_{\delta,s}=1$ and
\[
 a_{\delta,s}(t)\ge\min\{1,\delta/(L+\delta)\}>0
 \quad\hbox{for almost every }t,
\]
$\phi_{\delta,s}$ is an absolutely continuous increasing
homeomorphism of $[0,1]$ with Lipschitz inverse. Consequently
$R_{\delta,s}\gamma$ is well defined as an $H^1$ loop. It preserves
based endpoints, image, and length, and $R_{\delta,0}$ is the identity.
The change-of-variables formula gives
$E(R_{\delta,s}\gamma)=\int_0^1v^2/a_{\delta,s}$; in particular,
\begin{equation}\label{eq:based-13}
E(R_{\delta,1}\gamma)
 =(L+\delta)\int_0^1\frac{v^2}{v+\delta}
 \le L(L+\delta)<(L+\delta)^2.
\end{equation}
Moreover $E(R_{\delta,s}\gamma)\le E(\gamma)$ for all $s$.
To see this, $h(v)=v^2/(v+\delta)$ is increasing, so the covariance
inequality $\int vh(v)\ge(\int v)(\int h(v))$ proves the assertion
at $s=1$. Explicitly,
\[
 \int v^2=\int (v+\delta)h(v)
 \ge (L+\delta)\int h(v)=E(R_{\delta,1}\gamma).
\]
For each $t$, the reciprocal of the positive affine function
$a_{\delta,s}(t)$ is convex in $s$. Thus
$E(R_{\delta,s}\gamma)\le(1-s)E(\gamma)+sE(R_{\delta,1}\gamma)$,
which proves the estimate throughout the homotopy.

\smallskip\noindent\emph{Continuity in the Hilbert topology.}
The map $(\gamma,s)\mapsto R_{\delta,s}\gamma$ is jointly continuous
for the strong $H^1$ topology.
For an $H^1$-convergent sequence the speeds converge in $L^2$, the
lengths converge, and the maps $\phi$ converge uniformly. Locally their
derivatives have a common positive lower bound, so their inverses
converge uniformly. The reparametrized loops therefore converge
uniformly. The energy formula is continuous by convergence in measure,
the common positive denominator bound, and uniform integrability of
the squared speeds. In a fixed isometric Euclidean embedding this
identifies the weak $H^1$ limit and proves convergence of derivative
norms, hence strong $H^1$ convergence. No exact-arclength limit is taken.

Thus a length filling $B$ of an $F<a$ cycle $z$, supported below $a+D$,
can be replaced by $R_{\delta,1}B$, supported in $F<a+D+\delta$.
The prism $R_{\delta,s}z$ stays in $F<a$ and corrects its boundary back
to $z$. Since $\ell\le F$, the $F$ minimax is at least the length
minimax; applying $R_{\delta,1}$ to any representative and letting
$\delta\downarrow0$ gives the reverse inequality. The minimax values
therefore agree. This last passage compares numbers; it does not take
a limit of the reparametrization maps in the strong $H^1$ topology.

\smallskip\noindent\emph{Comparison with the quantitative filling model.}
Before applying Nabutovsky--Rotman's piecewise differentiable model to a compact
$H^1$ input, first polygonize that input if necessary and restore its
prescribed boundary by the collar below. If a cited quantitative theorem
supplies a compact family of Lipschitz
loops only in the uniform topology, use one common sufficiently fine
time subdivision and replace each subarc by its short geodesic chord.
Compactness gives equicontinuity, so all chords lie in convex balls.
The polygonal family is $H^1$-continuous and has no greater length.
On the prescribed $H^1$ boundary family the subdivision can be chosen
so that the standard polygonal homotopy stays within any prescribed
positive length slack. Indeed, polygonization converges uniformly in
strong $H^1$ on a compact strong-$H^1$ family: sequentially, if
$\gamma_i\to\gamma$ strongly and the meshes tend to zero, the polygons
converge uniformly, their energies are at most $E(\gamma_i)$, and weak
$H^1$ compactness together with the energy upper bound gives strong
convergence. The short-geodesic homotopy is then uniformly $H^1$-small.
Its prism corrects the polygonal boundary.

Preserving the specified relative class requires one additional
compatibility. Polygonize the entire shortening homotopy with one common time
mesh. It fixes the polygonized original boundary. Attach the same
boundary-restoration collar throughout. At its initial end the original
disk has compact strong-$H^1$ image, so its own polygonal approximation
is connected to it by the preceding small homotopy, including the
boundary trace. This identifies the collar-modified initial disk with
the prescribed disk relative to its exact boundary. The final shortened
disk therefore has the specified relative class in the strong-$H^1$
model. We need a length bound only on the final disk and the
boundary-restoration collar, not on the homotopy from an arbitrarily
long initial interior filling. Apply $R_{\delta,s}$ afterwards to pass
from this length estimate to the energy estimate.

We have now established both ingredients announced at the beginning
of the section. The cut-and-restore construction controls classes
which vanish globally, while the minimax comparison controls classes
which survive. We collect them in the form used by the next section.

\begin{proposition}[Uniform based-loop death and resonance]
\label{prop:based-uniform-persistence}
For a smooth Riemannian metric on a simply connected CROSS,
the following properties
hold with constants $D,C>0$ and $\lambda>0$,
uniformly in all degrees $j$ and levels $a$:
\begin{enumerate}
\item every globally null cycle in $\Omega_q^{F<a}M$ bounds in
      $\Omega_q^{F<a+D}M$;
\item the image of $H_j(\Omega_q^{F<a}M)$ in $H_j(\Omega_qM)$ is zero
      for $a<\lambda j-C$, and is surjective for $a>\lambda j+C$.
\end{enumerate}
The same statements hold with $F$ replaced by $\ell$.
The first conclusion uses Lemma~\ref{lem:based-cut-square}; the second
uses only the algebra and minimax estimates of the preceding subsections.
\end{proposition}

\begin{proof}
For length sublevels, the first assertion follows from the construction
leading to \eqref{eq:based-12}, including its degree-zero correction. Only the
finite set $C_0,\ldots,C_{\nu-1}$ occurs, so a single constant works
in every degree. The regularized change of clock and its boundary
prism give the same assertion for $F$ after enlarging that constant.

For the second assertion put $\lambda=\mu/\nu>0$.
The nonzero based homology groups are one-dimensional in degrees
$k\nu$ and $k\nu+e$. In the first row \eqref{eq:based-9} gives
$k\mu\le c(y^k)\le k\mu+Y$. In the second row, \eqref{eq:based-7} gives
\[
 k\mu\le c(uy^k)\le k\mu+Y+H,
 \qquad \lambda(k\nu+e)=k\mu+\lambda e.
\]
Thus every nonzero class has minimax within
$C=Y+H+\lambda e$ of $\lambda j$. A map into a one-dimensional
homology group has zero image before the minimax and is surjective
after it. The strict inequalities in the statement avoid any issue
of attainment. The zero-degree class is born at zero, and a zero
global group imposes no condition. Equality of length and $F$ minimax
values proves the two versions of the assertion.
\end{proof}

\section{Passing from based loops to free loops}
\label{sec:promotion}
\label{sec:finite-cover}

Closed geodesics are critical points on the free loop space, so the
based-loop estimates must be transferred there. Evaluation at the
initial point describes a free loop as a loop based at a varying
point. Over a small convex set, we can move that point to a fixed
basepoint by adding a path and its reverse. These extra paths, often
called whiskers, have uniformly bounded length.

Theorem~\ref{thm:promotion} implements this idea over a finite convex
cover of $M$. Its main issue is gluing: in a Mayer--Vietoris sequence,
a homology preimage need not have a representative at the same length
level as the original class. This is why we carry the filling estimate
and the spectral estimate together. We first prove the abstract union
lemma, then apply it to the evaluation cover. All homology in this
section has coefficients in $\F$.

For a space $X$ with increasing open sublevels $X^a$, write
$\mathrm{UD}(D)$, or \emph{uniform death with increment $D$}, if
\[
 \ker\bigl(H_q(X^a)\longrightarrow H_q(X)\bigr)
 \longrightarrow H_q(X^{a+D})
 \quad\hbox{is zero for every }q,a.
\]
Thus a class in a sublevel which vanishes in the whole space already
vanishes after the same increase $D$ in every degree. Write
$\mathrm{SR}(\lambda,C)$, or \emph{spectral resonance with slope
$\lambda$ and width $C$}, if the map $H_q(X^a)\to H_q(X)$ is zero
for $a<\lambda q-C$ and surjective for $a>\lambda q+C$. Constants are
independent of the degree $q$ and level $a$. Negative homology groups
are zero. For a subspace $A\subset X$, its filtration is always the
induced one, $A^a=A\cap X^a$.

The two estimates must be propagated together. Uniform death controls
fillings of classes already known to be zero, whereas spectral
restriction controls the level of a nonzero Mayer--Vietoris preimage.
That preimage has no preferred representative at the original level.

\begin{lemma}[Filtered union lemma]\label{lem:union}
Suppose $X=A\cup B$ is an open cover and its restrictions cover every
sublevel. If $A,B,I=A\cap B$ have UD and SR with a common slope
$\lambda>0$, then so does $X$.
\end{lemma}
\begin{proof}
We use the ordinary Mayer--Vietoris sequence and its natural maps from
the corresponding sequences of sublevels. The part relevant in degree
$q$ is
\[
 H_q(I)\longrightarrow H_q(A)\oplus H_q(B)
 \longrightarrow H_q(X)\xrightarrow{\delta}H_{q-1}(I).
\]
Over $\F$ the two inclusions in the first arrow have the same sign.
Set $D_0=\max(D_A,D_B)$ and $C_0=\max(C_A,C_B)$, taking all constants
nonnegative. All occurrences of $\epsilon>0$ below can be chosen as
small as necessary; the final added constant $1$ absorbs them.

\smallskip\noindent\emph{Uniform death.}
Let $z\in C_q(X^a)$ be a cycle null in $X$. By subdivision subordinate
to $A^a\cup B^a$, we may write
$z=z_A+z_B$ with $\partial z_A=\partial z_B=b\in C_{q-1}(I^a)$.
The subdivision prism stays in $X^a$ and will restore the original
cycle after this argument. Naturality of the connecting map gives
$[b]=\delta[z]=0$ in $H_{q-1}(I)$.

Uniform death in $I$ supplies a chain $w\in C_q(I^{a+D_I+\epsilon})$
with $\partial w=b$. Put $z'_A=z_A+w$ and $z'_B=z_B+w$. These are
cycles, and their sum is $z$ because $w+w=0$. Write
$\alpha=[z'_A]\in H_q(A)$ and $\beta=[z'_B]\in H_q(B)$.
Their sum maps to zero in $H_q(X)$, so exactness supplies a class
$\eta\in H_q(I)$ mapping to $(\alpha,\beta)$.

If $(\alpha,\beta)=0$, both pieces fill with an additional cost $D_0$.
Otherwise at least one is a nonzero global class represented below
$a+D_I+\epsilon$. The lower spectral restriction in that member gives
$\lambda q\le a+D_I+C_0+\epsilon$. The upper spectral restriction in
$I$ now represents the \emph{chosen} preimage $\eta$ by a cycle $v$
below
\[
 \lambda q+C_I+\epsilon
 \le a+D_I+C_0+C_I+2\epsilon.
\]
The cycles $z'_A+v$ and $z'_B+v$ have zero global classes, so uniform
death in $A$ and $B$ fills them after a further increase $D_0$.
Adding the two fillings gives a filling of $z$. For example,
\[
 D_X=D_I+C_0+C_I+D_0+1
\]
works. This is the step for which the lower spectral restriction is
essential: it bounds the level needed to choose $v$.

\smallskip\noindent\emph{The upper spectral restriction.}
Let $\xi\in H_q(X)$. Choose a cycle $b$ representing
$\delta\xi\in H_{q-1}(I)$ below
$\lambda(q-1)+C_I+\epsilon$. Exactness says that the images of $b$
in $A$ and $B$ are globally null. Uniform death in these two spaces
gives fillings $a_A,a_B$ at that level plus $D_0$. Their sum is a
cycle in $X$; let $\xi'=[a_A+a_B]$. Its connecting image is
$\delta\xi'=\delta\xi$.

The difference $\xi-\xi'$ therefore has a preimage in
$H_q(A)\oplus H_q(B)$. Represent its two components below
$\lambda q+C_0+\epsilon$ using the upper restrictions in $A$ and $B$.
Adding those cycles to $a_A+a_B$ represents $\xi$ below the larger
of the two bounds. Hence one can take
\[
 C_X^+=\max\{C_0,C_I-\lambda+D_0\}+1.
\]
\smallskip\noindent\emph{The lower spectral restriction.}
Suppose a cycle in $X$ is supported below $a<\lambda q-C_X^-$.
Its connecting cycle in $I$ has degree $q-1$ and the same length
bound. If $C_X^->C_I+\lambda$, then
$a<\lambda(q-1)-C_I$, so the connecting class is globally zero.
Uniform death in $I$ closes the two subdivided pieces at level
$a+D_I+\epsilon$, as in the first part of the proof. If also
$C_X^->C_0+D_I$, choose $\epsilon$ so small that this level is below
$\lambda q-C_0$. The two closed pieces then vanish in $H_q(A)$ and
$H_q(B)$ by their lower spectral restrictions. Their sum consequently
vanishes in $H_q(X)$. Thus one can take
\[
 C_X^-=\max\{C_I+\lambda,C_0+D_I\}+1.
\]
Use $C_X=\max(C_X^+,C_X^-)$. For $q=0$, the connecting group
$H_{-1}(I)$ is zero; all steps involving its representative or filling
are omitted, and the same bounds apply.
\end{proof}

The union lemma preserves the slope and enlarges only additive
constants. It can therefore be iterated over a finite cover, provided
the same hypotheses are known on every nonempty intersection.

\begin{corollary}[Finite-cover promotion]\label{cor:cover}
If every nonempty finite intersection in a finite open cover has UD and
SR with a common slope, the whole space has UD and SR with that slope.
\end{corollary}
\begin{proof}
We induct on the number $N$ of cover members, making the assertion
simultaneously for every cover with at most $N$ members. The assertion
for $N=1$ is the hypothesis. Suppose it is known for $N-1$ members,
and write $X=U_1\cup\cdots\cup U_N$. Set
\[
 A=U_1\cup\cdots\cup U_{N-1},\qquad B=U_N.
\]
The induction hypothesis applies to $A$. It also applies to
$A\cap B$, covered by the nonempty sets $U_i\cap U_N$ for $i<N$:
every multiple intersection of these sets is a multiple intersection
of the original cover. Thus $A$, $B$, and $A\cap B$ satisfy both
properties with slope $\lambda$, and Lemma~\ref{lem:union} applies.
Empty members are discarded. Only finitely many constants are
introduced, so the resulting constants remain independent of degree.
\end{proof}

For loop spaces, the required local identifications are supplied by
whiskers. The estimate below depends on the cover and the chosen
connecting paths, but its increase in length is independent of the
dimension of a chain. That is exactly the uniformity needed to apply
the preceding corollary.

\begin{theorem}[Based-to-free promotion]\label{thm:promotion}
Let $M$ be a closed connected Riemannian manifold. If $\Omega_qM$,
filtered by length, has UD and SR with slope $\lambda>0$,
so does $\Lambda M$.
\end{theorem}
\begin{proof}
Choose a finite cover of $M$ by sufficiently small strongly convex balls.
Every nonempty finite intersection $V$ is convex. Choose $x_V\in V$,
a fixed piecewise smooth path $\sigma:q\to x_V$, and the short
geodesic $\eta_x:x_V\to x$
for $x\in V$. Put $\sigma_x=\sigma\#\eta_x$. These paths have uniformly
bounded length, including over the finite collection of intersections.
Here $\operatorname{ev}_0:\Lambda M\to M$ evaluates a loop at its
initial point, and $\bar\sigma$ denotes a path with reversed
parametrization. On $E_V=\operatorname{ev}_0^{-1}(V)$ use the maps
\[
 \begin{split}
 F_V:E_V\longrightarrow\Omega_qM,\qquad
 &F_V(\gamma)=\sigma_{\gamma(0)}\#\gamma\#\bar\sigma_{\gamma(0)},\\
 G_V:\Omega_qM\longrightarrow E_V,\qquad
 &G_V(\omega)=\bar\sigma\#\omega\#\sigma.
 \end{split}
\]
The latter loop is based at $x_V$. Use fixed time slots for every
concatenation, so the maps are continuous in the strong $H^1$ topology.
Cancellation of a path followed by its reverse gives a homotopy
$F_VG_V\simeq\mathrm{id}$. For $G_VF_V$, first
cancel the fixed $\bar\sigma\#\sigma$ backtracks. The remaining
$\eta_x\#\gamma\#\bar\eta_x$ moves its basepoint along $\eta_x$
while the two whiskers shrink. Its evaluation remains in $V$ throughout.
Both maps and both homotopies increase length by a constant $K_V$;
the variable loop is traversed once. For example, the length added
by $F_V$ is at most $2\sup_{x\in V}\ell(\sigma_x)$.
The corresponding bounds for cancellation and shrinking whiskers
depend only on the same fixed path lengths and the diameter of $V$.
Taking their maximum gives $K_V$.

These homotopies make $F_V$ and $G_V$ inverse maps on global homology.
For the lower spectral restriction, a class represented in $E_V$ below
$a$ maps to a based class below $a+K_V$; it vanishes if
$a+K_V<\lambda q-C$, and hence the original class vanishes. For the
upper restriction, represent the corresponding based class below
$\lambda q+C+\epsilon$ and apply $G_V$. Thus the width $C+K_V$,
with harmless positive slack, works in $E_V$ and the slope is unchanged.

For uniform death, send a null cycle below $a$ to the based space below
$a+K_V$, fill it after an increase $D$, and send its filling back.
The resulting boundary is $G_VF_V$ of the original cycle. The controlled
identity homotopy supplies the remaining prism, so the original cycle
dies below $a+D+3K_V+\epsilon$.
Thus every $E_V$ has both properties with the same slope. The sets
$E_V$ form a finite cover whose intersections have the same form.
Every open length cut is covered by their restrictions, so subdivision
preserves the cut. Apply Corollary~\ref{cor:cover}.
\end{proof}

We will study nonconstant critical loops. It is consequently natural
to pass to homology relative to the constant loops. Evaluation provides
a projection which performs this passage without adding length; the
same finite cover also gives a degree-independent bound on global
Betti numbers.

\begin{corollary}[Relative groups and uniform ranks]\label{cor:relative-analytic}
For a simply connected CROSS, the conclusions of
Theorem~\ref{thm:promotion} hold also for the pairs $(\Lambda M,M)$.
Their ordinary global Betti numbers are bounded independently of degree.
\end{corollary}
\begin{proof}
Evaluation at zero gives a retraction of every positive sublevel onto
its constant loops; no deformation retraction is asserted.
If $i:M\to\Lambda M$ is the constant-loop inclusion, then
$\operatorname{ev}_0i=\mathrm{id}_M$. The chain endomorphism
\[
 P=1-i_*\operatorname{ev}_{0*}:C_*(\Lambda M)\longrightarrow
 C_*(\Lambda M)
\]
is therefore a projection onto $\ker\operatorname{ev}_{0*}$ with
kernel $C_*(M)$. It induces the chain isomorphism
$C_*(\Lambda M)/C_*(M)\cong\ker\operatorname{ev}_{0*}$ and preserves
every positive length cut, since its correction term consists of
constant loops. A relative cycle and a relative filling can thus be
represented in this direct summand without increasing length.

Applying the same projection to absolute representatives and fillings
transfers both the death estimate and the two spectral restrictions
to relative homology without filtration loss. The regularized
reparametrization in the preceding analytic section fixes constant
loops and preserves evaluation, so the same statement holds for
square-root energy.

Each degree of the based homology of a CROSS has dimension at most one.
The ordinary Mayer--Vietoris inequality
\[
 b_q(A\cup B)\le b_q(A)+b_q(B)+b_{q-1}(A\cap B)
\]
propagates a bound independent of $q$ through the same finite-cover
induction. Removing the constant-loop summand preserves this bound.
\end{proof}

\section{Equivariant cohomology of the loop space of a sphere}
\label{sec:spheretop}
\label{sec:sphere-topology}

We turn from filtration estimates to the topology needed for the
multiplicity argument. The relevant distinction is whether a
reflection-equivariant class comes from full $O(2)$-equivariant
cohomology, and whether such a lift is killed by a power of $U$.
Proposition~\ref{prop:sphere} gives these classes for a sphere.

The round metric is useful because its critical manifolds and their
indices are explicit. We first compute the contribution of each
critical manifold, then show that its generators and required
relations hold in the full loop space. This second step is essential:
a local calculation determines an associated graded module, whereas
the final argument uses multiplication and restriction of global
classes. Once established, those global statements are independent
of the metric used to calculate them.

Throughout this section $M=S^n$ with $n\geq3$, and
\[
 b=n-1,\qquad \nu=2b,\qquad a_r=(2r-1)b\quad(r\geq1).
\]
The case $S^2=\mathbb CP^1$ is included in
Theorem~\ref{thm:projective-equivariant-topology}.
Write $\eta_n\in H^n(S^n;\F)$ for the nonzero class and put
$e=\operatorname{ev}_0^*\eta_n\in H_R^n(\Lambda S^n;\F)$.
The evaluation map is $R$-equivariant, since reflection fixes the
parameter zero; no $G$-equivariance of evaluation is needed.

\subsection{The two representation-sphere sequences}

There are two natural circle bundles linking the acting groups. The
standard plane representation relates full equivariance to reflection;
the representation which doubles rotation angles also detects
half-period translation. Their Euler classes are different, and this
difference will distinguish odd from even iterates.

We use the groups $G=O(2)$, $R=\langle\rho\rangle$, and
$K=\langle\rho,\tau\rangle$ from the introduction. Thus $\rho$ reverses
the loop parameter and $\tau$ translates it by one half. For the
coefficient rings, choose coordinates
\[
 H^*(BG;\F)=\F[w,U],\quad |w|=1,\quad |U|=2,
 \qquad H^*(BK;\F)=\F[w,v],\quad |v|=1.
\]
Here $v$ is the character which is zero on $\rho$ and one on $\tau$;
the character $w$ is one on $\rho$ and zero on $\tau$. The standard
real two-plane representation restricts to the sum of the characters
$v$ and $v+w$. Consequently restriction sends $U$ to $v(v+w)$.

For an $H$-space $Y$, write $Y_H=EH\times_HY$ for its Borel
construction. For a $G$-pair $(Y,A)$, the unit circle of the standard representation
gives the circle bundle $Y_R\longrightarrow Y_G$, and likewise on
$A$. Its mod-two Euler class is $U$. The Gysin sequence of this pair of
circle bundles is
\begin{equation}\label{eq:sphere-1}
\cdots\to H_G^{j-2}(Y,A)\xrightarrow{U}H_G^j(Y,A)
 \xrightarrow{\mathrm{res}}H_R^j(Y,A)
 \xrightarrow{\pi_!}H_G^{j-1}(Y,A)\xrightarrow{U}\cdots.
\end{equation}
The map $\pi_!$ is integration along the circle, with coefficients in
$\F$. The same construction for the representation which doubles
rotation angles has stabilizer $K$. Its Euler class is zero. Indeed,
its restriction to $K$ has a fixed line, whereas
$\F[w,U]\longrightarrow\F[w,v]$, $U\mapsto v(v+w)$, is injective.
To see this injectivity, the highest power $U^k$ of a nonzero
polynomial produces a nonzero term of degree $2k$ in $v$.
The second Gysin sequence therefore takes the form
\begin{equation}\label{eq:sphere-2}
0\to H_G^j(Y,A)\to H_K^j(Y,A)
 \to H_G^{j-1}(Y,A)\to0.
\end{equation}
Both sequences are natural for maps of pairs, in particular for
inclusions of energy sublevels and for the Thom attachments below.

\subsection{The reflection decomposition and its lifts}

For the round metric, successive critical manifolds are indexed by
the number of times a great circle is traversed. Each gives the same
reflection module, with a shift in degree. Their full-group lifts
behave differently according to parity. The following proposition
states the global form of this distinction.

\begin{proposition}[Sphere decomposition]\label{prop:sphere}
There are homogeneous classes $X_r\in H_R^{a_r}(\Lambda S^n,S^n)$
such that
\begin{equation}\label{eq:sphere-3}
H_R^*(\Lambda S^n,S^n)
 =\bigoplus_{r\ge1}\F[w,e]/(w^n,e^2)\,X_r.
\end{equation}
Every $X_r$ has a $G$-lift $\zeta_r$. For odd $r$ it may be chosen
with $U^n\zeta_r=0$. For even $r$, $eX_r$ also has a $G$-lift.
Thus, if $W$ is the $\F$-span of $w^kX_r$ for odd $r$ and
$0\le k<n$, and $E$ is the sum of the even blocks, then
\[
 H_R^*(\Lambda S^n,S^n)=W\oplus eW\oplus E.
\]
Every class of $W$ has a $U$-torsion $G$-lift, and every class of $E$
has a $G$-lift.
\end{proposition}
\begin{proof}
We organize the proof around the passage from local to global
cohomology. The index gaps first give a range where a global supported
class is determined by its local term. We then compute the local
modules and prove that their connecting maps vanish. This supplies
the generators. Finally, the same degree ranges and the Gysin map
establish the relations, including a common $U$-annihilator at odd
stages.

\paragraph{The round filtration.}
For the round metric, the $r$th critical manifold $C_r$ of nonconstant
loops consists of the $r$-fold parametrized great circles and is
identified with $STS^n$. A point $(x,u)$ specifies initial position and
unit tangent vector. The normal Fourier modes of frequencies $|k|<r$
give Morse index $(2r-1)(n-1)=a_r$; the zero modes give precisely the
tangent space of $C_r$. Choose a regular sublevel $Z_r$ just above
$C_r$, and choose $Z_0$ below the first positive critical value. The
latter retracts equivariantly to the constant loops.

For $H=G,R$, put
\[
 F_r^H=H_H^*(\Lambda S^n,Z_{r-1}),\qquad
 L_r^H=H_H^*(Z_r,Z_{r-1}).
\]
A class in $F_r^H$ is called supported at stage $r$; its image in
$H_H^*(\Lambda S^n,S^n)$ vanishes on $Z_{r-1}$. Its image in $L_r^H$
will be called its local term. The equivariant Morse--Bott attachment
theorem and the mod-two Thom isomorphism give
\[
 L_r^H\cong H_H^{*-a_r}(C_r).
\]
All Thom maps are natural under subgroup restriction and ambient cup
product. Mod-two coefficients make these statements independent of
orientations of the negative bundles.

For clarity, the passage to the whole loop space is made before any
splitting is asserted. The long exact sequence of the triple has the
form
\[
 F_{r+1}^{H,j}\longrightarrow F_r^{H,j}
 \longrightarrow L_r^{H,j}\longrightarrow F_{r+1}^{H,j+1}.
\]
Successive local groups start in degrees $a_s\to\infty$. It follows
that $F_{r+1}^{H,j}=0$ for $j<a_{r+1}$. One can justify this directly
on the exhaustion: every singular simplex in the Borel space has
bounded energy, and restriction of singular cochains to a subcomplex
is surjective. Global cochains are therefore the inverse limit of the
finite-sublevel cochains. In each of degrees $j$ and $j-1$ the
cohomology tower is eventually constant, because all later Thom
groups start above $j+1$. Its derived inverse-limit term vanishes.
Thus the finite-stage vanishing passes to $F_{r+1}^{H,j}$.

In particular, $F_r^{H,j}\to L_r^{H,j}$ is injective for
$j<a_{r+1}$ and is an isomorphism for $j\leq a_{r+1}-2$.
After the connecting homomorphisms have been shown to vanish, every
local class will admit a supported lift, unique whenever its degree
is below $a_{r+1}$. We refer to $j-a_r$ as its \emph{offset}.

\paragraph{The local reflection module.}
Reflection sends $(x,u)$ to $(x,-u)$ and acts freely on $C_r$.
Consequently $(C_r)_R$ is homotopy equivalent to the real projective
tangent bundle $P(TS^n)$. Stable triviality of $TS^n$ and the
projective-bundle formula give
\[
 H_R^*(C_r)=\F[w,e]/(w^n,e^2).
\]
Here $w$ is the first Stiefel--Whitney class of the tautological real
line on $P(TS^n)$, and $e$ is pulled back from the base $S^n$.
Writing $x_r^R\in L_r^{R,a_r}$ for the Thom image of $1$, we obtain
the cyclic ambient module $L_r^R=\F[w,e]x_r^R/(w^n,e^2)$.

\paragraph{The odd and even full-group modules.}
Suppose first that $r$ is odd. The orbit space of great circles is
$\operatorname{Gr}_2(\mathbb R^{n+1})$. The ineffective stabilizer
of the $r$th parametrization is cyclic of odd order, whose classifying
space has no positive-degree mod-two cohomology. Thus
$H_G^*(C_r)$ is the cohomology of this Grassmannian, with $w,U$ the
Stiefel--Whitney classes of its tautological two-plane.
Define homogeneous polynomials $A_i\in\F[w,U]$ by
\[
 A_0=1,\quad A_1=w,\quad
 A_i=wA_{i-1}+UA_{i-2}\quad(i\geq2).
\]
They are the homogeneous components of $(1+w+U)^{-1}$. The
complement has rank $b$, so
\[
 H_G^*(C_r)=\F[w,U]/(A_{b+1},A_{b+2}).
\]
Generation by these characteristic classes can be seen from the
ordered-line flag space. It projects to $\mathbb RP^n$ with fiber
the projective space of the orthogonal complement of the first line,
so its cohomology is generated by the two line classes $s,t$.
Over the Grassmannian its projective-bundle basis is $1,s$, with
$s+t=w$ and $st=U$. Reducing a polynomial in $s,t$ to
$P(w,U)+sQ(w,U)$ and using uniqueness in that basis proves generation
on the Grassmannian; the complementary-bundle relations give the
displayed presentation. Let $x_r^G$ denote the local bottom Thom
class.

If $r$ is even, half-translation fixes every point of $C_r$. Hence
\[
 H_K^*(C_r)=\F[w,e]/(w^n,e^2)[v].
\]
Multiplication by $U=v(v+w)$ is injective: its leading coefficient
as a polynomial in $v$ is one, even though the coefficient ring has
zero divisors. Sequence \eqref{eq:sphere-2} embeds $H_G^*(C_r)$ in this module,
so $U$ is injective on $H_G^*(C_r)$ as well. Sequence \eqref{eq:sphere-1} now gives
\[
 H_G^*(C_r)/(U)\cong H_R^*(C_r).
\]
Choose homogeneous lifts of $1$ and $e$ through this quotient. They
generate over $\F[w,U]$: subtract a lift of the reduction modulo
$U$, divide the remaining class by $U$, and repeat in a degree two
lower. This process terminates in each degree. The corresponding
local Thom generators have offsets $0$ and $n$.

\paragraph{Vanishing of the connecting maps.}
For reflection, the connecting map is $\F[w,e]$-linear. It vanishes
on $x_r^R$ because $a_r+1<a_{r+1}$, and hence on the entire local
module. For $G$, odd local modules are generated at offset zero, and
even local modules at offsets zero and $n$. Since $n+1<2b$ when
$n>3$, the same degree argument applies.

For $n=3$ the even top generator has degree $a_r+3$, so its connecting
image has degree $a_r+4=a_{r+1}$. Restrict this image to the next
odd local module. Its target is the one-dimensional bottom group,
and restriction from $G$ to $R$ is injective on that group. Naturality
identifies its reflection restriction with the already vanishing
reflection connecting image. The next local term is therefore zero.
As $a_{r+1}<a_{r+2}$, the injectivity range proved above makes the
connecting image itself zero. In all cases the supported sequences
are consequently short exact in every degree.

Choose supported bottom lifts and denote their images in global
cohomology by $X_r$ for $R$ and $\zeta_r$ for $G$. Their degrees are
below the next index, so restriction sends $\zeta_r$ exactly to
$X_r$. The relation $e^2=0$ is an ambient relation. The relation
$w^nX_r=0$ follows from its local version, because $n<2b$ and the
supported restriction is injective in that degree. At even stages the
chosen top lift restricts exactly to $eX_r$ for the same reason.

These facts prove \eqref{eq:sphere-3} as a module identity. To prove independence,
take a nonzero finite relation and restrict to its earliest supported
stage; the coefficient of that stage must vanish in the local ring.
Remove it and continue. To prove spanning, successively subtract the
block representing the earliest nonzero local term of a class. This
terminates in every fixed degree because $a_r\to\infty$.

\paragraph{The global torsion relation at odd stages.}
Fix odd $r$. In the local full-group module, $A_bx_r^G$ is nonzero
because it restricts to $w^b x_r^R$. The recurrence and the relations
$A_{b+1}=A_{b+2}=0$ imply $UA_bx_r^G=0$. The reflection group in
offset $b+1=n$ is one-dimensional, generated by $e x_r^R$.
Exactness of \eqref{eq:sphere-1} therefore gives
\[
 \pi_!(e x_r^R)=A_bx_r^G.
\]
Let $\widetilde\zeta_r\in F_r^{G,a_r}$ be the supported class defining
$\zeta_r$. Applying the Gysin map in the supported pair and
restricting to the local term gives the preceding equality.
Both sides have offset $b<2b$, where supported restriction is
injective. Thus
\[
 \pi_!\bigl(e\,\mathrm{res}\widetilde\zeta_r\bigr)
       =A_b\widetilde\zeta_r,
 \qquad UA_b\widetilde\zeta_r=0.
\]
Likewise $A_{b+1}\widetilde\zeta_r=0$, since its offset $b+1$ is
strictly below $2b$. Finally the recurrence gives the polynomial
identity
\[
 A_i^2+A_{i-1}A_{i+1}=U^i\qquad(i\geq1).
\]
For $i=1$ this is $w^2+(w^2+U)=U$; substitution of the recurrence
propagates it from $i$ to $i+1$. At $i=b$, multiply the identity by
$U\widetilde\zeta_r$. Both terms on the left vanish by the two
relations just proved. Hence $U^{b+1}\widetilde\zeta_r=0$, and
passing to global cohomology gives $U^n\zeta_r=0$.
The argument uses only offsets below $2b$, so it includes $n=3$.
Multiplying these lifts by powers of $w$ proves the assertion for
$W$, and the two even-stage lifts prove the assertion for $E$.
\end{proof}

\section{Equivariant cohomology of projective loop spaces}
\label{sec:projective-topology}
\subsection{The reflection module and the image of full equivariance}
\label{sec:projective-equivariant-topology}

The sphere calculation suggests what to seek on a projective space:
odd critical blocks should occur in pairs related by an evaluation
class, and the two members should have different lifting behavior.
The cohomology generator of the projective space supplies that class.
Theorem~\ref{thm:projective-equivariant-topology} identifies the paired
module and determines the full restriction image, including its
$U$-torsion part.

The method again starts with the symmetric metric. Two additional
points require attention. In the complex case, a local relation reaches
the next critical index, so degree counting alone does not make it a
global relation. For the Cayley plane, the odd local calculation uses
the bundle of projective lines and its characteristic classes. We
include both arguments, since the final proof needs the resulting
relations in global cohomology.

All cohomology in this subsection has coefficients in $\mathbb F_2$.
Let $G=O(2)$ act on the free loop space $\Lambda M$ by reparametrization,
and let $R\cong C_2$ be generated by $t\mapsto-t$.
Write
\[
 A=H^*(BG)=\mathbb F_2[w,U],\qquad |w|=1,\quad |U|=2.
\]
Here $w$ and $U$ are the Stiefel--Whitney classes of the standard real
two-plane. Under restriction to $R$, $w$ remains $w$ and $U$ becomes zero.
Free-loop cohomology groups are relative to the constant loops; ambient
classes and critical-manifold cohomology rings are absolute. Products
of an ambient class and a relative class mean the usual relative cup
product.

For $M=\mathbb K P^m$, where $\mathbb K=\mathbb C,\mathbb H$, put
$d=\dim_{\mathbb R}\mathbb K$, $n=dm$, and
\[
 \nu=n+d-2,\qquad \lambda_r=d-1+(r-1)\nu\quad(r\geq1).
\]
For $M=\mathbb O P^2$, put $d=8$, $n=16$, $\nu=22$, and
$\lambda_r=7+22(r-1)$. These are the successive Morse indices for the
standard symmetric metric. The corresponding critical manifolds are
copies $C_r$ of the unit tangent bundle $SM$.

Let $\eta\in H^d(M;\mathbb F_2)$ be the nonzero generator. The
evaluation map at zero is $R$-equivariant, with trivial $R$-action on
$M$, so
\[
 h=\operatorname{ev}_0^*\eta\in H_R^d(\Lambda M;\mathbb F_2)
\]
is an ambient class. We also write $h$ for its restriction to critical
manifolds. The use of $h$ is confined to reflection-equivariant or
ordinary cohomology; a full-group lift of a class involving $h$ will
always be constructed explicitly.

We introduce the rings which will describe the answer. Reflection
forgets the orientation of a tangent vector, leading to the projective
tangent bundle and its ring $\mathcal R$. Full equivariance at odd
stages instead forgets the parametrization of the geodesic. Its
restriction image will be the subring denoted by $\mathcal S$.
The formulas below make the statement precise; the proof derives them
from these geometric models.

For the complex and quaternionic cases, let $a,z$ be formal variables
of degrees $d,2d$, respectively, and define homogeneous polynomials
$Q_j\in\mathbb F_2[a,z]$ by
\[
 Q_0(a,z)=1,\quad Q_1(a,z)=a,\quad
 Q_j(a,z)=aQ_{j-1}(a,z)+zQ_{j-2}(a,z)\quad(j\geq2).
\]
\begin{align*}
 P_m(w,h)&=\sum_{i=0}^m\binom{m+1}{i}h^i w^{d(m-i)},\\
 \mathcal R_m&=\mathbb F_2[w,h]/(h^{m+1},P_m(w,h)),\\
 \mathcal S_m&=\mathbb F_2[w,z]/
       (Q_m(w^d,z),Q_{m+1}(w^d,z)).
\end{align*}
Here $P_m$ has degree $dm$ and the binomial coefficients are reduced
modulo two. These formulas initially define graded quotient rings;
their geometric meaning will be established in the proof.
For the Cayley plane set
\[
 \mathcal R=\mathbb F_2[w,h]/(w^{24},h^2+w^8h+w^{16}),
 \qquad
 \mathcal S=\mathbb F_2[w]/(w^{24}),\qquad |h|=8.
\]
When no distinction is necessary we write $\mathcal R,\mathcal S$ also
in the complex and quaternionic cases.

Fix regular sublevels $Z_r$ of the symmetric metric just above $C_r$,
and take $Z_0$ below the first positive critical value. Then $Z_0$
equivariantly retracts to $M$. For $H=G,R$, write
\[
 F_r^H=H_H^*(\Lambda M,Z_{r-1}),\qquad
 L_r^H=H_H^*(Z_r,Z_{r-1}).
\]
The map $F_r^H\to H_H^*(\Lambda M,M)$ is induced by the inclusion
of pairs. A class in its image is \emph{supported at stage $r$}.
The restriction $F_r^H\to L_r^H$ records its contribution at the
$r$th critical level and will be called its local term. A homogeneous
class of degree $\lambda_r+q$ at this stage has \emph{offset} $q$.
We shall show that the supported maps are injective, which will allow
us to identify a supported class with its global image.

\begin{theorem}[Projective equivariant topology]
\label{thm:projective-equivariant-topology}
Let $M=\mathbb C P^m$ or $\mathbb H P^m$, $m\geq1$, or
$M=\mathbb O P^2$. There are homogeneous classes
$X_r\in H_R^{\lambda_r}(\Lambda M,M)$, supported at and above the
$r$th symmetric critical level, for which the following assertions hold.
\begin{enumerate}
\item The reflection-equivariant module is
\[
 B:=H_R^*(\Lambda M,M)=\bigoplus_{r\geq1}\mathcal R X_r
 \quad\text{over }\mathbb F_2[w,h].
\]
The sum is locally finite in total degree, and $w^{2n}B=0$.
\item Under $z\mapsto h(h+w^d)$ in the complex and quaternionic cases,
\[
 \mathcal R=\mathcal S\oplus h\mathcal S.
\]
Set
\[
 W=\bigoplus_{r\text{ odd}}\mathcal S X_r,
 \qquad E=\bigoplus_{r\text{ even}}\mathcal R X_r.
\]
Then $B=W\oplus hW\oplus E$, and
\[
 \operatorname{im}\bigl(H_G^*(\Lambda M,M)\longrightarrow B\bigr)
       =W\oplus E.
\]
\item Let $T_U$ be the subgroup of $H_G^*(\Lambda M,M)$ annihilated
elementwise by some power of $U$. Then
\[
 \operatorname{res}(T_U)=W,\qquad U^N T_U=0,
 \quad
 N=\begin{cases}
 2m& M=\mathbb C P^m,\\
 4m& M=\mathbb H P^m,\\
 16& M=\mathbb O P^2.
 \end{cases}
\]
In particular every class of $W$ has a full-group lift killed by this
same power of $U$, and every class of $E$ has a full-group lift.
\end{enumerate}
No splitting of $H_G^*(\Lambda M,M)$ as an $A$-module is asserted.
\end{theorem}

\begin{proof}
There are three tasks. The local Morse--Bott calculation gives the
rings in the statement and shows that every local class has a supported
lift. We then establish the relations among the reflection lifts,
including the possible extension in the complex case. Finally we
construct full-group lifts of the odd and even generators and determine
their $U$-torsion. The last task also proves that the listed lifts
exhaust the restriction image.

\paragraph{Supported groups and the precise degree range.}
The equivariant Morse--Bott attachment and the mod-two Thom isomorphism
give $L_r^H\cong H_H^{*-\lambda_r}(C_r)$.
In particular each local group starts in degree $\lambda_r$. The
Thom isomorphism is compatible with subgroup restriction and with
cup product by classes pulled back from the ambient Borel space.
The triple sequence is
\begin{equation}
\label{eq:supported-triple}
 F_{r+1}^{H,j}\longrightarrow F_r^{H,j}
 \longrightarrow L_r^{H,j}
 \xrightarrow{\delta_r}F_{r+1}^{H,j+1}.
\end{equation}
Before any collapse assertion, $F_{r+1}^{H,j}=0$ for $j<\lambda_{r+1}$.
Thus restriction is injective for $j<\lambda_{r+1}$ and is an isomorphism
for $j\leq\lambda_{r+1}-2$. These are distinct ranges: existence of
a lift also involves the group in degree $j+1$, whereas uniqueness
only involves the group in degree $j$.
To justify the infinite exhaustion, every singular chain in its
Borel space has bounded energy. Global relative cochains are therefore
the inverse limit of finite-stage relative cochains, with degreewise
surjective restriction maps. The Thom groups vanish below their indices,
so the cohomology towers in both degrees $j$ and $j-1$ are eventually
isomorphic once the subsequent indices exceed $j+1$. The corresponding
$\varprojlim^1$ term is zero. Finite-stage vanishing therefore proves
the asserted global vanishing before any collapse assertion.

\paragraph{The local reflection ring.}
Reflection sends a unit tangent vector to its negative. Its action on
$C_r=SM$ is free, and $H_R^*(C_r)=H^*(P(TM))$ for every $r$.
Here $P(TM)\to M$ is the bundle of unoriented real tangent lines, and
$w$ is the first Stiefel--Whitney class of its tautological line.
For $\mathbb C P^m$ and $\mathbb H P^m$ the tangent Stiefel--Whitney
class is $(1+h)^{m+1}$, giving precisely $\mathcal R_m$ by the real
projective-bundle formula. For the Cayley plane it is $1+h+h^2$, so
\[
 H^*(P(TM))=\mathbb F_2[w,h]/(h^3,w^{16}+hw^8+h^2)
            =\mathcal R.
\]
The Cayley tangent formula also follows from the Wu formula: the
positive Wu classes are $v_8=h$, $v_{16}=0$, and
$\operatorname{Sq}^8h=h^2$.
Each local Thom module is cyclic over the ambient ring
$\mathbb F_2[w,h]$. Its bottom generator is the Thom image of
$1\in H_R^0(C_r)$ and has degree $\lambda_r$,
and $\lambda_r+1<\lambda_{r+1}$. Consequently $\delta_r$ vanishes on
that generator. Ambient-linearity makes it zero on the whole module.
The same conclusion in degree $j-1$ also makes the first map in
\eqref{eq:supported-triple} injective. Thus that sequence is short
exact for $H=R$ in every degree.
In particular all supported inclusions into $B$ are injective.
The local bottom has a unique supported lift since
$\lambda_r<\lambda_{r+1}$. Let $X_r$ denote this lift and its image
in $B$.

\paragraph{Odd local full-group rings.}
A primitive geodesic in a complex or quaternionic projective space lies
in a unique projective line. That line is a round $S^d$, and its great
circles correspond to the unoriented real two-planes in the real bundle
$V=\operatorname{Herm}_0(\mathcal E)$ over
$\operatorname{Gr}_2(\mathbb K^{m+1})$, where $\mathcal E$ is the
tautological $\mathbb K$-plane and $\operatorname{Herm}_0$ denotes
the trace-free Hermitian endomorphisms. This Grassmannian parametrizes
projective lines in $M$; the bundle $\operatorname{Gr}_2(V)$ over it
parametrizes the unoriented geodesics within those lines.

Write $a,z$ for the first and second characteristic classes of
$\mathcal E$, using mod-two Chern classes when $\mathbb K=\mathbb C$
and mod-two symplectic classes when $\mathbb K=\mathbb H$. Their
degrees are $d,2d$. The real rank of $V$ is $d+1$, and
\[
 w(V)=1+a,\quad |a|=d,\qquad
 H^*(\operatorname{Gr}_2(\mathbb K^{m+1}))=
 \mathbb F_2[a,z]/(Q_m(a,z),Q_{m+1}(a,z)).
\]
The formula for $w(V)$ follows by splitting $\mathcal E$ into two
$\mathbb K$-lines. The total characteristic class of its orthogonal
complement is $(1+a+z)^{-1}$, whose components are the polynomials
$Q_j(a,z)$. The complement has $\mathbb K$-rank $m-1$, giving the
relations $Q_m=Q_{m+1}=0$; the Schubert basis proves that these are
all relations. The same projective-bundle and Leray--Hirsch argument
applies to the real Grassmann bundle.

The standard $G$-plane becomes the tautological real two-plane on
$\operatorname{Gr}_2(V)$. In this identification the rotation
parameter is the great-circle angle itself. If $\mathcal Q$ denotes
its real rank-$(d-1)$ complement, then
\[
 (1+w+U)w(\mathcal Q)=1+a.
\]
Equating the homogeneous terms above the rank of $\mathcal Q$ first
expresses $a$ as $A_d$ and then gives $w^{d-1}U=0$. In the two cases
this yields
\begin{align*}
 A_2&=w^2+U,& A_4&=w^4+w^2U+U^2,\\
 \mathcal O_m&=\mathbb F_2[w,U,z]/
 \bigl(w^{d-1}U,Q_m(A_d,z),Q_{m+1}(A_d,z)\bigr).
\end{align*}
For odd $r$ the ineffective isotropy is cyclic of order $r$, hence has
mod-two acyclic classifying space. The weight-$r$ representation
preserves the mod-two classes $w,U$ because $r$ is odd.
Thus $H_G^*(C_r)=\mathcal O_m$ for every odd $r$.
It is generated over $A$ by
\begin{equation}
\label{eq:odd-generators}
 z^j,\qquad 0\leq j\leq\left\lfloor\frac{m-1}{2}\right\rfloor.
\end{equation}
Indeed $Q_m$ for even $m$, and $Q_{m+1}$ for odd $m$, is monic in
$z^{\lfloor(m+1)/2\rfloor}$. Restriction to reflection is
\begin{equation}
\label{eq:restriction-z}
 w\longmapsto w,\qquad U\longmapsto0,\qquad
 z\longmapsto h(h+w^d).
\end{equation}
Geometrically, let $L$ be the tautological $\mathbb K$-line pulled
back from $M$ to $P(TM)$, and let $\ell$ be the tautological real
tangent-direction line. The projective line determined by a point
and a tangent direction pulls $\mathcal E$ back to
$L\oplus(L\otimes_{\mathbb R}\ell)$. The two characteristic roots
are $h$ and $h+w^d$, so their product is exactly the last formula in
\eqref{eq:restriction-z}.

For the Cayley plane the same great-circle construction gives
$\operatorname{Gr}_2(V)$ over the dual Cayley plane, with $V$ the
standard rank-nine $\operatorname{Spin}(9)$ bundle. Its characteristic
class is $w(V)=1+x$, $|x|=8$.
The geometric line and representation identifications are in
\cite[Section~3, Lemmas~1--2]{Wolf1963}. For the characteristic class,
let $p:F_4/\operatorname{Spin}(8)\to F_4/\operatorname{Spin}(9)$
be the incidence projection to the dual Cayley plane. With the
integral degree-eight generators $y,y'$ of
\cite[(2.1)--(2.3), Theorem~2.3, and (2.43)--(2.46)]{Yokotani1992},
the pullback of an integral lift of $x$ is $y$
and the Euler class of the standard eight-plane is $y+2y'$.
The pullback $p^*$ is injective in degree eight. Restricting $V$ as
the trivial line plus that eight-plane and reducing modulo two thus gives $w_8(V)=x$;
the other positive classes vanish by rank and degree.
The complement relations give
\begin{equation}
\label{eq:cayley-odd}
 H_G^*(C_r)=\mathbb F_2[w,U]/(w^7U,A_8^3),\qquad
 A_8=w^8+w^6U+w^4U^2+U^4,\qquad r\text{ odd}.
\end{equation}
Under restriction $x=A_8$ becomes $w^8$, not the evaluation class $h$.

\paragraph{Even local full-group modules.}
Let $r$ be even and write $\mathcal E_r=H_G^*(C_r)$.
The subgroup $K$ generated by reflection and half-translation is
$C_2\times C_2$, and half-translation acts trivially on $C_r$.
Use the coordinates $w,v$ on $H^*(BK)$ introduced in
Section~\ref{sec:sphere-topology}. The triviality of half-translation
and the freeness of reflection give
\[
 H_K^*(C_r)=\mathcal R[v],\qquad U\longmapsto v(v+w).
\]
Multiplication by this monic polynomial in $v$ is injective, including
when the coefficient ring has zero divisors. The weight-two real
$G$-plane has vanishing mod-two Euler class: restriction to $BSO(2)$
kills its Euler class modulo two, and restriction to $BR$ has a trivial
line summand. Its unit sphere has stabilizer $K$, so its Gysin sequence
injects $\mathcal E_r$ into $H_K^*(C_r)$. It follows that $U$ is regular
on $\mathcal E_r$. The standard-plane Gysin sequence then gives
\begin{equation}
\label{eq:even-U}
 \mathcal E_r/(U)\cong\mathcal R
 \quad\text{by restriction to }R.
\end{equation}
Choose homogeneous lifts of $1,h,\ldots,h^m$, or $1,h,h^2$ for
$\mathbb O P^2$, through \eqref{eq:even-U}. These generate
$\mathcal E_r$ over $A$: express the reduction of a class as an
$\mathbb F_2[w]$-linear combination of these generators, subtract the
chosen lifts, and write the remainder as $U$ times a class of degree
two lower. Repetition terminates because cohomology has no negative
degrees. This argument requires only a quotient presentation and
regularity of $U$, rather than a choice of an $A$-module splitting.

\paragraph{Collapse for the full group, including the complex boundary.}
The connecting map $\delta_r$ is $A$-linear and raises total degree
by one. A generator of offset $q$ has zero image whenever
$q+1<\nu$, since its target then lies below $\lambda_{r+1}$.
For the quaternionic case all odd generating offsets in
\eqref{eq:odd-generators} are at most $4m-4$ and all even offsets are at
most $4m=\nu-2$. For the Cayley plane the odd offset is zero and the
even offsets are at most $16<\nu-1$. Degree and $A$-linearity therefore
make every full-group $\delta_r$ vanish in these cases.
For the complex case the odd offsets are at most $2m-2=\nu-2$.
The even generators of offsets $2k$, $k<m$, are equally harmless.
Only the top even generator, of offset $2m=\nu$, requires attention.
Let $y$ be its connecting image, of degree $\lambda_{r+1}+1$.
Restriction of $y$ to the next, odd local module has offset one, so
it is a multiple of the class $w$ times the bottom Thom class.
The local class $w$ restricts to the nonzero class $w$ in
$\mathcal R_m$. On the other hand, the reflection restriction of
$y$ is zero, by naturality and the reflection calculation. Thus the
coefficient of this local class is zero. Finally
$\lambda_{r+1}+1<\lambda_{r+2}$, so the map from the supported
group containing $y$ to its local term is injective. Hence $y=0$.
This uses only the degree range already established before the
full-group connecting maps were considered.
We have proved in all cases that
\begin{equation}
\label{eq:supported-exact}
 0\longrightarrow F_{r+1}^H\longrightarrow F_r^H
 \longrightarrow L_r^H\longrightarrow0,
 \qquad H=G,R,
\end{equation}
in every degree. These all-degree injections will also transport
annihilator relations to later supported groups. In particular every
local class now has a supported lift; it is unique whenever its
offset is less than $\nu$.

\paragraph{Global reflection relations.}

The vanishing of connecting maps supplies global generators, but it
does not yet establish their local relations globally. A locally zero
product could become a class supported at a later critical level.
We now rule out this possibility. Most products lie below the next
index; the complex case has one boundary relation to check separately.
The ambient relation $h^{m+1}=0$, or $h^3=0$ for the Cayley plane,
already holds on the whole loop space. The other local reflection
relation has offset $4m=\nu-2$ for $\mathbb H P^m$ and offset
$16<\nu$ for the Cayley plane, so it holds on $X_r$ by injectivity
below the next index.
For $\mathbb C P^m$ its offset is exactly $\nu$. Thus
\eqref{eq:supported-exact} permits only
\[
 P_m(w,h)X_r=cX_{r+1},\qquad c\in\mathbb F_2.
\]
Indeed, the left side has zero local term at stage $r$. Its unique
antecedent at stage $r+1$ has offset zero, and restriction identifies
that supported group with its one-dimensional local bottom.
Forget equivariance and write $\theta_r$ for the ordinary image of
$X_r$. The symmetric energy is perfect in ordinary mod-two cohomology,
and each $\theta_r$ is nonzero; this is the cohomological perfection and
bottom-class construction of \cite[Theorem~14.2]{GH2009}.
For odd $m$, $P_m(0,h)=(m+1)h^m=0$. The displayed relation therefore
becomes $0=c\theta_{r+1}$, and $c=0$ because $\theta_{r+1}\ne0$.
For even $m$ write $\mathbb H_k(\Lambda M)=H_{k+2m}(\Lambda M)$
for the shifted grading and use the integral complex-projective
loop-ring computation \cite[Theorem~3]{CJY2004}.
The torsion in that integral ring has order dividing $m+1$, which is
odd in the present case. Thus coefficient
change to $\mathbb F_2$ introduces no Tor term and gives
\[
 \mathbb H_*(\Lambda\mathbb C P^m;\mathbb F_2)=
 \mathbb F_2[a,x,y]/(a^{m+1},x^2,xa^m,ya^m),
 \quad |a|=-2,\quad |x|=-1,\quad |y|=2m.
\]
Multiplication by $a$ is cap product with the evaluation class $h$.
Every odd-degree homology class contains $x$, so the relation
$xa^m=0$ says that cap product with $h^m$ kills every such class.
For any odd-degree cohomology class $\alpha$ and any homology class
$c$ of degree $|\alpha|+2m$, the Kronecker pairing gives
\[
 \langle h^m\alpha,c\rangle=\langle\alpha,h^m\mathbin{\cap}c\rangle=0.
\]
Over a field this proves $h^m\alpha=0$. The constant-loop inclusion
has evaluation as a retraction, so its cohomology restriction is
surjective; since $\mathbb C P^m$ has no odd cohomology, the relative
and absolute odd-degree groups are naturally isomorphic. Thus the
same vanishing applies relatively. Since
$|\theta_r|=1+2m(r-1)$ is odd, $h^m\theta_r=0$; forgetting
equivariance in the proposed relation again forces $c=0$.

All defining relations of $\mathcal R X_r$ now hold globally. Local
restriction identifies that block faithfully with its local module.
For independence, restrict a relation to its earliest nonzero stage:
the coefficient of $X_r$ must be zero in the corresponding local
module, hence in the displayed global block. Repeating eliminates
all terms of the relation. For spanning, subtract the block with the
same first local term as a given homogeneous class, and apply
\eqref{eq:supported-exact} to continue at the next stage. This
procedure terminates because only finitely many $\lambda_r$ lie
below a fixed degree. We obtain
\begin{equation}
\label{eq:actual-reflection-splitting}
 B=\bigoplus_{r\geq1}\mathcal R X_r.
\end{equation}
Thus the local relations hold in the module itself, and there are no
extensions between the displayed blocks. Collapse alone would have
identified only the associated graded module.
Since $P(TM)$ has dimension $2n-1$, $w^{2n}\mathcal R=0$;
\eqref{eq:actual-reflection-splitting} proves the stated global bound.

\paragraph{The paired algebra.}

We next explain the two summands in each odd block. In the complex
and quaternionic cases the useful variable is the quadratic expression
$h(h+w^d)$. With this choice, the reflection ring is a free extension
of the odd restriction ring with basis $1,h$. Thus pairing by $h$ is
an explicit algebraic feature of the computed ring.
Put $a=w^d$ and $z=h(h+a)$. Directly from the recurrence,
\begin{align*}
 Q_m(a,h(h+a))&=P_m(w,h),\\
 Q_{m+1}(a,h(h+a))+hQ_m(a,h(h+a))&=(h+a)^{m+1},\\
 aQ_m(a,h(h+a))&=(h+a)^{m+1}+h^{m+1}.
\end{align*}
For a direct verification put $b=h+a$. The recurrence gives
\[
 Q_m(a,hb)=\sum_{i=0}^{m}h^i b^{m-i}.
\]
Multiplication by $h+b=a$ telescopes to $h^{m+1}+b^{m+1}$,
and adding $hQ_m$ to $Q_{m+1}$ leaves $b^{m+1}$. Expanding
the sum in $a,h$ gives coefficient $\binom{m+1}{i}$ on
$h^i a^{m-i}$ by the hockey-stick identity, proving the first
identity as well. These arguments take place in the polynomial
ring over $\mathbb F_2$ and require no division by $a$.
To see that they identify the defining ideals, impose
$z=h(h+a)$. The relations $P_m=0$ and $h^{m+1}=0$ imply
$Q_m=Q_{m+1}=0$ by the second and third identities. Conversely
$Q_m=Q_{m+1}=0$ imply $(h+a)^{m+1}=0$ by the second identity,
and then $h^{m+1}=0$ by the third. Therefore
\[
 \mathcal R_m\cong
 \mathcal S_m[h]/(h^2+w^dh+z),
\]
a free rank-two extension with basis $1,h$. In particular the displayed
copy of $\mathcal S_m$ is embedded, and is the image of the odd local
restriction. The Cayley presentation makes the same rank-two assertion
immediate. This proves the abstract decomposition $B=W\oplus hW\oplus E$.

\paragraph{Lifting annihilator relations.}

It remains to control powers of $U$ on full-group lifts. A local
annihilator relation may acquire a later-supported error, just as in
the reflection calculation. The following observation lets us detect
that error at the next stage and remove it using an injective
multiplication there.
We will use the following consequence of \eqref{eq:supported-exact}.
Suppose $\alpha\in F_r^{G,j}$ has zero local term, so there is a
unique $\beta\in F_{r+1}^{G,j}$ mapping to $\alpha$. If a polynomial
$b\in A$ satisfies $b\alpha=0$, injectivity of
$F_{r+1}^G\to F_r^G$ gives $b\beta=0$ in the full product degree.
Let $\overline\beta\in L_{r+1}^{G,j}$ be the next local term. Then
$b\overline\beta=0$. Consequently, if multiplication by $b$ is
injective on this local degree, $\overline\beta=0$. If also
$j<\lambda_{r+2}$, the initial degree range makes
$F_{r+1}^{G,j}\to L_{r+1}^{G,j}$ injective and proves
$\beta=\alpha=0$. The degree restriction concerns $j$, the degree
of the possible extension; it does not concern $j+|b|$.

\paragraph{Complex odd lifts and their uniform annihilator.}
For odd $r$, let $Z_{r,j}$ be the supported lift of $z^j$ in
\eqref{eq:odd-generators}. Its offset is at most $\nu-2$, so it is
unique and its reflection restriction is exactly
$[h(h+w^2)]^jX_r$.
The locally zero product $wUZ_{r,j}$ has offset $4j+3$.
If this offset is below $\nu$, injectivity of supported restriction
makes the product zero globally.
The sole exception is odd $m$, $j=(m-1)/2$, where the offset is
$\nu+1$. Its later-supported antecedent has next-even local degree one.
Reflection kills the original product, and reflection supported
inclusions are injective; hence it kills that antecedent as well.
By \eqref{eq:even-U} the kernel of local restriction is $U\mathcal E_{r+1}$,
which is zero in degree one. The local antecedent is zero, and its
original degree is below $\lambda_{r+2}$. Therefore
\begin{equation}
\label{eq:CP-wU}
 wUZ_{r,j}=0
\end{equation}
also in this boundary case, including $m=1$.

Set $N_j=2m-1-2j>0$. The Grassmann base has real dimension $4m-4$,
whereas $a^{N_j}z^j$ has degree $4m-2$. Consequently
$A_2^{N_j}z^j=0$ locally. On $Z_{r,j}$ the possible extension into
the next even local module has degree $q=(4m-2)-\nu=2m-2$, and the original offset
$4m-2$ is below $2\nu=4m$.
Both defining relations of $\mathcal R_m$ belong to
$\mathbb F_2[w^2,h]$. Writing $\xi=w^2$, one therefore has
\[
 \mathcal R_m=D\oplus wD,\qquad
 D=\mathbb F_2[\xi,h]/(h^{m+1},\widetilde P_m(\xi,h)),
\]
where $D$ is concentrated in even degrees and $\widetilde P_m$ is
obtained from $P_m$ by replacing $w^{2k}$ by $\xi^k$.
Multiplication by $w$ is injective
on every even-degree part of $\mathcal R_m$.
Regularity of $U$ and induction after reducing modulo $U$ give the
same even-degree injectivity on $\mathcal E_{r+1}$: if $wx=0$, its
reduction is zero, so $x=Uy$; cancel $U$ and continue in lower degree.
At the final step the degree is negative and the class is zero.
Thus $wU$ is injective on degree $q$. By \eqref{eq:CP-wU} and the
lifting observation, the possible extension vanishes. We have proved
$A_2^{N_j}Z_{r,j}=0$ globally.
Modulo $wU$, every mixed monomial in $(w^2+U)^{N_j}$ vanishes.
Thus $A_2^{N_j}=w^{2N_j}+U^{N_j}$ modulo that ideal, and multiplying
the resulting relation on $Z_{r,j}$ by $U$ gives
\[
 U^{N_j+1}Z_{r,j}=0,
 \qquad\text{and hence}\qquad U^{2m}Z_{r,j}=0.
\]

\paragraph{Quaternionic odd lifts and their uniform annihilator.}
For odd $r$ and $0\leq j\leq\lfloor(m-1)/2\rfloor$, choose the
unique supported class $Z_{r,j}$ lifting $z^j$.
Its restriction is $[h(h+w^4)]^jX_r$.
The product $w^3UZ_{r,j}$ has offset $8j+5\leq4m+1<\nu$ and is zero.
Set $N_j=2m-1-2j$, which is positive and odd.
The Grassmann base has dimension $8m-8$, while $a^{N_j}z^j$ has
degree $8m-4$. Thus $A_4^{N_j}Z_{r,j}$ is locally zero.
For $m=1$ its offset four is already below $\nu=6$.
For $m\geq2$ its possible next-even obstruction has degree
$q=4m-6$, and its original offset $8m-4$ is below $2\nu$.
The first relation of $\mathcal R_m$ has degree $4m$.
Multiplication by $w^3$ is therefore injective on its degrees at most
$4m-4$: both the source and the target still lie below the first
defining relation. Reduction modulo regular $U$ lifts this injectivity to
$\mathcal E_{r+1}$ in the same range. In particular $w^3U$ is injective
in degree $q$. The lifting observation kills the obstruction.
For odd $N$, reduction modulo $w^3U$ yields
\[
 A_4^N=w^{4N}+w^2U^{2N-1}+U^{2N}.
\]
Apply this with $N=N_j$ to $Z_{r,j}$. Multiplication by $U$ eliminates
the term $w^{4N_j}U Z_{r,j}$ because $w^3U Z_{r,j}=0$, leaving
\[
 U^{2N_j+1}Z_{r,j}=w^2U^{2N_j}Z_{r,j}.
\]
Multiplying by $w^2$ gives $w^2U^{2N_j+1}Z_{r,j}=0$.
Multiplying the preceding equality by $U$ therefore gives
\[
 U^{2N_j+2}Z_{r,j}=0,
 \qquad\text{so uniformly}\qquad U^{4m}Z_{r,j}=0.
\]

\paragraph{Cayley odd lifts and their uniform annihilator.}
Fix odd $r$, and let $\zeta_r$ be the unique supported lift of the
local bottom. Its product with $w^7U$ has offset nine, strictly less
than $\nu=22$, and is zero. The locally zero $A_8^3\zeta_r$ has
offset $24<2\nu$. Its only possible extension is therefore in the
next even local module, in degree $24-\nu=2$.
Multiplication by $U$ takes this local degree to four, and
multiplication by $w^7$ is injective there modulo $U$, because its
target degree eleven is below the first reflection relation in
degree sixteen. Regularity of $U$ and the preceding reduction
argument give injectivity of $w^7U$ on degree two upstairs. Since
$w^7U\zeta_r=0$, the lifting observation proves
$A_8^3\zeta_r=0$ globally.
Modulo $w^7U$,
\[
 A_8^3=w^{24}+w^6U^9+w^4U^{10}+U^{12}.
\]
On multiplying the relation on $\zeta_r$ by $U$, the first term
vanishes and we obtain
\[
 U^{13}\zeta_r=w^6U^{10}\zeta_r+w^4U^{11}\zeta_r.
\]
Multiplication by $w^4$ and $w^6$ respectively gives
$w^4U^{13}\zeta_r=w^6U^{13}\zeta_r=0$, since every resulting term
is divisible by $w^7U\zeta_r$. Multiplying the displayed equality
by $U^3$ now gives $U^{16}\zeta_r=0$.

\paragraph{Exact even lifts and the entire restriction image.}
Choose supported full-group lifts of the even local generators
$1,h,\ldots,h^m$ (or $1,h,h^2$ for the Cayley plane).
All restrict exactly to $h^kX_r$ by the low-degree range, except
possibly the complex generator $h^m$ at offset $\nu$.
Its restriction can differ from $h^mX_r$ only by $cX_{r+1}$.
Since $r+1$ is odd, the full-group lift $Z_{r+1,0}$ of $X_{r+1}$
exists. Adding $cZ_{r+1,0}$ normalizes the even lift to restrict
exactly to $h^mX_r$. Its local leading term is unchanged.

The chosen odd and normalized even lifts generate the whole global
full-group module over $A$. Indeed, subtract a lift of the first
nonzero local term using \eqref{eq:supported-exact}, then repeat at
the next stage. For each homogeneous total degree this procedure
terminates because $\lambda_r\to\infty$.
Their reflection images generate exactly $W\oplus E$, since all
positive powers of $U$ restrict to zero and the odd $z^j$ generate
$\mathcal S$ over $\mathbb F_2[w]$.
This proves the asserted equality of restriction images, not just
the existence of selected lifts.

\paragraph{The entire $U$-torsion subgroup.}
Every chosen odd lift is killed by the stated common $U^N$.
Conversely let $\xi$ be a homogeneous $U$-torsion class and examine
its first nonzero local term. That term cannot be at an even stage,
where $U$ is regular. At an odd stage subtract an $A$-combination of
the chosen odd lifts with the same local term. The difference remains
$U$-torsion and is supported later. Repeating in fixed total degree
terminates. Thus the odd lifts generate precisely $T_U$ over $A$,
and $U^NT_U=0$. Their exact reflection restrictions show
$\operatorname{res}(T_U)=W$. This completes the proof.
\end{proof}

\paragraph{Low-dimensional examples.}
The general presentations become small rings in the first few cases.
These examples also indicate where the boundary-degree arguments in
the proof enter. For $\mathbb C P^1$, the odd ring is
$\mathbb F_2[w,U]/(wU,w^2+U)$ and the reflection ring is
$\mathbb F_2[w,h]/(w^2,h^2)$; the degree-one restriction used above
is nonzero and the conclusion gives $U^2T_U=0$.
For $\mathbb C P^2$ these rings are respectively
$\mathbb F_2[w,U]/(wU,(w^2+U)^3)$ and
$\mathbb F_2[w,h]/(w^6,h^2+w^2h+w^4)$, with $U^4T_U=0$.
For $\mathbb C P^3$ the odd generators have offsets zero and four.
The top $wU$ relation has offset seven, exceeding the index gap six;
its vanishing is exactly the degree-one reflection-kernel argument in
the proof.
For $\mathbb H P^2$ the rings are
$\mathbb F_2[w,U]/(w^3U,A_4^3)$ and
$\mathbb F_2[w,h]/(w^{12},h^2+w^4h+w^8)$, with $U^8T_U=0$.
For $\mathbb H P^3$ the relation $A_4^3$ on both odd generators
sharpens the general bound $U^{12}$ to $U^8$. For the Cayley plane the odd ring has total
rank 108 and top degree 30, agreeing with the Grassmann-bundle model.

All sublevels in this calculation belong to the symmetric metric.
Their role is to determine the Borel modules and restriction maps,
which depend only on $M$ and the reparametrization action. The sublevel estimates for the arbitrary metric are supplied separately
by the analytic part of the argument.


\section{The finite-prime contradiction}
\label{sec:fp-contradiction}
\label{sec:finite-prime}

We now explain why the analytic and topological results force new
geodesic images. Suppose, to the contrary, that only finitely many
prime images occur. Each critical level then consists of finitely many
iterated phase circles. A class of degree at least two vanishes on a
phase circle, so the evaluation class $h$ kills the local contribution
at every crossing. On the other hand, the global calculation pairs
classes by multiplication by $h$. The proof must bring these local
and global statements to the same finite sublevel.

The quantity which makes this possible is the number of reflection
classes at a finite cut which do not extend from the full loop space.
We call these classes transient. We prove that their total dimension
is bounded and cannot decrease when the length cut is doubled. At a
cut where equality holds, multiplication by $h$ must be injective on
the visible odd summand. A global odd class at its first nonzero
restriction is killed by $h$, contradicting this injectivity.

We state the argument as a criterion. This makes clear which parts
use the metric-dependent filtration estimates and which use the
global equivariant calculation. The proof includes the local estimates
for degenerate critical orbits and allows the global summands to
occupy overlapping degrees.

Let $(M,g)$ be a closed smooth Riemannian manifold of dimension $n\geq2$.
Write $\Lambda=H^1(S^1,M)$, let $M\subset\Lambda$ denote the constant
loops, and put
\[
 E(\gamma)=\int_0^1|\dot\gamma(t)|_g^2\,dt,
 \qquad F=\sqrt E,\qquad X_a=\{F<a\}\quad(a>0).
\]
All homology and cohomology groups in this section have coefficients in
$\F$. A cut $a$ is called regular when $a^2$ is a regular energy value.
Let $G=O(2)$ act by reparametrization, let $\rho(t)=-t$ and
$\tau(t)=t+1/2$, and set
\[
 R=\langle\rho\rangle,\qquad K=\langle\rho,\tau\rangle\cong C_2\times C_2.
\]
For a group $H$ and an $H$-pair $(Y,Z)$, put
$Y_H=EH\times_HY$ and $Z_H=EH\times_HZ$, where $EH$ is a
contractible free $H$-space. The notation $H_H^*(Y,Z)$ means
$H^*(Y_H,Z_H;\F)$. We use the direct sum,
rather than the product, of its homogeneous groups. Set
\[
 B_a^j=H_R^j(X_a,M),\qquad B_\infty^j=H_R^j(\Lambda,M),
 \qquad B_a=\bigoplus_j B_a^j.
\]
The class $w\in H^1(BG)$ is the first Stiefel--Whitney class of the
standard plane and $U\in H^2(BG)$ its second class. The same letter $w$
denotes its restrictions. Let $v\in H^1(BK)$ be the character nonzero
on $\tau$ and zero on $\rho$. Thus
\begin{equation}\label{eq:fp-coefficient-rings}
 H^*(BG)=\F[w,U],\qquad H^*(BK)=\F[v,w],
 \qquad U|_K=v(v+w).
\end{equation}
If $\eta\in H^d(M)$ with $d\geq2$, evaluation at zero defines the
ambient class $h=\operatorname{ev}_0^*\eta\in H_R^d(\Lambda)$.
Indeed evaluation is $R$-equivariant when $R$ acts trivially on $M$.
Its cup action on any relative group is the action of this absolute
class after restriction to the relevant Borel pair.

\begin{theorem}[A finite-prime obstruction]\label{thm:fp-obstruction}
Assume the following properties.
\begin{enumerate}
\item There is $D>0$, independent of $a$ and $q$, such that
\begin{equation}\label{eq:fp-UD}
 \ker\bigl[H_q(X_a,M)\longrightarrow H_q(\Lambda,M)\bigr]
 \subset
 \ker\bigl[H_q(X_a,M)\longrightarrow H_q(X_{a+D},M)\bigr].
\end{equation}
\item There are $\lambda>0$ and $C_{\rm sp}\geq0$ such that, for every
$q\geq0$, the map $H_q(X_a,M)\to H_q(\Lambda,M)$ is zero if
$a<\lambda q-C_{\rm sp}$ and is surjective if
$a>\lambda q+C_{\rm sp}$. Each $H_q(\Lambda,M)$ is finite dimensional.
\item There is an integer $N_w\geq1$ such that $w^{N_w}B_\infty=0$.
\item There are graded subspaces $W,E_0\subset B_\infty$ and an ambient
evaluation class $h$ of degree $d\geq2$ such that
\begin{equation}\label{eq:fp-global-decomposition}
 B_\infty=W\oplus hW\oplus E_0.
\end{equation}
Every homogeneous class of $W$ is the restriction of a class in
$H_G^*(\Lambda,M)$ killed by some power of $U$, and every homogeneous
class of $E_0$ is the restriction of a class in $H_G^*(\Lambda,M)$.
\item There are nonzero homogeneous classes $Z_s\in W$ whose degrees
$q_s$ tend to infinity and have bounded gaps: for some $b<\infty$,
$q_{s+1}>q_s$ and $q_{s+1}-q_s\leq b$ for every $s$.
\end{enumerate}
Then $(M,g)$ has infinitely many geometrically distinct prime closed
geodesics.
\end{theorem}

The first two hypotheses control when ordinary classes appear and how
quickly null classes disappear. The next two provide the global
reflection module and its full-group lifts. The last hypothesis
ensures that, in a sufficiently long interval of degrees, there is an
odd class available for the first-appearance argument. The spectral
condition applies to all classes in a degree; it does not require that
the degree have rank one. Likewise the direct sum in
\eqref{eq:fp-global-decomposition} is a statement about subspaces,
even when their degree supports overlap.

We prove the criterion in the following order. Local iteration theory
first gives uniform crossing bounds. Restriction to half-period-fixed
loops then yields a linear bound on a finite-cut defect. Combining it
with the filtration estimates bounds the total number of transients.
The final doubling comparison produces the contradiction. Throughout,
mean indices are allowed to vanish, and only degreewise finiteness of
global homology is assumed.

\subsection{Gysin sequences and local geometry}

The group inclusions $R\subset K\subset G$ give the comparison maps
used throughout the proof. We first record their exact sequences and
then construct local crossing models on which these maps are induced
by the same inclusions of spaces. This will allow the subsequent rank
counts to retain both the ordinary grading and the subgroup restriction.

The sphere of the standard $G$-plane is $G/R$. The plane obtained by
doubling rotation angles has sphere $G/K$ and zero mod-two Euler class.
To see the last assertion, its restriction to $K$ has a fixed line,
whereas the map $H^*(BG)\to H^*(BK)$ in
\eqref{eq:fp-coefficient-rings} is injective. One can also see this
injection from the projective bundle of the standard plane, whose total
space is $BK$: its two line characters are $v$ and $v+w$.
The relative Thom sequences consequently give, for every $G$-pair,
\begin{align}
 \cdots\longrightarrow H_G^{j-2}(Y,Z)
 &\xrightarrow{\ U\ }H_G^j(Y,Z)
 \xrightarrow{\operatorname{res}_R}H_R^j(Y,Z)
 \xrightarrow{\pi}H_G^{j-1}(Y,Z)
 \xrightarrow{\ U\ }H_G^{j+1}(Y,Z)\longrightarrow\cdots,
 \label{eq:fp-standard-gysin}\\
 0\longrightarrow H_G^j(Y,Z)
 &\xrightarrow{\operatorname{res}_K}H_K^j(Y,Z)
 \xrightarrow{\pi}H_G^{j-1}(Y,Z)\longrightarrow0.
 \label{eq:fp-weight-two-gysin}
\end{align}
The sign lines give
\begin{align}
 \cdots\longrightarrow H_R^{j-1}(Y,Z)
 &\xrightarrow{\ w\ }H_R^j(Y,Z)\longrightarrow H^j(Y,Z)
 \longrightarrow H_R^j(Y,Z)
 \xrightarrow{\ w\ }H_R^{j+1}(Y,Z)\longrightarrow\cdots,
 \label{eq:fp-sign-R}\\
 \cdots\longrightarrow H_K^{j-1}(Y,Z)
 &\xrightarrow{\ v\ }H_K^j(Y,Z)\longrightarrow H_R^j(Y,Z)
 \longrightarrow H_K^j(Y,Z)
 \xrightarrow{\ v\ }H_K^{j+1}(Y,Z)\longrightarrow\cdots.
 \label{eq:fp-sign-K}
\end{align}
The first of these applies to an $R$-pair and the second to a $K$-pair.
All maps are natural for maps of pairs and for compatible subgroup
restrictions, and all connecting maps satisfy the projection formula
for ambient absolute classes.

These sequences apply to singular Borel pairs without a cofibration
assumption. For clarity, if $D(\xi),S(\xi)$ are the disk and sphere
bundles of a real bundle over $Y_H$, restrict them over $Z_H$ and use
the chain subspace
\[
 C_*(S(\xi))+C_*(D(\xi)|_{Z_H})\subset C_*(D(\xi)).
\]
The two summands intersect in $C_*(S(\xi)|_{Z_H})$. Restriction of
cochains is surjective onto a chain subspace, since coefficients form
a field. The two resulting short exact cochain sequences, together
with cup product by a Thom cocycle, give the relative Thom sequence
and the displayed naturalities. The absolute Thom isomorphism over a
general base follows by pullback to its singular CW replacement and
the disk and sphere fibrations.

Assume from now on that there are finitely many prime images, and
choose a unit-speed prime representative $c_i$ of each image,
$1\leq i\leq p$. Write $\ell_i>0$ for its prime length. A
nonconstant closed geodesic has a least
positive period, so it is an iterate of a prime. Two unit-speed primes
with the same image coincide up to phase and reversal. Indeed finitely
many embedded arcs cover the first parameter circle; Baire's theorem
provides an interval of the second mapping into one such arc. Away
from the arc endpoints their tangent lines agree, and uniqueness for
the geodesic equation proves the assertion. A phase circle cannot
coincide with its reverse: an equality $\gamma(t)=\gamma(a-t)$ would
give $\dot\gamma=0$ at a fixed parameter of this reflection.
Thus each unoriented label $(i,m)$ consists of two disjoint oriented
phase circles, and the positive critical $F$-values are $m\ell_i$.
Here a phase circle is the orbit of a parametrized closed geodesic
under translations of its parameter. We call a positive critical
value an \emph{event}; several labels may occur at the same event.
In particular,
\begin{equation}\label{eq:fp-event-count}
 \#\{\hbox{labels below }a\}=O(1+a),\qquad
 \#\{\hbox{labels in }[a,a+b]\}=O(1+b).
\end{equation}
The constants depend on the fixed prime list. More explicitly, there
are at most $\sum_i\lfloor a/\ell_i\rfloor$ labels below $a$, and
at most $\sum_i(2+b/\ell_i)$ labels in an interval of width $b$.

We record the analytic models used below. A smooth embedding of $M$
into Euclidean space identifies $\Lambda$ with a closed Hilbert
submanifold of the ambient $H^1$ space. Its induced strong metric is
complete and invariant under the rotations and reflections in $O(2)$.
An arbitrary change of speed need not preserve this metric or the
energy. The energy is smooth
and satisfies Palais--Smale. In fact an energy-bounded sequence has a
weakly $H^1$ and uniformly convergent subsequence. Testing the energy
derivative on the tangent projection of its difference from the limit
shows convergence of its velocity norms; uniform positivity of the
metric then gives strong $H^1$ convergence. The weak Euler--Lagrange
equation gives precisely the smooth periodic geodesics and constants.
The bounded vector field
\[
 -\nabla E\big/(1+\|\nabla E\|^2)^{1/2}
\]
has a complete invariant flow, strictly decreasing energy away from
the critical set. Palais--Smale bounds its energy drop away from zero
on a critical-point-free bounded band. Small sublevels retract
equivariantly to $M$, by the positive normal Hessian along the
constants and the mean-zero Poincar\'e inequality. It follows that
regular open and closed sublevels are naturally equivalent, and that
a strip without critical values induces equivalences compatible with
fixed loci.

For a fixed energy bound $a^2$ and a sufficiently large even $N$, use
the polygon manifold
\[
 P_N=\{(x_0,\ldots,x_{N-1})\in M^N:
                 d_g(x_j,x_{j+1})<\rho_0\},\qquad
 f_N=N\sum_{j=0}^{N-1}d_g(x_j,x_{j+1})^2,
\]
where $\rho_0>0$ is smaller than the injectivity radius and indices
are cyclic. Taking $a/\sqrt N<\rho_0$ makes the closed sublevel compact.
Replace the centered part of each time cell by its minimizing geodesic
and let this part grow to the entire cell. Cauchy--Schwarz proves that
energy does not increase. Strong $H^1$ continuity at zero replacement
width follows from absolute continuity of the velocity integrals;
at positive width it follows from smooth dependence of short
geodesics on their endpoints. This gives an energy-nonincreasing
deformation onto polygons, commuting with reflection and all grid
shifts. In particular it is $K$-equivariant. No continuous rotation
equivariance of this polygon deformation is needed.

The polygon critical points are the smooth closed geodesics: its first
variation imposes matching velocities at vertices. The splitting of
a variation into its piecewise Jacobi interpolation and a field
vanishing at vertices is orthogonal for the index form. For sufficiently
small cells the latter summand is positive by the Dirichlet Poincar\'e
inequality, so polygon and Hilbert indices and nullities agree.
Regular polygon sublevels are compact smooth manifolds with boundary,
and hence have finite CW models. If $k$ divides $N$, the fixed
manifold of the shift by $1/k$ is the repeated copy of $P_{N/k}$ and
\begin{equation}\label{eq:fp-iteration-energy}
 f_N|_{P_{N/k}}=k^2f_{N/k}.
\end{equation}

The polygon model reduces local questions to finite dimensions. We
now use iteration to make the estimates uniform over all iterates of
a prime. The index stays within a bounded distance of a linear
function of the iteration number. More subtly, every null space is
obtained by repeating one of finitely many root null spaces. The
second fact is what will control degenerate local groups.

\begin{lemma}[Indices and finitely many root kernels]\label{lem:fp-indices}
Let $I_{i,m}$ be the periodic index of $c_i^m$. For each $i$ there are
$\Delta_i\geq0$ and an integer $d_i\geq1$ such that
\begin{equation}\label{eq:fp-index-estimates}
 |I_{i,m+k}-I_{i,m}-I_{i,k}|\leq n,
 \qquad |I_{i,m}-m\Delta_i|\leq n.
\end{equation}
Every null Jacobi field at $c_i^m$ is a repetition of a null field at
$c_i^{\gcd(m,d_i)}$. For even $m$,
\begin{equation}\label{eq:fp-halved-index}
 |I_{i,m}-2I_{i,m/2}|\leq n.
\end{equation}
\end{lemma}
\begin{proof}
Work at prime speed on integer time intervals. In the $(m+k)$-period
form domain require the values at times $m$ and $0$ to agree. In the
direct sum of the $m$- and $k$-period domains require their initial
values to agree. The restricted forms identify, and each restriction
has codimension at most $n$. The two possible index drops lie in
$[0,n]$, proving the first estimate. The sequences $I_{i,m}+n$ and
$I_{i,m}-n$ are respectively subadditive and superadditive, with the
same normalized limit $\Delta_i$: their normalized terms differ by
$2n/m$, and the subadditive and superadditive limit formulas apply.
Applying their inequalities to
multiples of a fixed $m$ gives
$(I_{i,m}-n)/m\leq\Delta_i\leq(I_{i,m}+n)/m$.
The indices are nonnegative, hence so is $\Delta_i$. Taking $k=m/2$
in the first estimate proves \eqref{eq:fp-halved-index}.

For the nullity statement, let $P_i$ be the real Jacobi Cauchy-data
monodromy over one prime
period, and let $d_i$ be the least common multiple of the orders of
its root-of-unity eigenvalues. Put $d_i=1$ if there are none. If
$s=\gcd(m,d_i)$, the same eigenvalues satisfy $z^m=1$ and $z^s=1$.
On a contributing complex Jordan block each polynomial has a simple
root, so its kernel is the eigenspace, with no further generalized
eigenvectors. Thus
\[
 \ker(P_i^m-1)=\ker(P_i^s-1).
\]
Uniqueness for the Jacobi equation identifies the fields as repeated
root fields. This computation is over $\mathbb R$ or $\mathbb C$;
it involves no division by iteration multiplicities in $\F$.
\end{proof}

\paragraph{Ordinary local groups and index pairs.}
At a nondegenerate critical point, local Morse data are determined by
a negative disk. For a degenerate critical circle, the tangential
part can be more complicated. We retain that part in an index pair,
which records a neighborhood and the portion through which descending
trajectories leave it. Adding nondegenerate normal directions then
changes the local cohomology only by the negative rank. This is the
comparison we need for repeated geodesics.

We use the following standard finite-dimensional Conley-index facts.
A compact isolated invariant set has an isolating block and an index
pair; the pointed quotient of an index pair is independent of the
pair up to homotopy; and it is invariant under continuation through
a common isolating neighborhood. For index-pair independence and
continuation we use Salamon's published paper
\cite[Lemma~4.7, Theorems~4.10 and~6.7]{Salamon1985}.
For the block existence and the homological continuation and block-pair
identification, see also \cite[Proposition~2.3, Theorems~5.5 and~7.1]{RV2014}.
These statements are applied to finite-dimensional polygon manifolds.
Index-pair homotopy invariance is independent of coefficients, and we
use singular cohomology over $\mathbb F_2$ throughout. No change of
coefficients from an asserted integral homology formula is needed.

Here an index pair $(N,L)$ is a compact pair for which the invariant
set lies in $\operatorname{int}(N\setminus L)$, $L$ is positively
invariant relative to $N$, and every forward trajectory leaving $N$
first meets $L$. The maximal invariant set in
$\overline{N\setminus L}$ is required to be the isolated set in question.
For good index pairs its cohomological index is $H^*(N,L)$.

\begin{lemma}[Normal stabilization of ordinary local groups]
\label{lem:normal-stabilization}
Let $g$ be smooth on a finite-dimensional manifold $Z$, and let
$C$ be a compact isolated critical set, all of whose points have
the same value $c$. Let $E^-,E^+\to Z$ be real vector bundles on a
neighborhood of $C$, with metrics, and let $b=\operatorname{rank}E^-$.
For
\[
 f(z,u,y)=g(z)-|u|^2+|y|^2,
\]
the ordinary local crossing groups at the zero section over $C$
are isomorphic to those of $g$ at $C$, shifted by $b$.
More precisely, if $(P,A)$ is a sufficiently small isolating
block and exit pair for a descending gradient flow of $g$, the
local cohomology of $f$ is
\[
 H^j\bigl(D(E^-|_P),S(E^-|_P)\cup D(E^-|_A)\bigr)
       \cong H^{j-b}(P,A).
\]
Neither $C$ nor $g$ is required to be nondegenerate.
\end{lemma}
\begin{proof}
\emph{Isolation and continuation.}
Take a compact neighborhood containing $C$ in its interior and no
other critical points. For any smooth vector field that vanishes
on $C$ and strictly decreases $g$ elsewhere in this neighborhood,
its maximal invariant set is $C$. Indeed, a complete trajectory
remaining in that compact neighborhood has nonempty alpha and omega
limit sets. The value of $g$ is constant on each limit set.
Invariance and strict decrease force both limit sets to lie in $C$.
Both limiting values are therefore $c$. Monotonicity then forces
$g$ to be constant on the whole trajectory, which must lie in $C$.
The same argument applies to $f$ near its zero section over $C$.
It also applies to every convex interpolation of two vector fields
strictly decreasing the same function. Consequently a fixed compact
neighborhood with no other critical points is isolating throughout
such an interpolation. This checks the continuation condition, rather
than assuming that a change of descending flow preserves local groups.

\emph{A disk-bundle index pair.}
Choose a small smooth isolating block $P$ for a descending gradient
flow $\psi$ of $g$ and write $A$ for its closed exit set. Use metric
connections on the two bundles and the lifted flow
\[
 \dot z=-\nabla g(z),\qquad
 \nabla_tu=u,\qquad \nabla_ty=-y.
\]
It strictly decreases $f$ off $C$, since
$\dot f=-|\nabla g|^2-2|u|^2-2|y|^2$.
Choose fixed positive fiber radii small enough that the following
pair lies in the normal-form neighborhood; rescale them to one in
the notation:
\[
 B=D(E^-|_P)\mathbin{\times_P}D(E^+|_P),\qquad
 B^-=(B|_A)\cup
       \bigl(S(E^-|_P)\mathbin{\times_P}D(E^+|_P)\bigr).
\]
We check the index-pair conditions explicitly. The base exits $P$
only through $A$. The negative fiber radius increases and the
positive fiber radius decreases, so a forward trajectory can leave
$B$ only through $B^-$. Points of $B^-$ stay in that set for as
long as their forward trajectory stays in $B$; at a negative sphere
they immediately exit. The maximal invariant set in $B$ has $u=0$,
since otherwise its forward radius grows without bound, and $y=0$,
since otherwise its backward radius grows without bound. Its base
therefore lies in $\operatorname{Inv}(P)=C$.
This set lies in the interior away from $B^-$. Thus $(B,B^-)$ is
an index pair. The block collars, together with radial collars,
make it a good pair; corners do not affect the relative construction.

\emph{Removing the positive directions and taking Thom.}
The homotopy $(z,u,y)\mapsto(z,u,(1-t)y)$ preserves $B^-$ and
retracts the pair onto
$(D(E^-|_P),S(E^-|_P)\cup D(E^-|_A))$.
The relative Thom isomorphism gives the displayed degree shift.
It applies to the actual bundle, including a nontrivial bundle
around $C$; mod-two coefficients remove the orientation system.
Index-pair independence and the continuation verified above identify
these groups with the original descending flow's local groups.

For clarity, the crossing pairs used in this paper are such index
pairs. For a labeled flow lens $(N,Q)$ of
Lemma~\ref{lem:fp-lenses}, the function decreases at its upper
boundary and $f_T$ decreases at the exit $Q$. Hence a forward
trajectory leaves $N$ only through $Q$, and the critical set lies
strictly between the two regular cuts. The preceding strict-Lyapunov
argument identifies its maximal invariant set with that critical
set. The collars in that lemma make this a good index pair.
Thus the index-pair comparison computes the actual ordinary
crossing groups, not a separately defined invariant.
\end{proof}

We apply normal stabilization to a repeated geodesic. The root-kernel
lemma places every degenerate direction in the fixed manifold of the
repetition action; the remaining normal directions are nondegenerate.
The resulting comparison concerns the whole phase circle, so its
negative normal bundle may be nontrivial. This is the reason for the
bundle form of the preceding lemma.

\begin{proposition}[Whole-orbit local Thom comparison]
\label{prop:local-thom-realization}
Under the finite-prime assumption, let $s=\gcd(m,d_i)$.
The ordinary local groups at one oriented phase circle satisfy
\[
 \mathcal C^j(c_i^m)
   \cong\mathcal C^{j-(I_{i,m}-I_{i,s})}(c_i^s).
\]
Here $\mathcal C^j(c_i^m)$ denotes $H^j(N,Q)$ for a sufficiently
small ordinary index pair of that oriented phase circle, equivalently
its labeled local crossing group. This is an assertion about whole-orbit
ordinary local groups, with mod-two coefficients. It makes no assertion that a chosen isomorphism
intertwines restriction to a half-period fixed locus.
\end{proposition}
\begin{proof}
Put $k=m/s$ and choose a polygon order $N$ divisible by $k$ and
sufficiently fine for both the root and iterate. The fixed submanifold
$Z=P_{N/k}$ has restricted energy $g_0=k^2f_{N/k}$.
Choose a cyclically invariant metric and an equivariant tubular
neighborhood of $Z$. The normal derivative vanishes along $Z$:
the gradient there is fixed by the cyclic action and is tangent to
its fixed submanifold. The Hessian splits into tangential and normal
parts at the critical circle. Lemma~\ref{lem:fp-indices} puts the
entire kernel in the tangential part, so the normal Hessian is
nonsingular. Its negative rank is $b=I_{i,m}-I_{i,s}$.

The parameterized Morse lemma in the normal variables gives
$f=g_0(z)-|u|^2+|y|^2$ on the negative and positive normal bundles
over a small neighborhood of the whole circle. To see that this
does not require a trivialization around the circle, write the
normal Taylor term as $\langle A(z,v)v,v\rangle$. Set
$A_0=A(z,0)$, $J=\operatorname{sign}(A_0)$, and
$B_0=|A_0|^{-1/2}A|A_0|^{-1/2}$. Near the zero section $JB_0$
is near the identity and has a smooth real square root $S$.
It is $J$-self-adjoint, and $S^*JS=B_0$. The fiber map
$v\mapsto S|A_0|^{1/2}v$ is a local diffeomorphism, uniformly
near the compact circle. These operations are intrinsic bundle
operations and give the claimed normal form.

Finitely many prime images imply that each of the phase circles
is an isolated critical set in its polygon neighborhood. All its
points have the same energy. Lemma~\ref{lem:normal-stabilization}
therefore applies and gives the shift by $b$. Repetition is a
homeomorphism onto $Z$; its inverse restricts to $[0,1/k]$ and
rescales time. Dividing $g_0$ by $k^2$ changes the positive flow
speed and the energy cuts, but not the local index. This identifies
the base local groups with those at the root circle.

Shrinking isolating neighborhoods or changing regular cuts replaces
one index pair by another for the same isolated set. For polygon
refinement, use a common finer order divisible by both orders.
The embedded coarse polygons have zero derivative in the inserted
vertex variables, since each short segment is a minimizing geodesic.
The additional fixed-endpoint variations on each short time cell
have positive Hessian by the Dirichlet estimate already established
above; their normal stabilization has negative rank zero. The same
lemma identifies the refined and coarse local groups. This also
explains why no choice of polygon order changes the asserted shift.
\end{proof}

This comparison gives the dimensions and degree supports of the
ordinary local groups. The doubling argument also needs a specific
map: restriction to loops fixed by half-period translation. To retain
that map, we construct the ambient and fixed local pairs from the same
flow. Their literal inclusion, rather than a chosen Thom isomorphism,
will define every restriction used below.

\begin{lemma}[A common local comparison]\label{lem:fp-lenses}
At a sufficiently tight crossing of a positive critical value there
are local pairs indexed by the oriented phase circles. They compute
the ordinary and equivariant crossing groups, and the actual
restriction to $\operatorname{Fix}(\tau)$ is obtained by literal
intersection of these pairs with that fixed manifold. The comparison
commutes with all relative maps, connecting maps and ambient cup
products.
\end{lemma}
\begin{proof}
The local pairs will be defined by a single descending flow. This is
needed because an isomorphism between abstract local groups alone
does not identify the restriction map to the half-period fixed locus.
Choose one sufficiently fine even polygon model and one complete
$K$-invariant strict descending flow $\phi$ for the original polygon
energy $f$. The gradient is tangent to $\operatorname{Fix}(\tau)$;
there a critical point of the restricted energy is precisely a full
critical point, since the full gradient is fixed. For regular cuts
$a_-<c<a_+$ in energy and a time $T>0$ define
\begin{align*}
 X&=\{f\leq a_+\},&X_-&=\{f\leq a_-\},&f_T&=f\circ\phi_T,\\
 Y&=\{x\in X:f_T(x)\leq a_-\},&
 N&=\{x\in X:f_T(x)\geq a_-\},&Q&=N\cap Y.
\end{align*}
Thus $N$ consists of points that remain above the lower cut for time
$T$, and $Q$ is their exit at time $T$. We call $(N,Q)$ the
\emph{flow lens}. The notation $a_\pm$ in this construction denotes
energy values; the corresponding length cuts are $\sqrt{a_\pm}$.
The inclusions
\begin{equation}\label{eq:fp-lens-diagram}
 (X,X_-)\longrightarrow(X,Y)\longleftarrow(N,Q)
\end{equation}
induce relative cohomology isomorphisms. For the first, $\phi_T$ is
an inverse as a map of pairs, with the intermediate flow as homotopy:
the flow preserves $X$, $X_-$ and $Y$. For the second, the finite
closed cover $X=N\cup Y$ gives the quotient homeomorphism
$N/Q\cong X/Y$ by closed pasting. Both pairs are good pairs.
Indeed the derivative of $f_T$ along the flow is negative at
$f_T=a_-$, so transverse first hitting gives a continuous neighborhood
deformation to $Y$ in $X$, and, from the upper side stopped at
equality, to $Q$ in $N$. Empty collapse sets are treated by adjoining
a disjoint basepoint.

Here are the pair and quotient details. For $0\leq t\leq T$,
$\phi_t(Y)\subset Y$ because
$f(\phi_T\phi_t x)\leq f(\phi_Tx)$; also
$\phi_t(X_-)\subset X_-$ and $\phi_t(X)\subset X$.
Consequently the homotopies comparing the inclusion and $\phi_T$
are homotopies of the indicated pairs, not just of their ambient spaces.
At $f_T=a_-$ the derivative in the flow direction is strictly negative:
otherwise $\phi_Tx$ would be a critical point at the regular level $a_-$.
Compactness of this level in the polygon sublevel makes a transverse
collar uniform. Stopping the flow at its first hit of $f_T=a_-$ gives
the neighborhood deformation of each collapse set. On $X$ it is
stationary on $Y$; on $N$ it is stationary on $Q$.
Both stopping rules commute with $K$.
The resulting good-pair quotient argument is therefore equivariant.
After Borel construction the same collars give relative excision;
one may equivalently apply it on finite-dimensional free $K$-models
and pass degreewise to the Borel groups. Compactness of the infinite
Borel construction is not used.

The complete comparison with the fixed locus is the inclusion diagram
\[
\begin{array}{ccccc}
(X^\tau,X_-^\tau)&\longrightarrow&(X^\tau,Y^\tau)
  &\longleftarrow&(N^\tau,Q^\tau)\\
\downarrow&&\downarrow&&\downarrow\\
(X,X_-)&\longrightarrow&(X,Y)&\longleftarrow&(N,Q).
\end{array}
\]
Every square commutes as a diagram of spaces. The pair homotopies
and collars restrict to the top row because the invariant flow is
tangent to the fixed submanifold. Naturality of the long exact
sequence and of cup product is consequently inherited from maps of
pairs. Replacing a vertical arrow by an arbitrary isomorphism of
local cohomology groups would not have this consequence.

Choose disjoint neighborhoods $\mathcal W_\alpha$ of the critical
circles, permuted by $K$, and smaller neighborhoods
$\mathcal V_\alpha$ separated by positive distance from the
complements of $\mathcal W_\alpha$. On a fixed isolating band the
drop outside the $\mathcal V_\alpha$ has a positive lower bound.
A sufficiently long segment remaining in the band therefore enters
some $\mathcal V_\alpha$. Since the flow is bounded, crossing the
annulus to $\mathcal W_\alpha^c$ takes a positive amount of time and
costs a fixed positive amount of energy. Now choose $a_+-a_-$ smaller
than all these costs. A segment beginning below $a_+$ and ending
above $a_-$ cannot cross such an annulus either before or after its
visit to $\mathcal V_\alpha$. It stays in one $\mathcal W_\alpha$.
Consequently $N$ and $Q$ split into clopen labeled pairs.

Every construction is equivariant. The finite closed-cover argument
and the neighborhood homotopies also hold after Borel construction;
no Borel space is being assumed compact. Restricting the same flow to
the fixed manifold gives exactly
\[
 N^\tau=N\cap\operatorname{Fix}(\tau),\qquad
 Q^\tau=Q\cap\operatorname{Fix}(\tau).
\]
Thus both the full and fixed comparison diagrams consist of the
same inclusion maps. This proves compatibility with the relative
sequences, subgroup restrictions and ambient cup products.
In particular, one first intersects $(N,Q)$ with the fixed manifold
and only then passes to relative cohomology. No choice of an
isomorphism between independently constructed local models is involved.
\end{proof}

The two local constructions now give complementary information.
Normal stabilization bounds the size and degree range of each
crossing; the common flow preserves the maps and the ambient cup
action. The next lemma collects these facts and proves the local
vanishing of $h$ that will ultimately contradict global pairing.

\begin{lemma}[Uniform local groups and evaluation cup product]
\label{lem:fp-local}
Under the finite-prime assumption, there are constants $A,B<\infty$, depending only on the metric and
the finite prime list, with the following properties. For an
unoriented label $(i,m)$ let $P_{i,m}$ be its local pair in one
orientation. Then
\begin{equation}\label{eq:fp-local-support}
 \sum_j\dim H^j(P_{i,m})\leq B,
 \qquad H^j(P_{i,m})=0\quad\text{if }|j-I_{i,m}|>A.
\end{equation}
The reflection block is induced from this ordinary pair; its $K$-block
is induced from the $\langle\tau\rangle$-pair in that orientation.
Odd $m$ has no fixed label. For even $m$, the ordinary fixed pair has
total rank at most $B$ and support $|j-I_{i,m/2}|\leq A$.
The actual restriction respects labels and the naturalities of
Lemma~\ref{lem:fp-lenses}. Cup product by $h$ vanishes on the whole
reflection crossing, including a crossing with simultaneous labels.
\end{lemma}
\begin{proof}
\emph{Orientations and fixed labels.}
The two orientations are exchanged freely by $\rho$, and $\tau$ preserves
each orientation circle. It fixes the latter pointwise exactly when
$m$ is even. When $m$ is odd choose its neighborhood disjoint from
the fixed locus. The assertions about induced blocks and their
restriction maps now follow from the common local pairs.

\emph{Reduction to finitely many ordinary local groups.}
An ordinary labeled pair can be
isolated without changing its local function: add small nonnegative
bumps constant near the other critical circles at the same level and
supported away from the selected neighborhood. On their compact
transition regions the original strict drop is bounded away from
zero, so sufficiently small bumps leave that flow strict and create
no new critical points. Common tight cuts and the trapping argument
identify the isolated crossing with the original labeled pair.
Its ordinary groups may therefore be calculated with a different
strict flow, while all restriction maps continue to use the original
flow of Lemma~\ref{lem:fp-lenses}.
Only the dimensions and degree supports will be computed in the
following normal form. The common-flow construction continues to
specify every homomorphism used later in the proof.

Put $s=\gcd(m,d_i)$ and $b_-=I_{i,m}-I_{i,s}$.
Proposition~\ref{prop:local-thom-realization}, proved by the
index-pair construction in Lemma~\ref{lem:normal-stabilization}, gives
\begin{equation}\label{eq:fp-normal-shift}
 \dim H^j(P_{i,m})
  =\dim H^{j-(I_{i,m}-I_{i,s})}(P_{i,s}).
\end{equation}
The relative Thom map used there pulls back a relative cocycle on
$(P,A)$ and cups with the negative bundle's Thom class. Mod-two
coefficients allow this bundle to be nonorientable. Repetition
identifies the root phase circle homeomorphically with the repeated
phase circle; no division by the iteration multiplicity is involved.
Only dimensions and supports are transferred by this comparison.
Every half-period restriction map continues to be the actual map
in Lemma~\ref{lem:fp-lenses}. In particular, no compatibility between
a chosen Thom isomorphism and that restriction is assumed or used.
Forgetting this distinction would invalidate the later defect count:
the Thom comparison supplies ranks and degree supports, while the
common-flow inclusion supplies the homomorphism.

\emph{Uniformity over all iterates.}
There are only finitely many root labels $(i,s)$ with $s\mid d_i$.
Choose fixed regular isolating cuts for their levels and one coarse
polygon model adequate for all those cuts. Its regular sublevels
are compact manifolds with boundary, so all root crossing groups
have finite total rank and uniformly bounded degree support. Each
root label is a direct summand of its crossing by
Lemma~\ref{lem:fp-lenses}. At any iterate event choose a finer
polygon order divisible by every needed $k$ times this coarse order;
for each even label impose the analogous divisibility for its half
iterate as well. Polygon retractions and regular-strip deformations
identify the fine and coarse root groups. Equation
\eqref{eq:fp-normal-shift} now proves \eqref{eq:fp-local-support}
with uniform constants. Applying the same argument on
$\operatorname{Fix}(\tau)$, whose energy is four times the root
energy, gives the assertion for even labels. No step divides by
$\Delta_i$; primes of zero mean index are included.

\emph{Vanishing of the evaluation cup product.}
Take a sufficiently small $R$-invariant
tubular neighborhood of each pair of reversed phase circles.
Its components are interchanged by reflection, and its reflection
Borel construction retracts to one phase circle. Every cohomology
class of degree at least two is zero there; in particular the ambient
evaluation class $h$ restricts to zero. Choose the original common
lens inside this neighborhood (and intersect its translates to keep
the choices $K$-invariant). Cup by $h$ is then zero on its relative
reflection groups. The inclusion isomorphisms in
\eqref{eq:fp-lens-diagram} transfer this statement to the actual
crossing. A finite disjoint union of such neighborhoods handles
simultaneous labels.
This uses the degree assumption $d\geq2$ precisely at the phase
circle: its cohomology vanishes in that range, independently of the
index or nullity of the geodesic.
\end{proof}

For completeness, we clarify the relation between the index-pair
construction and a disk-bundle quotient. A radius which shrinks to
zero on the exit set gives the usual relative Thom quotient. This is
a statement about quotient spaces; comparison with a flow lens still
requires the continuation argument above.

\begin{lemma}[Collapsing a variable-radius disk bundle]
\label{lem:variable-radius-thom}
Let $(P,A)$ be a compact Hausdorff pair with $A$ closed, let $E\to P$
be a finite-rank real vector bundle with a continuous norm, and let
$r:P\to[0,\infty)$ be continuous with zero set exactly $A$.
Set
\[
 D_r=\{(z,u):|u|\leq r(z)\},\qquad
 S_r=\{(z,u):|u|=r(z)\}.
\]
Scaling $(z,u)\mapsto(z,r(z)u)$ induces a homeomorphism
\[
 D(E)/(S(E)\cup D(E)|_A)\ \cong\ D_r/S_r.
\]
Empty collapse sets are interpreted by adjoining a basepoint.
\end{lemma}
\begin{proof}
The scaling map is continuous and sends the entire collapsed set
to $S_r$. Away from the collapsed set $r(z)>0$, and its inverse
is division by $r(z)$. Thus it induces a continuous bijection on
the quotients. The source is compact; the target is Hausdorff
because $D_r$ is compact Hausdorff and $S_r$ is closed. The bijection
is a homeomorphism. The same reasoning covers rank zero with the
stated basepoint convention. Division by a vanishing radius is
never used to define a map on the uncollapsed disk bundles.
\end{proof}

\begin{remark}[Why a quotient calculation alone is insufficient]
For a split flow with $u(t)=e^tu$ and $y(t)=e^{-t}y$, the expression
$r(z)=e^{-T}\sqrt{g_0(\psi_Tz)-a_-}$ defines a variable radius over
a base lens. Lemma~\ref{lem:variable-radius-thom} identifies its
bundle quotient. It does not by itself deform a full flow lens
onto that quotient: intermediate points could leave the normal-form
neighborhood. The proof of Proposition~\ref{prop:local-thom-realization}
avoids this inference. It constructs a disk-bundle index pair directly,
checks its exit set, and uses continuation in a common isolating
neighborhood to compare it with the original lens. The actual
fixed-locus restrictions are computed separately by the literal
inclusions of Lemma~\ref{lem:fp-lenses}.
\end{remark}

\subsection{Finite-cut restriction and its defect}
\label{subsec:fp-defect}

Half-period-fixed loops at a cut are precisely repetitions of loops
at half that cut. Restriction to this fixed locus therefore compares
two different length levels. Equivariant localization makes the
comparison exact in high Borel degree. For the proof we need more:
a bound on the failure of exactness in the remaining degrees, uniform
as the cut grows. We first identify the actual restriction map and
then introduce a quantity which can be estimated one crossing at a
time.

For $a>0$, define the repetition map
\[
 d_2:X_a\longrightarrow X_{2a}^{\tau},\qquad
 d_2(\gamma)(t)=\gamma(2t).
\]
It is a homeomorphism, with inverse obtained by restricting a
half-period-invariant loop to $[0,1/2]$ and rescaling the parameter.
It preserves the reflection action, and the half-period generator
$\tau$ acts trivially on its target. The restriction considered below
is induced by the inclusion
$i_a:(X_{2a}^{\tau},M)\hookrightarrow(X_{2a},M)$.
We first prove that this particular map is an isomorphism after
inverting $v$, and then measure its failure to be an isomorphism
before localization.

\begin{lemma}[Localization with the actual restriction map]
\label{lem:fp-localization}
If $2a$ is regular, then $a$ is regular and the restriction
\[
 H_K^*(X_{2a},M)[v^{-1}]
       \longrightarrow H_K^*(X_{2a}^{\tau},M)[v^{-1}]
\]
is an isomorphism. It is an actual isomorphism in all sufficiently
high degrees, with a threshold depending on $a$. The localization
assertion also holds on crossing pairs, label by label. There is a
constant $A_1$ such that the actual $K$-group of any label $(i,m)$
has uniformly bounded dimension in every degree, vanishes below
$I_{i,m}-A_1$, and restricts isomorphically to its stable fixed target
above $I_{i,m}+A_1$. For odd $m$ this target is zero; for even $m$
its root degrees lie in $|l-I_{i,m/2}|\leq A_1$.
\end{lemma}
\begin{proof}
\emph{Localization.}
A critical loop at cut $a$ doubles to one at $2a$, proving regularity
of $a$. Repetition identifies $X_a$ with $X_{2a}^{\tau}$ as an
$R$-space, and $\tau$ acts trivially on the latter. Thus
\begin{equation}\label{eq:fp-fixed-decomposition}
 H_K^J(X_{2a}^{\tau},M)
       =\bigoplus_{l\leq J}v^{J-l}B_a^l.
\end{equation}
The same formula holds for fixed crossing pairs.

Choose a compact regular even-polygon sublevel through the ambient
cut and double it across an invariant boundary collar to obtain a
closed smooth $K$-manifold $D_0$. Inclusion and folding make the
original restriction arrow a retract of the arrow for $D_0$.
Let $O_0$ be an invariant tubular neighborhood of $D_0^\tau$.
Its compact complement is covered by finitely many invariant open
sets $O_1,\ldots,O_s$, each admitting an equivariant map to $K/H_i$
with $\tau\notin H_i$. Such a neighborhood is obtained by taking a
small stabilizer-invariant neighborhood disjoint from the translates
by other cosets and then saturating it. The possible stabilizers are
$1$, $R$ and $\langle\tau\rho\rangle$. The class $U=v(v+w)$
restricts to zero on each, hence lifts to $H_K^2(D_0,O_i)$.
Multiplying these $s$ lifts with a class relative to $O_0$ yields
a class relative to $\bigcup_{i=0}^sO_i=D_0$, and is therefore zero.
Thus $U^sH_K^*(D_0,O_0)=0$. The pair sequence proves localization
at $U$ for the restriction to the fixed locus. Passing to the
retract and then using the exact sequences of regular sublevel
pairs proves the same result relative to constants and on crossings.
For individual labels, use the direct decomposition and actual
labelwise restriction in Lemma~\ref{lem:fp-local}.

\emph{From $U$ to $v$.}
On a local $K$-block, $w$ is zero: that block is induced from the
$\langle\tau\rangle$-pair in one orientation. Finitely many crossing
exact sequences show that $w$ is nilpotent on every finite-cut
relative $K$-group. The same holds for the fixed groups, using their
reflection local blocks. Localization at $U$ is therefore equivalent
to localization at $v$. After inverting $v$, the factor $1+w/v$ in
$U=v^2(1+w/v)$ is a unit; after inverting $U$, $(v+w)/U$ is an
inverse for $v$.

\emph{Uniform local thresholds.}
Apply the sign-line sequence to the $\langle\tau\rangle$-pair of
one orientation. Lemma~\ref{lem:fp-local} gives a uniform bound
for every equivariant degree by the total ordinary rank, vanishing
below $I_{i,m}-A$, and invertibility of multiplication by $v$ above
$I_{i,m}+A$ after a fixed enlargement. On an even fixed block,
\eqref{eq:fp-fixed-decomposition} and
\eqref{eq:fp-halved-index} show that its target has also stabilized
in that range. Both source and target multiplication maps are then
invertible. Clear a power of $v$ in the localized restriction
isomorphism and invert these multiplication maps; the restriction
itself is an isomorphism. The odd target is empty, so its stabilized
source is zero. This proves the uniform local assertions.

\emph{A high-degree isomorphism at each finite cut.}
At a finite cut, $B_{2a}$ has finite total dimension. It starts at
zero near the constants and is built by finitely many reflection
crossings, each of finite total rank. The sign sequence
\eqref{eq:fp-sign-K} consequently makes multiplication by $v$ an
isomorphism in all sufficiently high ambient degrees. Formula
\eqref{eq:fp-fixed-decomposition} gives the same statement on the
fixed side. The preceding clearing argument proves the high-degree
assertion for this cut.
\end{proof}

To make the finite-cut maps explicit, also write
$i_\infty:(\Lambda^\tau,M)\hookrightarrow(\Lambda,M)$.
For each integer $J\geq0$, repetition identifies the lower row in
the following commutative
diagram with reflection cohomology, together with the polynomial
class $v$:
\begin{equation}\label{eq:fp-finite-cut-diagram}
 \begin{array}{ccc}
 H_K^J(\Lambda,M)&\longrightarrow&H_K^J(X_{2a},M)\\
 \big\downarrow i_\infty^*&&\big\downarrow i_a^*\\
 \displaystyle\bigoplus_{l\leq J}v^{J-l}B_\infty^l
 &\longrightarrow&
 \displaystyle\bigoplus_{l\leq J}v^{J-l}B_a^l.
 \end{array}
\end{equation}
The horizontal arrows are restrictions of pairs. In particular, the
bottom arrow restricts each coefficient separately; it is not an
arbitrary identification of vector spaces of the same dimension.

For a regular paired cut $(2a,a)$ put
\[
 V_a=\bigoplus_l B_a^l,\qquad
 A^J(a)=H_K^J(X_{2a},M).
\]
There are two gradings in this comparison. A power of $v$ raises the
Borel degree without changing the underlying reflection class on the
half-length space. To count the latter classes, we separate that
power from its coefficient.

The integer $J$ is the total Borel degree. The integer $l$ in
\eqref{eq:fp-fixed-decomposition} is the degree of the coefficient
in $B_a$, which we call its \emph{root degree}. Define a linear
normalization on the fixed target by
\[
 \nu_J\left(\sum_{l\leq J}v^{J-l}b_l\right)
       =\sum_{l\leq J}b_l,\qquad b_l\in B_a^l.
\]
The actual restriction becomes a linear map
\[
 \epsilon_J=\nu_J\circ i_a^*:A^J(a)\longrightarrow V_a,
 \qquad E_J=\operatorname{im}\epsilon_J,
 \qquad t_J=\dim\ker\epsilon_J.
\]
Every vector of $E_J$ has root degrees at most $J$, but $E_J$ need
not be a graded subspace. In these coordinates multiplication by
$v$ leaves root coefficients unchanged, whereas multiplication by
$w$ acts by the root $w$ and raises root degree by one.
Thus, for $u\in A^J(a)$,
\[
 \epsilon_{J+1}(vu)=\epsilon_J(u),\qquad
 \epsilon_{J+1}(wu)=w\epsilon_J(u).
\]
These identities retain the total degree on the left even when the
vector on the right has components in several root degrees.
A single kernel or cokernel is not enough to measure the comparison:
its image may mix several root degrees. We must keep track both of
what the image misses in low root degrees and of the image vectors
lying entirely in high root degrees. Two neighboring Borel degrees
are needed because the crossing exact sequence has a connecting map.
The following four-term quantity records this information.

Let $P_j$ project to root degrees at most $j$, and define
\begin{align}
 D_{j,q}&=\dim P_jV_a-\dim P_jE_q,
 &H_{j,J}&=\dim(E_J\cap\ker P_j),\nonumber\\
 \Phi_j(a)&=D_{j,2j}+D_{j,2j-1}+H_{j,2j}+t_{2j},
 &&j\geq1.
 \label{eq:fp-potential}
\end{align}
These four terms are nonnegative. For each fixed cut they vanish
for all sufficiently large $j$, by Lemma~\ref{lem:fp-localization}.
Here $D_{j,q}$ is the deficit of the projected image in root degrees
at most $j$; the projection $P_jE_q$ may differ from
$E_q\cap P_jV_a$. The number $H_{j,J}$ counts image vectors supported
entirely above $j$, and $t_J$ counts the kernel of restriction before
normalization. Keeping both projected images and high-degree vectors
is necessary because $E_J$ need not split as a sum of root-degree
subspaces.

\begin{proposition}[Summed finite-cut defect]\label{prop:fp-budget}
There is a constant $C_{\rm def}$ such that for every regular $2a>0$,
\begin{equation}\label{eq:fp-budget}
 \sum_{j\geq1}\Phi_j(a)\leq C_{\rm def}(1+a).
\end{equation}
\end{proposition}
\begin{proof}
The estimate follows by adding the changes at individual events. We
will show that an event can increase $\Phi_j$ only when $2j$ is
within a uniformly bounded distance of one of its indices. The
remaining values of $j$ contribute a nonpositive change.

Grow the ambient cut through its positive critical events and use
primes for the groups just after one event. Denote its actual local
$K$-group by $C^q$ and the total fixed root reflection crossing
by $W_0$. Thus $A^q$ and $V$ refer to the cuts before the event,
whereas $A'^q$ and $V'$ refer to the cuts after it. Exactness and
naturality give the comparison of rows
\begin{equation}\label{eq:fp-event-rows}
 \begin{array}{ccccccccc}
 A^{q-1}&\longrightarrow&C^q&\longrightarrow&A'^q
       &\longrightarrow&A^q&\longrightarrow&C^{q+1}\\
 \downarrow\epsilon_{q-1}&&\downarrow\lambda_q&&
 \downarrow\epsilon'_q&&\downarrow\epsilon_q&&
 \downarrow\lambda_{q+1}\\
 V&\xrightarrow{\delta}&W_0&\xrightarrow{f}&V'
       &\xrightarrow{i}&V&\xrightarrow{\delta}&W_0.
 \end{array}
\end{equation}
The bottom row is obtained by inverting $v$ and normalizing in
neighboring degrees. It is the total graded fixed sequence: $f,i$
preserve root degree and $\delta$ raises it by one. The repeated
symbols in the row do not assert ungraded periodicity.
For example, the displayed $\delta$ is the sum of maps
$V^l\to W_0^{l+1}$ over root degrees $l$. The maps
$\epsilon_q$ retain the indicated total degree $q$, even though
their images in $V$ may have several root-degree components.

\emph{Events outside the index window.}
Fix $j\geq1$ and put $J=2j$. Call $j$ far from the event if every
participating iterate satisfies $|I_{i,m}-2j|>C_0$, where $C_0$ is a
sufficiently large fixed constant. The local assertions and index
comparison give
\begin{equation}\label{eq:fp-far-event}
 (W_0)^j=(W_0)^{j+1}=0,\qquad
 \lambda_{J-1},\lambda_J,\lambda_{J+1}
       \text{ injective with image }P_jW_0.
\end{equation}
Indeed a block below the window is already stable in all three
actual degrees and has all root degrees below $j$. A block above
the window has zero actual source and no root degree at most $j$.

For this far event abbreviate $P=P_j$, $L=PV$, $L'=PV'$, and put
\[
 A_0=f(PW_0),\qquad r_0=\operatorname{rank}(\delta|_L),
 \qquad r_q=\operatorname{rank}(\delta|_{PE_q})
       \quad(q=J-1,J).
\]
There is no connecting term across the projection edge, since the
two root degrees in \eqref{eq:fp-far-event} vanish. Truncating the
fixed exact sequence therefore gives
\begin{equation}\label{eq:fp-low-exact}
 \dim PW_0=\dim A_0+r_0,
 \qquad \dim L'-\dim L=\dim A_0-r_0.
\end{equation}

Mixed low and high components require retaining the actual maps.
For $z\in E_q$, commutativity in \eqref{eq:fp-event-rows} puts
$\delta z$ in $PW_0$. The term $\delta(Pz)$ is low, because
$(W_0)^{j+1}=0$, and $\delta((1-P)z)$ is high. Their supports
are disjoint, so
\begin{equation}\label{eq:fp-mixed-delta}
 \delta((1-P)z)=0,\qquad \delta z=\delta(Pz).
\end{equation}
This support argument is the point at which the absence of a
root-degree splitting of $E_q$ must be taken into account. We project
a restricted actual class and retain its high component until the
connecting homomorphism is shown to kill that component.
If $x\in PE_q\cap\ker\delta$, choose $z\in E_q$ with $Pz=x$
and then an actual representative in $A^q$. Its actual connecting
image restricts to zero by \eqref{eq:fp-mixed-delta}; injectivity
of $\lambda_{q+1}$ makes that connecting image zero. Actual
exactness lifts the representative to $A'^q$. Hence
\[
 i(PE'_q)=PE_q\cap\ker\delta.
\]
The kernel of this projected map is exactly $A_0$. One inclusion
follows from fixed exactness. Conversely every element of $PW_0$
has an actual local lift by \eqref{eq:fp-far-event}, whose image
in $A'^q$ restricts to the desired vector in $A_0$.
It follows that
\begin{equation}\label{eq:fp-deficit-change}
 D'_{j,q}-D_{j,q}=-(r_0-r_q),\qquad q=J-1,J.
\end{equation}
The actual connecting ranks are $r_{J-1}$ and $r_J$, again by
injectivity of the local restrictions. Thus
\[
 \dim A'^J-\dim A^J=\dim PW_0-r_{J-1}-r_J.
\]
Indeed exactness first removes the image of
$A^{J-1}\to C^J$, of rank $r_{J-1}$, and then removes the image
of $A^J\to C^{J+1}$, of rank $r_J$. The local identification in
\eqref{eq:fp-far-event} gives $\dim C^J=\dim PW_0$.
Since $H_{j,J}+t_J=\dim A^J-\dim PE_J$, this identity and
\eqref{eq:fp-low-exact} give
\[
 (H_{j,J}+t_J)'-(H_{j,J}+t_J)=r_0-r_{J-1}.
\]
Adding this to the two identities
\eqref{eq:fp-deficit-change} yields the exact formula
\begin{equation}\label{eq:fp-far-change}
 \Phi'_j-\Phi_j=-(r_0-r_J)\leq0.
\end{equation}
The preceding-degree deficit in \eqref{eq:fp-potential} is what
cancels the term $r_{J-1}$.

\emph{Events inside the index window.}
For a cutoff not far from the event, let $M_0$ bound its total
fixed local dimension and all the relevant actual local dimensions.
It is uniform, since at most one iterate of each prime occurs at
an event. The kernel of $i:PE'_q\to PE_q$ has dimension at most
$M_0$. Its cokernel is bounded by the surjections
\[
 A^q/\operatorname{im}A'^q\longrightarrow E_q/iE'_q
       \longrightarrow PE_q/i(PE'_q),
\]
whose first source has dimension at most $\dim C^{q+1}\leq M_0$.
Thus $|\dim PE'_q-\dim PE_q|\leq M_0$. Fixed exactness similarly
gives $|\dim L'-\dim L|\leq M_0$, and actual exactness gives
$|\dim A'^J-\dim A^J|\leq M_0$. The identity
\[
 \Phi_j=2\dim L-2\dim PE_J-\dim PE_{J-1}+\dim A^J
\]
therefore bounds $|\Phi'_j-\Phi_j|$ by $6M_0$.
For clarity, each label contributes at most $\lfloor C_0\rfloor+2$
integers satisfying $|I_{i,m}-2j|\leq C_0$. Since an event contains
at most $p$ labels, at most $p(\lfloor C_0\rfloor+2)$ integers
fail the far condition. Each event increases the sum of
potentials by at most one fixed constant. They start at zero below
the first positive event, are unchanged on regular strips, and
have finite support at each cut. The number of events below $2a$
is at most $\sum_i\lfloor2a/\ell_i\rfloor$ by
\eqref{eq:fp-event-count}. Summing proves \eqref{eq:fp-budget}.
\end{proof}

\subsection{A uniform bound for transient reflection classes}
\label{subsec:fp-transients}

The defect estimate controls how much information can be lost in
passing to the half-length fixed locus. We now combine it with the
ordinary filtration estimates. The aim is to bound the reflection
classes at a finite cut which are not restrictions of global classes.
Precisely, we call their quotient space the space of
\emph{transient classes}. This cohomological definition concerns
extension from the full space; an individual transient need not have
a canonical continuation to a larger cut.

The proof first locates possible transients near finitely many lines
in the degree--length plane. Away from these lines, multiplication by
$w$ would propagate any transient through a long interval of degrees.
Such a string would force quadratic growth in the summed defect,
contrary to the linear bound just proved.

Let
\[
 Q_a^j=\operatorname{coker}(B_\infty^j\longrightarrow B_a^j),
 \qquad T(a)=\sum_j\dim Q_a^j,
 \qquad \alpha=\lambda^{-1},\quad\beta_i=\Delta_i/\ell_i.
\]
These sums are finite at every regular cut. List the distinct slopes
as
\[
 \{0,\alpha,\beta_1,\ldots,\beta_p\}
       =\{s_0<s_1<\cdots<s_N\},\qquad s_0=0.
\]
The slope $\alpha$ describes where global homology appears, and the
$\beta_i$ describe where iterates of the individual primes contribute
local groups. Between consecutive slopes there is a widening interval
with no local index contributions near the cut. We call these
intervals degree gaps and leave a fixed margin at their ends.

For consecutive slopes $s<s'$ and a fixed margin $c>0$, write
\[
 I_{s,s'}(a)=\{j\in\mathbb Z_{\geq0}:sa+c<j<s'a-c\}.
\]
A gap is called present when it lies below $\alpha a$, and absent
when it lies above $\alpha a$; $\alpha$ is one of the listed slopes,
so these exhaust the gaps.

\begin{lemma}[Comparison on degree gaps]\label{lem:fp-gaps}
There is a fixed $c$ such that for every regular $a>0$:
\begin{enumerate}
\item $B_\infty^j\to B_a^j$ is injective on every present gap;
\item $B_\infty^j\to B_a^j$ is zero on every absent gap;
\item $w:Q_a^j\to Q_a^{j+1}$ is an isomorphism whenever both
$j,j+1$ belong to the same gap;
\item $B_a^j=0$ if $j>s_Na+c$.
\end{enumerate}
\end{lemma}
\begin{proof}
\emph{Ordinary comparison away from index bands.}
Choose a regular $a_+\in(a+D,a+D+1)$. An ordinary local group at
an event $m\ell_i\in[a,a_+]$ can affect degree $j$ or its required
neighbors only if $|j-I_{i,m}|\leq C_1$, for one fixed $C_1$.
By \eqref{eq:fp-index-estimates}, this implies
\begin{equation}\label{eq:fp-strip-index}
 |j-\beta_i a|\leq C_1+n+\beta_i(D+1).
\end{equation}
Outside these fixed-width bands the crossing sequences make
$H_j(X_a,M)\to H_j(X_{a_+},M)$ injective. If a source class is
globally zero, \eqref{eq:fp-UD} kills it by $a+D<a_+$, so it was
already zero. Thus $H_j(X_a,M)\to H_j(\Lambda,M)$ is injective
outside the bands.

Combine this with the spectral hypothesis. Below
$j<\alpha a-C_2$ that map is also surjective, hence an isomorphism.
Above $j>\alpha a+C_2$ it is zero, and its injectivity forces its
source to vanish. Field duality gives
\begin{equation}\label{eq:fp-ordinary-comparison}
 \begin{cases}
 H^j(\Lambda,M)\longrightarrow H^j(X_a,M)
        \text{ is an isomorphism},&j<\alpha a-C_2,\\
 H^j(X_a,M)=0,&j>\alpha a+C_2,
 \end{cases}
\end{equation}
when $j$ is outside all bands in \eqref{eq:fp-strip-index}.
Here $C_2$ incorporates $\alpha C_{\rm sp}$ and any chosen finite
low-degree range. These are all-class statements; no restriction
on the number of global classes in degree $j$ has been used.

\emph{Reflection comparison below the spectral line.}
In a present gap put $P_a^j=H_R^j(\Lambda,X_a)$. Two neighboring
ordinary isomorphisms in \eqref{eq:fp-ordinary-comparison} imply
$H^j(\Lambda,X_a)=0$ by triple exactness. The sign sequence
\eqref{eq:fp-sign-R} makes $w:P_a^j\to P_a^{j+1}$ an isomorphism
on the corresponding interior. Allow $N_w$ further degrees below
the degree being considered. Every class of $P_a^j$ is then
divisible by $w^{N_w}$, and its image in $B_\infty^j$ is zero,
since $w^{N_w}B_\infty=0$. Applying this in degrees $j,j+1$ to
the triple $(\Lambda,X_a,M)$ gives
\begin{equation}\label{eq:fp-present-sequence}
 0\longrightarrow B_\infty^j\longrightarrow B_a^j
      \longrightarrow P_a^{j+1}\longrightarrow0.
\end{equation}
The quotient identification commutes with $w$. It proves the
injection and the isomorphisms on $Q_a$ in a present gap.

\emph{Reflection comparison above the spectral line.}
In an absent gap, the second part of
\eqref{eq:fp-ordinary-comparison} and
\eqref{eq:fp-sign-R} make $w:B_a^j\to B_a^{j+1}$ an isomorphism
in the interior. Allow $N_w$ further degrees above. Then
$w^{N_w}$ is injective on the degree in question, whereas the
global image is killed by it. That image is zero, so $Q_a^j=B_a^j$
and the same $w$-isomorphisms hold. Choose $c$ larger than all
ordinary band widths plus $2N_w+6$; this simultaneously supplies
all the neighboring degrees just used, uniformly in $a$ and $j$.

Every reflection local contribution below $a$ has
$j\leq I_{i,m}+A\leq\beta_i a+n+A$. Starting from zero near
the constants and summing the crossing sequences therefore gives
$B_a^j=0$ for $j>s_Na+c$, after a final fixed enlargement of $c$.
\end{proof}

\begin{lemma}[Long gaps contain no transients]\label{lem:fp-no-strings}
For all sufficiently large regular $a$ with $2a$ regular,
$Q_a^j=0$ on every gap $I_{s,s'}(a)$.
\end{lemma}
\begin{proof}
Suppose that one transient survives in a long degree gap. The
$w$-isomorphisms of Lemma~\ref{lem:fp-gaps} then produce a long
string of independent root-degree components. We use the defect
bound twice: first to lift a fixed proportion of that string to
actual $K$-classes, and then to count their high-degree restrictions
in a fixed proportion of the possible target degrees.

Fix a consecutive pair $s<s'$, and set
$A_0=sa$ and $r_0=(s'-s)a$. Since there are finitely many distinct
slopes, one fixed $\delta>0$ satisfies $r_0\geq\delta a$ for all
these gaps. Suppose a transient is nonzero somewhere in this gap.
Lemma~\ref{lem:fp-gaps} propagates it, by $w$-isomorphisms, to
every integer in the gap. For large $a$ choose
\[
 l_0=\lceil A_0+0.10r_0\rceil,
 \qquad y_{l_0}\in B_a^{l_0}
\]
with nonzero quotient class, and put $y_l=w^{l-l_0}y_{l_0}$ up to
$l\leq\lfloor A_0+0.90r_0\rfloor$. All their quotient classes
are nonzero. Let
\[
 \mathcal L=\mathbb Z\cap[A_0+0.15r_0,A_0+0.20r_0],
 \qquad q=\lfloor0.50r_0\rfloor,
 \qquad W_1=\operatorname{span}\{y_l:l\in\mathcal L\}.
\]
The vectors $y_l$ are independent because their root degrees are
distinct. The map $w^q$ is injective on $W_1$ for the same reason
and because every resulting quotient class is still nonzero.
In particular, this assertion concerns the actual vectors in $B_a$:
any relation between them would give a relation between their
nonzero homogeneous quotient classes in $Q_a$.

At ambient cut $2a$, Proposition~\ref{prop:fp-budget} implies
\begin{equation}\label{eq:fp-two-term-budget}
 \sum_{k\geq1}(D_{k,2k}+H_{k,2k})\leq C_{\rm def}(1+a).
\end{equation}
Average $D_{k,2k}$ over the integers
$k\in[A_0+0.30r_0,A_0+0.35r_0]$. There are a fixed positive
multiple of $r_0$ such integers for large $a$, and $r_0\geq\delta a$.
For at least one of them $D_{k,2k}\leq C_3$, with $C_3$ independent
of $a$. Since $W_1\subset P_kV_a$,
\[
 \dim(W_1\cap P_kE_{2k})
      \geq\dim W_1-D_{k,2k}\geq c_1r_0
\]
for a fixed $c_1>0$ and sufficiently large $a$. Choose a basis
$z_1,\ldots,z_\nu$ of this intersection and actual lifts
$u_\mu\in H_K^{2k}(X_{2a},M)$. Their normalized restrictions have
the form $z_\mu+c_\mu$, where every root component of $c_\mu$
has degree greater than $k$.
The vectors $z_\mu$ need not be homogeneous in root degree: the
intersection was formed after projecting $E_{2k}$. This is why the
following argument separates two intervals of root degrees instead
of assigning a single degree to each lift.

Multiply the actual classes by $w^q$. The selected components
$w^qz_\mu$ remain independent and lie in root degrees between
$A_0+0.65r_0-O(1)$ and $A_0+0.70r_0$. The other components,
if nonzero, have root degrees greater than
$k+q\geq A_0+0.80r_0-1$. Hence they cannot cancel the selected
components. All root components exceed $A_0+0.64r_0$ for large
$r_0$. For every integer
$k'\in[A_0+0.61r_0,A_0+0.64r_0]$, the actual-degree inequality
\[
 2k+q\leq2A_0+1.20r_0<2A_0+1.22r_0\leq2k'
\]
permits multiplication by the nonnegative power
$v^{2k'-(2k+q)}$. This multiplication changes no normalized root
coefficient. The restrictions of
\[
 v^{2k'-(2k+q)}w^qu_\mu\in H_K^{2k'}(X_{2a},M)
\]
are therefore independent and supported entirely above $k'$. Thus
$H_{k',2k'}\geq c_1r_0$ for each of at least $c_2r_0$ choices
of $k'$, for some fixed $c_2>0$. The left side of
\eqref{eq:fp-two-term-budget} is at least
$c_1c_2r_0^2\geq c_1c_2\delta^2a^2$, contradicting its linear
upper bound. All rounding errors are absorbed by the strict
separations in these intervals for sufficiently large $a$.
There are only finitely many gaps, so one common lower threshold
works for all of them.
All constants have been chosen independently of the gap and the
cut: take the minimum separation of the finite slope set for
$\delta$, use the common defect constant for $C_3$, and then
increase the lower bound on $a$ to absorb that fixed $C_3$ and all
integer endpoints.
\end{proof}

\begin{proposition}[Bounded total transients]\label{prop:fp-transients}
There is a constant $C_T$ such that $T(a)\leq C_T$ for every
positive regular cut $a$.
\end{proposition}
\begin{proof}
The long-gap lemma confines transients to finitely many bands of
bounded width. To bound their total dimension we treat high degrees
by index growth and the remaining fixed range by uniform death.
First suppose $a$ is large and $2a$ is regular. By
Lemmas~\ref{lem:fp-gaps} and \ref{lem:fp-no-strings}, transients
can occur only within distance $c$ of the finitely many degrees
$s_ia$. There are a uniformly bounded number of integer degrees
in these bands. It remains to bound their dimensions uniformly.

For high degrees the local reflection bounds give
\begin{equation}\label{eq:fp-high-rank}
 \dim B_a^j\leq C_4\sum_{i=1}^p
 \#\{m\geq1:m\ell_i<a,\ |j-m\Delta_i|\leq A+n\}.
\end{equation}
This follows from the crossing exact sequences starting with zero
near the constants. If $\Delta_i>0$, the count is bounded by
$2+2(A+n)/\Delta_i$, independent of $a$ and $j$. If $\Delta_i=0$
it is zero for $j>A+n$. Thus all sufficiently high degree groups
$B_a^j$ have a common bound.

For the fixed low-degree interval, again choose a regular
$a_+\in(a+D,a+D+1)$. Uniform death and factorization through
global homology give equality of kernels
\[
 \ker[H_j(X_a,M)\to H_j(\Lambda,M)]
      =\ker[H_j(X_a,M)\to H_j(X_{a_+},M)].
\]
The right kernel is bounded by the sum of ordinary local ranks
in a cut interval of width $D+1$. By \eqref{eq:fp-event-count},
at most $\sum_i(2+(D+1)/\ell_i)$ labels occur there. Uniform
local ranks bound that kernel. Its complementary global-image
dimension is at most the finite number $\dim H_j(\Lambda,M)$.
Consequently each fixed ordinary degree $H^j(X_a,M)$ is bounded
uniformly in $a$. From \eqref{eq:fp-sign-R},
\[
 \dim B_a^j\leq\dim B_a^{j-1}+\dim H^j(X_a,M),
 \qquad B_a^{-1}=0.
\]
Induction, including degree zero, bounds $B_a^j$ on the entire
fixed low-degree interval. This also controls every zero-mean-index
prime, regardless of how many of its iterates occur below $a$.
Together with \eqref{eq:fp-high-rank}, these estimates bound the
sum of transient dimensions in all the remaining bands.

Cuts in a bounded interval have only finitely many crossing
stages. Finally, if $a$ is regular but $2a$ is not, choose a nearby
$a'$ in the same regular interval such that $2a'$ is regular.
Critical values and their halves are locally finite. The regular
strip equivalence is natural with respect to global restriction,
so it preserves both $B_a$ and its global image. Thus
$T(a)=T(a')$, proving the assertion for every regular cut.
\end{proof}

\subsection{Doubling and the first appearance of a class}
\label{subsec:fp-activation}

It remains to turn boundedness of transients into an obstruction.
The Gysin sequences compare full-group cohomology at twice a length
cut with reflection cohomology at the original cut. Separating free
and $U$-torsion parts gives a dimension identity; the global paired
classes then determine when equality can occur.

For a regular cut put
\[
 V(a)=\dim\operatorname{im}(B_\infty\longrightarrow B_a),
 \qquad T(a)=\dim B_a-V(a).
\]
The next calculation uses full $G$-equivariance only on the original
Hilbert sublevels, not on polygon models.
The number $V(a)$ is the dimension of the visible global image,
whereas $V_a$ used earlier is the entire vector space $B_a$. Thus
$T(a)$ is exactly the cokernel dimension already defined in
Section~\ref{subsec:fp-transients}. The finite-dimensional polygon
models established finiteness and local bounds; the Gysin sequences
in the next two results are applied directly to the original
$G$-invariant sublevels.

\begin{lemma}[Free rank and torsion at a finite cut]
\label{lem:fp-free-rank}
Suppose $2L>0$ is regular and set
$H=\bigoplus_jH_G^j(X_{2L},M)$. Then $H$ is a finitely generated
graded $\F[U]$-module. Its torsion submodule $\mathcal T$ is finite
dimensional and killed by one power of $U$. If
\[
 t=\dim(\mathcal T/U\mathcal T),\qquad
 \varrho=\operatorname{rank}_{\F[U]}(H/\mathcal T),
\]
then
\begin{equation}\label{eq:fp-rank-identity}
 \dim B_{2L}=\varrho+2t,\qquad \varrho=\dim B_L.
\end{equation}
\end{lemma}
\begin{proof}
\emph{Finite generation and torsion.}
The standard-plane sequence \eqref{eq:fp-standard-gysin} injects
$H/UH$ into the finite-dimensional space $B_{2L}$. Choose a finite
homogeneous basis of $H/UH$ and lift it to $H$. These lifts generate:
subtract the appropriate combination modulo $U$, write the remainder
as $U$ times a class two degrees lower, and induct on the nonnegative
degree. This proves finite generation.

A polynomial with nonzero constant coefficient acts injectively on
$H$. The lowest nonzero homogeneous component of a nonzero element
survives in its product with that polynomial. Factoring an arbitrary
nonzero annihilating polynomial as $U^sp(U)$ with $p(0)\neq0$
therefore proves that every torsion element is killed by a power
of $U$. Since $\F[U]$ is Noetherian, $\mathcal T$ has finitely
many generators and one power kills them all; it is consequently
finite dimensional. Moreover
\[
 UH\cap\mathcal T=U\mathcal T,
 \qquad \ker(U|_H)=\ker(U|_{\mathcal T}),
 \qquad\dim\ker(U|_H)=t.
\]
For the first equality, $Uu$ torsion implies $u$ torsion; the last
equality is rank-nullity on $\mathcal T$.

\emph{The reflection dimension.}
Set $F_0=H/\mathcal T$. Lift a homogeneous basis of $F_0/UF_0$.
Degree induction again proves generation. The lifts are free: in
a homogeneous relation, divide by the smallest power of $U$
appearing, using injectivity of $U$ on $F_0$, and then reduce
modulo $U$. Independence of the selected basis rules out the
nonzero terms with that smallest power. Thus $F_0$ has a
homogeneous free basis of size $\varrho$, and
\[
 \dim(H/UH)=\varrho+t.
\]
Write $H^j$ for the degree-$j$ part of the module $H$ in the next
display. For each $j$, the standard-plane sequence gives the short
exact sequence
\begin{equation}\label{eq:fp-torsion-short-exact}
 0\longrightarrow H^j/UH^{j-2}\longrightarrow B_{2L}^j
 \longrightarrow
 \ker\bigl(U:H^{j-1}\longrightarrow H^{j+1}\bigr)
 \longrightarrow0.
\end{equation}
Summing this sequence in all degrees gives
$\dim B_{2L}=(\varrho+t)+t$, because its other piece is the shifted
kernel of multiplication by $U$.

\emph{Identification of the free rank.}
For two sufficiently large adjacent degrees $D_0,D_0-1$, torsion
is absent and each free generator contributes once, according to
its parity. Hence
\[
 \varrho=\dim H_G^{D_0}(X_{2L},M)
             +\dim H_G^{D_0-1}(X_{2L},M)
       =\dim H_K^{D_0}(X_{2L},M)
\]
by \eqref{eq:fp-weight-two-gysin}. The actual high-degree
restriction isomorphism of Lemma~\ref{lem:fp-localization} and
\eqref{eq:fp-fixed-decomposition} identify the last dimension
with $\dim B_L$. This proves the second identity without any
surjectivity assertion for full-group restriction to the fixed
locus.
\end{proof}

The free-rank identity accounts for reflection classes at half the
cut. The remaining contribution comes from $U$-torsion. A visible
class in $W$ forces such a torsion contribution through its global
lift, whereas its $h$-multiple accounts for the paired reflection
class. Comparing these contributions proves the inequality below.
The equality condition is the decisive part: it prevents a nonzero
visible $W$-class from being killed by $h$.

\begin{proposition}[Doubling and its equality condition]
\label{prop:fp-doubling}
For every regular $2L>0$,
\begin{equation}\label{eq:fp-doubling}
 T(2L)\geq T(L).
\end{equation}
If equality holds, multiplication by $h$ is injective on the
subspace $\operatorname{im}(W\to B_{2L})$, and that subspace has
zero intersection with its image under $h$.
\end{proposition}
\begin{proof}
We separate the visible global classes according to
\eqref{eq:fp-global-decomposition}, but take dimensions only after
restriction to the finite cut. In this way relations between
different global summands, including relations in the same degree,
are retained in the count.

All spaces in the following rank count are finite dimensional,
although they may occupy overlapping homogeneous degrees. At cut
$2L$ let
\begin{align*}
 O&=\operatorname{im}(W\longrightarrow B_{2L}),&N&=\dim O,\\
 V_o&=\dim(O+hO),&
 V_e&=\dim\left(
 \operatorname{im}(B_\infty\to B_{2L})/(O+hO)\right).
\end{align*}
By \eqref{eq:fp-global-decomposition}, the last quotient is
generated by restrictions of $E_0$, and
\begin{equation}\label{eq:fp-visible-pairs}
 V(2L)=V_o+V_e,\qquad V_o\leq2N.
\end{equation}
The operation on $O$ is well defined: if a global $W$-class
restricts to zero, its product with the ambient class $h$ also
restricts to zero.

\emph{Torsion forced by the visible $W$-classes.}
Choose homogeneous global $W$-classes whose restrictions are a
basis of $O$, and take their full-group $U$-torsion lifts. At the
finite cut these lifts belong to $\mathcal T\subset H$, in the
notation of Lemma~\ref{lem:fp-free-rank}. Their classes modulo
$UH$ are independent, because $H/UH\to B_{2L}$ is injective and
their reflection restrictions are independent. They therefore span
an $N$-dimensional subspace of $\mathcal T/U\mathcal T$. Write
\begin{equation}\label{eq:fp-extra-torsion}
 t=N+J_0,\qquad J_0\geq0.
\end{equation}

\emph{Free rank forced by the remaining visible classes.}
Choose homogeneous $E_0$-classes whose restrictions form a basis
modulo the entire subspace $O+hO$, and take their global full-group
lifts. Their finite classes modulo $UH$ are independent jointly
with the chosen $N$ classes from $W$: any relation would restrict to a
relation among their reflection classes, and reduction modulo
$O+hO$ first kills all coefficients from $E_0$. The projection
\[
 H/UH\longrightarrow F_0/UF_0,\qquad F_0=H/\mathcal T,
\]
has kernel $\mathcal T/U\mathcal T$ of dimension $N+J_0$.
Its intersection with the selected $E_0$-span has dimension at
most $J_0$, because the independent $W$-span already occupies
$N$ kernel dimensions. Consequently the $E_0$-projections have
rank at least $V_e-J_0$.
If $V_e\leq J_0$, the lower bound eventually required for $V(L)$
is automatic. In all cases the following argument applies to a
maximal independent family of projections, whose size is at least
$\max\{0,V_e-J_0\}$.

Choose a maximal homogeneous independent family among these $E_0$
projections. The corresponding global lifts can be multiplied by
powers of $U$ so that their finite images lie in two common large
degrees $D_0,D_0-1$, according to parity. Independence persists in
$F_0$. Indeed any homogeneous relation among these powers, after
cancelling the smallest power of $U$ and reducing modulo $U$,
would contradict the chosen independence in $F_0/UF_0$.
Thus the sum of the finite ambient global-image ranks in these
two degrees is at least $V_e-J_0$.

\emph{Comparison with the half-length cut.}
We compare these images with the fixed locus. The subspace
$X_{2L}^{\tau}$ is $G$-invariant, since $\tau$ is central in $G$.
Under repetition it is $X_L$, with doubled rotation parameter
and the original reflection parameter. In high degree the map
\begin{equation}\label{eq:fp-full-fixed-injection}
 H_G^q(X_{2L},M)\longrightarrow H_G^q(X_{2L}^{\tau},M)
\end{equation}
is injective: its composite with fixed $G$-to-$K$ restriction
factors as the ambient $G$-to-$K$ injection followed by the
actual high-degree $K$ restriction isomorphism of
Lemma~\ref{lem:fp-localization}. Therefore the ambient global
images just found retain their dimensions on the fixed locus.

For $H'=G,K$ define the fixed global-image spaces
\[
 I_{H'}^q=\operatorname{im}\left[
 H_{H'}^q(\Lambda^\tau,M)
       \longrightarrow H_{H'}^q(X_{2L}^{\tau},M)\right].
\]
By naturality they contain those restricted ambient global
images. Apply \eqref{eq:fp-weight-two-gysin} to the fixed global
pair and to the fixed finite pair. The global push map is
surjective. Hence the finite push, restricted to images, gives
a surjection
\[
 \pi:I_K^{D_0}\longrightarrow I_G^{D_0-1}.
\]
The finite $G$-to-$K$ injection embeds $I_G^{D_0}$ into its
kernel, so
\begin{equation}\label{eq:fp-fixed-image-inequality}
 \dim I_G^{D_0}+\dim I_G^{D_0-1}\leq\dim I_K^{D_0}.
\end{equation}
To see the surjectivity on images directly, represent a class of
$I_G^{D_0-1}$ by a global fixed $G$-class, lift that class by the
surjective global push map, and restrict the lift to the finite
fixed pair. Naturality gives its prescribed image under the finite
push map. This uses the exact sequence before taking images; it
does not assert that taking images preserves exactness in general.
On fixed $K$-spaces, global restriction is the coefficientwise
map $B_\infty[v]\to B_L[v]$. Since $B_L$ has finite total
dimension, for sufficiently large $D_0$ the last dimension in
\eqref{eq:fp-fixed-image-inequality} is exactly $V(L)$. Choose
$D_0$ large enough for all previous comparisons. This choice is made
for the fixed cut $L$, after choosing the finite family of lifts.
No high-degree stabilization threshold uniform in $L$ is required.
We obtain
\begin{equation}\label{eq:fp-even-image-bound}
 V(L)\geq V_e-J_0.
\end{equation}

\emph{The inequality and its equality case.}
Subtract the global-image dimensions from
\eqref{eq:fp-rank-identity}, and use
\eqref{eq:fp-visible-pairs}--\eqref{eq:fp-even-image-bound}:
\begin{align*}
 T(2L)-T(L)
 &=2(N+J_0)-V_o-V_e+V(L)\\
 &\geq J_0+2N-V_o\geq0.
\end{align*}
If equality holds, necessarily $V_o=2N$. But
\[
 V_o=\dim(O+hO)=N+\dim(hO)-\dim(O\cap hO),\qquad
 \dim(hO)\leq N.
\]
Thus $\dim(hO)=N$ and $O\cap hO=0$. The first equality says
precisely that $h:O\to hO$ is injective. This proves both the
doubling inequality and its equality condition, without assigning
indicators to potentially overlapping degree towers.
\end{proof}

\begin{proof}[Proof of Theorem~\ref{thm:fp-obstruction}]
All estimates are now in place. Since the transient count is a bounded
integer, choose an interval where it is maximal. Doubling forces
equality throughout its dilated regular intervals, and hence
injectivity of $h$ there. The bounded-gap sequence supplies a global
class which first becomes visible inside one of these intervals.
Local vanishing of $h$ at that crossing yields the contradiction.

\emph{Intervals on which equality holds under doubling.}
Suppose there are only finitely many prime images. By
Proposition~\ref{prop:fp-transients}, $T$ is a bounded nonnegative
integer on the regular positive cuts. It is constant on each
interval without critical values. Its maximum is attained on a
nonempty bounded open regular interval $I\subset(0,\infty)$.
If $x$ is regular, every iterated half $x/2^k$ is regular: a
critical loop at a half cut would iterate to one at $x$.
For every regular $x\in2^kI$, iterating
\eqref{eq:fp-doubling} gives
\[
 T(x)\geq T(x/2)\geq\cdots\geq T(x/2^k)=\max T.
\]
All these inequalities are equalities. In particular, for $k\geq1$
the equality condition of Proposition~\ref{prop:fp-doubling} holds
at every regular point of $2^kI$.

\emph{A class visible at the upper cut and invisible at the lower cut.}
Choose regular points $x_1<x_2$ in separated lower and upper
quarters of $2^kI$. Their distance tends to infinity with $k$.
The interval
\[
 (\alpha x_1+c,\ \alpha x_2-c)
\]
of possible degrees has length tending to infinity. Remove the
fixed-width bands about each $\beta_i x_1$ and $\beta_i x_2$,
and the fixed low-degree range used in Lemma~\ref{lem:fp-gaps}.
These exclusions have a total length bounded independently of
$k$. A sequence of integers with gaps at most $b$ supplies
arbitrarily many $q_s$ in the growing interval, whereas these
bands exclude only a bounded number. For large $k$, choose a
remaining $q_s$ and set $Z=Z_s\in W$.
For this choice, $q_s$ lies above the ordinary spectral line at
$x_1$ and below it at $x_2$, and it avoids the local index bands
at both cuts. These are exactly the two alternatives distinguished
by Lemma~\ref{lem:fp-gaps}.

At $x_1$ its degree is either in an absent gap or above the last
band, so $Z$ restricts to zero by Lemma~\ref{lem:fp-gaps}. At
$x_2$ it is in a present gap, where global reflection restriction
is injective, so its restriction is nonzero. Visibility of a
fixed global class is monotone with increasing cuts: restriction
to a smaller sublevel factors through restriction to a larger
one. There are finitely many events between $x_1$ and $x_2$.
Choose the first event at which $Z$ becomes nonzero, and choose
tight regular cuts $x_-<x_+$ around it, both inside $2^kI$.

\emph{The first crossing contradicts the equality case.}
The post-event class $Z|_{X_{x_+}}$ is nonzero and has zero
restriction to $X_{x_-}$. Exactness therefore places it in the
image of the actual relative crossing map
\[
 H_R^{q_s}(X_{x_+},X_{x_-})\longrightarrow B_{x_+}^{q_s}.
\]
By Lemma~\ref{lem:fp-local}, multiplication by $h$ is zero on
that relative crossing, including all simultaneous labels.
Naturality of the ambient cup action now gives
\[
 h\,Z|_{X_{x_+}}=0.
\]
On the other hand $Z|_{X_{x_+}}$ belongs to
$\operatorname{im}(W\to B_{x_+})$, and multiplication by $h$
is injective there because $x_+$ lies in the maximal dilated
interval. This contradicts the nonvanishing of $Z|_{X_{x_+}}$.
The finite-prime assumption is impossible.
\end{proof}

\begin{remark}\label{rem:fp-scope}
The proof uses neither nondegeneracy nor positivity of every mean
index. Its local estimates use only finitely many root kernels,
and its low-degree argument explicitly includes all zero-mean-index
iterates. The local comparison is supplied by
Proposition~\ref{prop:local-thom-realization}. The global inputs
needed from a particular manifold are
the all-class spectral estimate, uniform additive death,
the bounded $w$-annihilator, the actual paired decomposition and
full-group lifts, and the bounded-gap sequence in $W$ stated in
Theorem~\ref{thm:fp-obstruction}.
\end{remark}


\section{Proof of the main theorem}\label{sec:main-proof}

We finish by checking the hypotheses of the finite-prime criterion
for each manifold in Theorem~\ref{thm:main}. The filtration hypotheses
come from the estimates for the given metric in
Sections~\ref{sec:based-analytic} and~\ref{sec:finite-cover}.
The remaining hypotheses concern global classes and restriction maps;
these were computed using symmetric metrics in
Sections~\ref{sec:sphere-topology} and~\ref{sec:projective-topology}.
Since those global data do not depend on the metric, they apply to
the metric under consideration.

\begin{proof}[Proof of Theorem~\ref{thm:main}]
Fix a smooth Riemannian metric $g$ on one of the stated manifolds $M$.
For $M=S^2=\CP^1$, the conclusion also follows directly from
Bangert~\cite{Bangert1993} together with Franks~\cite{Franks1992}.
Thus the classical surface case does not depend on any new chain
construction in this manuscript. The argument below provides the
uniform route, with the low-dimensional projective calculation retained
as a consistency check.
The controlled cuts exist by
Proposition~\ref{hyp:controlled-cut-realization}, and
Proposition~\ref{prop:based-uniform-persistence} gives uniform death and
spectral resonance on the based loop space, with a positive slope.
Theorem~\ref{thm:promotion} transfers both estimates to the free loop
space. By Corollary~\ref{cor:relative-analytic}, they also hold relative
to the constant loops. The degreewise global homology groups have a
uniform finite bound. In particular, the first two hypotheses of
Theorem~\ref{thm:fp-obstruction} hold.

For $M=S^n$ with $n\ge3$, apply Proposition~\ref{prop:sphere},
and take $h=e$ and $E_0=E$ in the obstruction theorem.
Its odd-stage summand and its multiple by the evaluation class give the
two summands denoted by $W$ and $hW$ in
Theorem~\ref{thm:fp-obstruction}; the even-stage summand gives the third.
The evaluation class has degree $n\ge2$. The classes of $W$ have
$U$-torsion lifts, the even classes have full-group lifts, and a fixed
power of $w$ annihilates the relative reflection module. The odd bottom
degrees form an arithmetic progression with gap $4(n-1)$ and hence
have bounded gaps.
Thus all the remaining hypotheses of that theorem are satisfied.

For complex and quaternionic projective spaces and for the Cayley plane,
the corresponding statements are supplied by
Theorem~\ref{thm:projective-equivariant-topology}, with $E_0=E$.
Here the evaluation
class has degree $2$, $4$, or $8$, respectively. The actual global
relations and lifting statements proved there give the required
annihilation and restriction properties. Again, the odd bottom degrees
form an arithmetic progression with gap $2\nu$, where $\nu$ is the
index increment defined in that section. The case $S^2$ is included by
$S^2=\CP^1$; $S^4=\HP^1$ is consistent with either calculation.

If $g$ had finitely many prime images,
Proposition~\ref{prop:local-thom-realization} and
Lemma~\ref{lem:fp-local} would give the uniform local bounds used in
Theorem~\ref{thm:fp-obstruction}. The obstruction theorem therefore excludes the finite-prime
possibility and proves the assertion.
\end{proof}

\subsection*{Acknowledgements and computational assistance}
OpenAI Codex was used  in developing the draft arguments, some of the proof strategies and calculations, 
checking selected algebraic consequences. We take  responsibility for the human verification and exposition.

Earlier drafts consulted unpublished work by \v{C}adek--Moravec,
Clivio, and Cieliebak--Oancea--Shelukhin. The finite-cut and doubling
strategy was developed from the September 2026 OpenAI sphere
manuscript~\cite{OpenAI2026}; the bibliography identifies the archived
revision used here.
These manuscripts are acknowledged as sources of ideas and strategy.
Their conclusions are not invoked as external theorems in the proofs
above. The claim of logical independence concerns the stated proofs
and their published or accepted inputs; it does not imply historical
independence or originality of the strategy.


\end{document}